\documentclass[3p,twocolumn,12pt]{elsarticle}

\usepackage{amssymb}
\usepackage{amsmath}
\usepackage{xcolor}
\usepackage{graphicx}
\usepackage{algorithm}
\usepackage{algorithmic}
\usepackage{bm}
\usepackage{siunitx}

\journal{Nuclear Physics B}

\begin{document}

\begin{frontmatter}



\title{Anatomical Connectivity from Tractography-free Diffusion Tensor Imaging and Its Application to EEG Brain Imaging}


\author{Joonas Lahtinen\\ \today} 

\affiliation[1]{
    organization={Tampere University},
    addressline={Korkeakoulunkatu 1}, 
    city={Tampere},
    postcode={33720}, 
    state={Pirkanmaa},
    country={Finland}
}

\begin{abstract}
This proof-of-concept paper demonstrates how the Kalman filter, as a brain imaging method, can be enhanced by incorporating raw diffusion tensor imaging (DTI) data without tractography, in addition to EEG recordings. An efficient algorithm is designed to streamline the otherwise tedious process of inferring connectivity between brain regions from DTI data. The DTI-based evolution model is applied to the Kalman filter and the Standardized Kalman filter, and these are compared with their conventional counterparts using the random walk evolution model. The numerical experiments use synthetic somatosensory and auditory evoked potentials. The results show that the DTI-modelled Kalman filter can obtain deep activity in the brainstem and thalamus. Moreover, the connection model limits the estimated spreads to the correct or nearest anatomical brain region.

\end{abstract}




\end{frontmatter}



\section{Introduction}
The diffusion process is one of the most widely used stochastic models that is in finance \cite{MONSALVECOBIS20113073}, computational neuroscience \cite{MARINOV2014169}, cellular transport \cite{FENZ2012260}, and social networks \cite{BakenovaKamila2025IDMi}, to name a few. An advanced and modern magnetic resonance neuroimaging (MRI) modality, Diffusion tensor imaging (DTI), has become increasingly utilized in clinical routines, especially in brain trauma research. 

DTI  is particularly valuable for visualizing tissues with organized, fibrous structures, such as the white matter tracts of the brain or the muscle fibers of the heart. In these tissues, water molecules tend to diffuse more freely along the length of the fibers (axial diffusion) and are more restricted when moving perpendicular to the fiber orientation (radial diffusion) \cite{DTIhandbook}.

At the molecular level, water undergoes random motion—known as Brownian motion—in all spatial directions. When this diffusion occurs equally in every direction, it is isotropic, as in water in an unrestricted environment. Conversely, in tissues with directional structure, diffusion becomes anisotropic, meaning water moves preferentially along the underlying architecture, much like water flowing through a bundle of straws. The degree of anisotropy reflects how directionally constrained the diffusion is. Because water diffusion patterns vary with tissue composition and organization, DTI enables the noninvasive assessment of these properties—most notably, the mapping of water movement along neuronal axons --- by capturing and quantifying their anisotropic behavior.

Since water preferentially diffuses along the length of nerve fibers, DTI leverages this property to map the complex network of neural connections \cite{Skudlarski2008DTIbrainconnectivity,Greicius2009RDTIconnectivity,KowalczykNatalia2018DTIbrainconnectivity}. Unfortunately, DTI images are inherently blurry and noisy \cite{TournierJacques-Donald2011DTIandBeyond,GuptaVikash2013DTIblur,Chan2014HighResDTI}, which poses challenges for its usage and reliability. However, if the uncertain measurements are treated probabilistically, we can draw conclusions appropriate to the certainty the imaging provides.

Kalman filtering \cite{GalkaAndreas2004KalmanEEG} can reliably track human brain activity; as a special type of Bayesian filter, it yields estimates analytically without further numerical methods. The Kalman filter provides a recursive filtering algorithm that produces a spatiotemporal estimate of the brain state at a given time point using measurements from the initial moment up to the current time step. It allows us to use the high temporal resolution of EEG to draw inferences about human brain activity and its evolution. An advanced version of the Kalman filter dedicated to EEG brain imaging is called the Standardized Kalman Filter (SKF) \cite{Lahtinen2024SKF}. It seeks brain states from the Z-score field, testing the hypothesis that brain activity is localized to a specific area. This allows us to detect brain activity even far from the sensors, which is not possible with a basic Kalman filter due to its strong bias toward superficial sources of brain activity. However, applying the method to different source scenarios and brain and electrode models has proven challenging; thus, many variations to its state formulation and parameter tuning have been proposed \cite{Dilshanie2026SKF,Lahtinen2026CRSKF,piispa2025DSKF}. The main difficulty arises from the so-called evolution model, a stochastic forward equation used as a prior model of brain activity dynamics. In the state-of-the-art, the evolution is modeled as a simple random walk, with the variances of certain brain regions being tuned. While this has turned out to be a working approach in some situations, it is not applicable in others, where signals undetectable by scalp EEG travel across connected brain areas, causing separate activity to appear at a distance. This kind of spatial "jump" is nearly impossible to model without prior knowledge of possible connections between brain regions.

In this paper, we use DTI-derived information to establish probabilistic connections between anatomical brain regions. An efficient algorithm is developed to determine the flow proportions for different brain regions, along with the diffusion process as determined by DTI. The flow information is then turned into a stochastic evolution model that is incorporated into a Kalman filter \cite{GalkaAndreas2004KalmanEEG} and Standardized Kalman filter \cite{Lahtinen2024SKF}. The methodology is demonstrated with a numerical example considering a causal activity component system appearing in early somatosensory evoked potentials, namely components from P14 to P22/N22 \cite{Deiber1986P22,Morioka1991N16,DesmedtOzaki1991P20N20,NoelPierre1996OoNa,Valeriani2014P16}, and simulated auditory evoked potentials according to \cite{SchergMichael1989AEP,Picton2011AEP,PaavoRonniMstr2022}. The paper aims to answer two main questions: (1) can DTI be used without pre-processing and tractography, and (2) could a connectivity-based transformation model be beneficial for EEG source imaging via Kalman filtering. Based on the results, the answer to (1) is positive; however, it is not time-feasible with the pipeline used, and the answer to (2) is also positive. 

\section{Methods}
\subsection{Diffusion Tensor Imaging}
The diffusion tensor $\sigma$ of the DTI is computed by fitting the diffusion-weighted image to six directions using the Stejskal-Tanner equation:
\begin{equation}
    I_k=I_0\exp\left(-b\hat{\bf v}_k^T\sigma\hat{\bf v}_k\right),
\end{equation}
where $I_k$ is the signal intensity measured in the direction given by the unit vector $\hat{\bf v}_k$, $I_0$ is the baseline intensity, $b$ is the diffusion weighting factor. By forming a matrix out of all the directions as
\begin{equation}
    V=\begin{bmatrix}
        \hat{\bf v}_1 & \cdots &\hat{\bf v}_6
    \end{bmatrix}
\end{equation}
and by denoting
\begin{equation}
    {\bf f}=\begin{bmatrix}
        \log(I_0/I_1)\\\vdots \\ \log(I_0/I_6),
    \end{bmatrix}
\end{equation}
we can formulate a linear system
\begin{equation}
    {\bf f} = b\cdot\mathrm{diag}\left(V^T\sigma V\right)=b(V\odot V)^T\mathrm{vec}(\sigma),
\end{equation}
where $\odot$ denotes the Khatri-Rao product \cite{khatri1968product}, and $\mathrm{vec}$ is the vectorization operator turning a matrix into a vector as follows:
\begin{equation}
    \mathrm{vec}(\sigma)=\begin{bmatrix}
        \sigma_{11} & \cdots & \sigma_{13} & \sigma_{21} & \cdots &\sigma_{33}
    \end{bmatrix}^T.
\end{equation}
This equation can be solved for the diffusion tensor $\sigma$ in practice by the least squares method.

\subsection{Kalman filter}

The observation model in spatiotemporal electric source imaging can be expressed as
\begin{equation}\label{eq:frwrdmodel}
    {\bf y}_t=L{\bf s}_t+{\bf r}_t,\quad \textnormal{for }t=1,\dots, T
\end{equation}
where ${\bf y}_t\in\mathbb{R}^m$ denotes EEG observations measured by scalp electrodes, $L\in\mathbb{R}^{m\times n}$ is the lead field or system matrix, ${\bf s}_t\in \mathbb{R}^n$ denotes the unknown brain activity expressed in vector form, and ${\bf n}_t$ is observation noise that is assumed to follow a Gaussian $\mathcal{N}({\bf 0},R_t)$ in this filtering scheme.

In a Bayesian sense, the Kalman filter provides a posterior for each state ${\bf x}_t$ given all the observations up to step $t$,  which is denoted ${\bf y}_{1:t}$. The Kalman filter works recursively by computing the predictive "prior" distribution $p({\bf s}_t\mid {\bf y}_{1:t-1})$ via the Chapman-Kolmogorov equation:
\begin{equation}
\begin{split}
    p({\bf s}_t\mid {\bf y}_{1:t-1}) &= \int_{\mathbb{R}^n} p({\bf s}_t\mid {\bf s}_{t-1})\\
    &\times p({\bf s}_{t-1}\mid {\bf y}_{1:t-1})\, \mathrm{d}{\bf s}_{t-1},
\end{split}
\end{equation}
where $p({\bf s}_t\mid {\bf s}_{t-1})$ is the transition probability, and the recursivity comes from the posterior distribution of the previous state $p({\bf s}_{t-1}\mid {\bf y}_{1:t-1})$. The current posterior state is then again obtained using Bayes' rule: $p({\bf s}_t\mid {\bf y}_{1:t})\propto p({\bf y}_{1:t}\mid {\bf s}_t)p({\bf s}_t\mid {\bf y}_{1:t-1})$, where the likelihood model $p({\bf y}_{1:t}\mid {\bf s}_t)$ comes from Eq. (\ref{eq:frwrdmodel}). 

Kalman filter assumes that the states follow a linear dynamic model
\begin{equation}
    {\bf s}_t=A_t{\bf s}_{t-1}+{\bf q}_t,
\end{equation}
where $A_t$ is called the transition matrix and the process is affected by process noise ${\bf q}_t$ following a Gaussian $\mathcal{N}({\bf 0},Q_t)$.

The trackability of the Kalman filter is one of its advantages. This allows us to compute the posterior statistics (mean $\hat{\bf s}_{t\mid t}$ and covariance $P_{0\mid 0}$) recursively \cite{SarkkaSimo2013}:
\begin{algorithm}[H]
\caption{Kalman filtering algorithm}
\begin{algorithmic}[1]
\STATE \textbf{Input:} Prior statistics: $\hat{s}_0$ and $P_{0\mid 0}$
\FOR{$t = 1$ to $T$}
    \STATE $\hat{\bf s}_{t\mid t-1}=A_t\hat{\bf s}_{t-1\mid t-1}$
    \STATE $P_{t\mid t-1}=A_tP_{t-1\mid t-1}A_t^T+Q_t$
    \STATE $S_t=LP_{t\mid t-1}L^T+R_t$
    \STATE $K_t = P_{t\mid t-1}L^TRS_t^{-1}$
    \STATE $\hat{\bf s}_{t\mid t}=\hat{\bf s}_{t\mid t-1}+K_t({\bf y}_t-L\hat{\bf s}_{t\mid t-1})$
    \STATE $P_{t\mid t}=P_{t\mid t-1}-K_tS_tK_t$
   \ENDFOR
\end{algorithmic}\label{Algo:KF}
\end{algorithm}

In this paper, we have a static transition matrix $A$ that is either the identity matrix for the conventional random walk evolution model or is derived using DTI as described later.

\subsection{Standardized Kalman filter}
Standardization is a technique developed in \cite{PascualMarqui2002} and later formalized as a random-field approach based on Z-scores \cite{Lahtinen2024Onbias}, since these scores have an order-preserving transformation to goodness-of-fit measures that indicate how well a feasible solution matches the data \cite{HoltershinkenOnEquivalence}. Because the Kalman filter estimate is the Gaussian mean, the transformation to time-variant Z-scores via 
\begin{equation}
    W_t=\mathrm{Diag}\left(P_{t\mid t-1}^{1/2}L^TS_t^{-1}LP_{t\mid t-1}^{1/2}\right)^{-1/2}P_{t\mid t-1}^{-1/2}.
\end{equation}
By requiring that each standardized Kalman estimate can be written as $\hat{\bf z}_t=W_t\hat{\bf s}_t$ for every step $t$ up to $T$, the evolution must be of the form
\begin{equation}
    \hat{\bf z}_t=W_tAW_{t-1}^{-1}\hat{\bf z}_{t-1}+W_t{\bf q}_t,
\end{equation}
where transition matrix $A$ is assumed static. For the Kalman filter algorithm (Algorithm \ref{Algo:KF}), this means that the transition matrix for Standardized Kalman filter is necessarily time-variant
\begin{equation}
    \tilde{A}_t = W_tAW_{t-1}^{-1}.
\end{equation}
Also, the process noise covariance is now $W_tQ_tW_t^T$.

\subsection{DTI diffusion process}
While the connectivity matrix for different brain regions can be formed using DTI fiber tracking \cite{conturo1999fiber,mori2002fiber,yap2010fiber}, the crossing points of the white matter fibers can be resolved by deconvolution-based algorithms or Q-Ball imaging \cite{Tournier2008DTIcrossingfibers}, and the connected fiber ends can be sought in the anatomical brain regions of interest. Needless to say, this process as a whole is time-consuming, which is undesirable, as it needs to be done for each subject separately. Moreover, a couple of pathologies emerge through the workflow: (1) Fiber tracking is usually done with the Euler method \cite{ChenBin2007DTIEulerMethod}, and due to the noise in the data, the end result varies highly depending on the chosen step size. Even a smaller step size may not yield better results, as it can cause, e.g., the track to get stuck or enter an endless loop. (2) As the step direction is determined by the eigenvector corresponding to the largest eigenvalue, it filters out many valid possible tracks to the target as an independent point estimate done each step. For example, in a case where two or three eigenvalues are closely equal, picking one of these principal directions is highly arbitrary.

We note that the diffusion tensor $\sigma\colon \mathbb{L}^3\to\mathrm{Sym}(3)$ is a discrete map that provides a symmetric $3\time3$-matrix for each node of the regular lattice $\mathbb{L}^3$ on $\mathbb{R}^3$. Here, we model the DTI diffusion process as
\begin{equation}\label{eq:StochasticEq}
    \mathrm{d}\tilde{\bf x}_{t}=\sigma(\tilde{\bf x}_t)\mathrm{d}\tilde{\bf W}_t
\end{equation}
without a drift, where $\tilde{\bf W}_t\in\mathbb{R}^3$ is a Wiener process. Interpreting this from the point of view of Riemann manifolds, we can call the equation an index-raising operation. In general, this means that $\sigma(\tilde{\bf x}_t)$ can be thought as a (1,1)-tensor; we can make a geometrical connection by stating that the spatially dependent diffusion tensor induces a Riemannian metric tensor $g_\sigma=\sigma^{\dagger}$, as proved in \ref{app:FP-DTg}. It should be noted that a proper metric is non-degenerate, which cannot be guaranteed in numerical settings.

M. Capitaine also noted a similar connection in 2000 \cite{Capitaine2000} for the elliptic diffusion process. Shifting the framework from a stochastic process to a Riemannian manifold lets us use new tools to derive the connections. First of all, we can observe that (\ref{eq:StochasticEq}), as a relatively complex process in a flat space, is just a simple random walk in the manifold $(\mathcal{M},g_\sigma)$. Using the Fokker-Planck equation in a Riemannian context, we can derive the probability distribution for a signal traveling along with diffusion to point ${\bf x}$, when the process starts at ${\bf x}_0$: 
\begin{equation}
    p(\mathbf{x},t\mid \mathbf{x}_0)=\frac{1}{Z\sqrt{t^3}}\exp\left(-\frac{d_\sigma(\mathbf{x},\mathbf{x}_0)^2}{2t}\right),
\end{equation}
where $d_\sigma(\mathbf{x},\mathbf{x}_0)$ is the geodesic distance from point ${\bf x}_0$ to ${\bf x}$, and $Z$ is the normalization constant at $t=1$. This probability is then used to calculate the hitting probabilities for surfaces of the initial region that are facing those of the target region. This way, we are avoiding the use of point estimates that determine the most likely paths for diffusion (tractography) and use the DTI-provided information fully and as-is. Here we utilize the practical notion that the curvature from the thalamus to the neocortex, or between neocortical regions, is particularly low, i.e., the geodesic path curves only slightly. This allows us to compute the hitting probability as the probability for the volumetric-like region $\Omega = \lbrace \left(\mathbf{x},[T(\mathbf{x},\mathbf{x}_0),\infty]\right)\colon \mathbf{x}\in\mathcal{S}\rbrace$, where $\mathcal{S}$ is the surface of the target region and $T(\mathbf{x},\mathbf{x}_0)$ denotes the initial hitting time. The initial hitting time is defined as the first physically plausible moment a signal from the initial point to the target has moved through the whole geodesic path $\bm{\gamma}_\sigma({\bf x},{\bf x}_0)\colon \mathbb{R}^3\times \mathbb{R}^3\to\mathbb{R}^3$. The initial hitting time is defined as follows:
\begin{equation}\label{eq:Td-connection}
    \frac{1}{T(\mathbf{x},\mathbf{x}_0)}=\frac{1}{d_\sigma(\mathbf{x},\mathbf{x}_0)^2}\int_{\bm{\gamma}_\sigma(\mathbf{x},\mathbf{x}_0)}\hat{\bf n}^T\sigma \hat{\bf n}\,\mathrm{d}s,
\end{equation}
where the path integral, scaled by the path length, can be interpreted as the average velocity. Due to the initial hitting time, our formulations take into account the radial diffusion of matter. The radial diffusion serves as an effective model for the branching of the fiber bundle. This provides the equation (derived in \ref{app:probFormula}) for the brain activity proportion arriving at the target:
\begin{equation}
\begin{split}
    P_\infty-P_T&\propto \left(1-\exp\left(-\frac{h^2}{8T(\mathbf{x},\mathbf{x}_0)}\right)\right)\\
    &\times\mathrm{erf}\left(\frac{d_\sigma(\mathbf{x},\mathbf{x}_0)}{\sqrt{2T(\mathbf{x},\mathbf{x}_0)}}\right),
\end{split}
\end{equation}
where $h$ is the diameter of a voxel. The formulation of the proportion assumes that we examine a streamline that occupies a single voxel at each time point, since the diffusion tensor is defined discretely for each voxel. 

To calculate the proportions, we need to know the geodesic paths. Multiple geodesic pathfinding algorithms exist, such as the Fast Marching Method \cite{Sethian1996FMM,Mirebeau2018FMM}, the A* search algorithm \cite{HartPeter1968Astar}, and a geodesic path can be obtained by solving the heat equation in the domain \cite{Crane2013GeodesicInHeat}. However, considering the given formatting of the diffusion tensor and its size $256^3$, all the aforementioned methods are computationally too taxing, and since we will stick with a classical shooting method category \cite{Noakes1998GeodesicShooting,Bryner2017GeodesicShooting}. Namely, to make the algorithm as light as possible, we will split the problem into subproblems considering two previously optimized ends of a Euclidean path and aim to find a middle point that minimizes the time spent on the path consisting of that middle point; i.e., we solve the following minimization problem 
\begin{equation}\label{eq:geodesicMin}
\scalebox{1.1}{$\min\limits_{{\bf x}_m\in P}\left\lbrace \frac{\left\|\mathbf{x}_b-\mathbf{x}_m\right\|_2+\left\|\mathbf{x}_a-\mathbf{x}_m\right\|_2}{\left\|\sigma({\bf x}_m)(\mathbf{x}_b-\mathbf{x}_m)\right\|_2+\left\|\sigma({\bf x}_a)(\mathbf{x}_a-\mathbf{x}_m)\right\|_2}\right\rbrace.$}
\end{equation}

As the problem of finding the optimal path consisting of $N$ line segments should be solved accurately as a system of $N-1$ variables, for the shooting algorithm that divides the problem into subproblems consisting of two linear paths, we check the path for outliers, i.e., suspiciously slow sections, and optimize those separately. The error of this fast algorithm is upper-bounded by $Cd_\sigma(\mathbf{x},\mathbf{x}_0)^2(1/2)^{\sqrt{n}}$ after $n$ steps and for some constant $C$, as proved in \ref{app:shootingConverg}. 

The full path-finding algorithm reads as follows:

\begin{algorithm}[H]
\caption{Shooting algorithm for geodesics}
\begin{algorithmic}[1]
\STATE \textbf{Input:} End points of the path ${\bf x}_{a}$ and ${\bf x}_b$, diffusion tensor $\sigma$, number of steps $N\in\mathbb{N}$, correction frequency $K\in\mathbb{N}$.
\FOR{$n = 1$ to $N$}
    \STATE Compute the mid point ${\bf x}_m =({\bf x}_a+{\bf x}_b)/2$
    \STATE Define the plane $P\perp ({\bf x}_b-{\bf x}_a)$ that goes through ${\bf x}_m$.
    \STATE Find optimal ${\bf x}_m^*$ by solving (\ref{eq:geodesicMin}) and add it to the point set $\mathcal{G}$.
    \STATE Find subsequent points ${\bf x}_a,{\bf x}_b\in \mathcal{G}$ with longest distance Euclidean from each other.
    \IF{$K$ divides $n$}
    \STATE Find subsequent points ${\bf x}_a,{\bf x}_m,{\bf x}_b\in \mathcal{G}$ that have higher value for object function of Eq. (\ref{eq:geodesicMin}) than 95\% of similar segments, and recalculate ${\bf x}_m$.
    \ENDIF
   \ENDFOR
\end{algorithmic}\label{Algo:shootingMethod}
\end{algorithm}

In this paper, we use 12 iteration steps and outlier reduction every fourth step.

In the implementation, the connection from a brain region $\mathcal{R}_a$ to $\mathcal{R}_b$ is computed by first determining the opposing surfaces of the inspected regions; let us denote those as $\mathcal{S}_a$ and $\mathcal{S}_b$. From each voxel $i$ that touches the surface $\mathcal{S}_a$, we calculate the hitting probability $\pi_{ij}$ to voxel $j$ overlapping with $\mathcal{S}_b$. Then for any source position ${\bf x}_a\in \mathcal{R}_a$ and ${\bf x}_b\in \mathcal{R}_b$, the corresponding element of the transition matrix $A$ is computed as the sum over all $\pi_{ij}$ and scaling it in such a way that 
\begin{equation}
    \sum_{u=1}^nA_{uv}=1,\quad \textnormal{for }v=1,\dots,n,
\end{equation}
which means that there is no activity loss.

\section{Experiments}
To demonstrate the impact of DTI-based evolution modeling, we generated two distinct synthetic time-dependent brain activity scenarios reflecting real evoked potentials in somatosensory (SEP) and auditory (AEP) cases.

We model synthetic SEP data as subsequent Gaussian pulses appearing at 14, 16, 20, and 22 \unit{\milli\second}. All four components have either deep activity in the brainstem or the thalamus. 20 and 22 \unit{\milli\second} have also cortical contributor. A synthetic setup of this kind is used to explore hierarchical Bayesian methods \cite{RezaeiA2021} and is frequently used with SKF \cite{Lahtinen2024SKF,Lahtinen2026CRSKF,Dilshanie2026SKF}. 

The deep 14 \unit{\milli\second} component has been attributed to the pons, specifically the medial lemniscus pathway, based on investigations of epilepsy patients with lesions in this area \cite{NoelPierre1996OoNa}. The 16 \unit{\milli\second} component is elicited by a somatosensory volley propagating along the medial lemniscus pathway \cite{Morioka1991N16,HsiehCL1995N16,Valeriani2014P16}. At 20 \unit{\milli\second}, activation occurs tangentially directed in the somatosensory cortex \cite{DesmedtOzaki1991P20N20,allison1991cortical,buchner1995somatotopy,FuchsManfred1998SEP}, while concurrent thalamic activity has been observed in the ventral posterolateral region of the thalamus \cite{GotzTheresa2014TIPa}. The cortical and radial 22 \unit{\milli\second} component is localized variably across participants, appearing in either Brodmann area 1 or 4 \cite{buchner1995somatotopy}. Simultaneous activation of the left thalamus is also observed at this time \cite{PapadelisChristos2011BaBa}.

With SEP data, we focus on the estimated activity spread qualitatively, based on the estimation plotted on the brain model, and quantitatively by computing the normalized amplitude time series, since the standardized and non-standardized methods do not provide estimates in the same units. Key aspects are how well the sources are localized and how cortical and deep activity amplitudes relate.

The AEP experiment in this paper is motivated by \cite{PaavoRonniMstr2022}, where a total of six dipolar generators were used to describe the auditory stimulus. The activity starts in the insula of the left hemisphere. This triggers activity in the superior temporal gyrus and simultaneously travels to the right insula, which, when activated, causes activity in the right superior temporal gyrus. Superior temporal gyri will cause activity in the transverse temporal gyri of the corresponding side. The activity causality and the waveform are modeled using the Jansen-Rit neural mass model with causal couplings \cite{BabajaniFeremiAbbas2010CoupledJRmodel,SubramaniyamNarayanPuthanmadam2024NMM}. The AEP setup is particularly challenging because activity alternates between far-apart regions in different hemispheres. It is observed in practice \cite{LuckaFelix2012HBMsLORETA} and shown mathematically \cite{Lahtinen2024Onbias,HoltershinkenOnEquivalence} that sLORETA, a static standardized methodology, is not able to estimate activity of this kind correctly.
The previous numerical result with KF and SKF also shows poor localization \cite{PaavoRonniMstr2022}.

In the experiments, we conduct the same qualitative and quantitative analysis as with synthetic SEP data. In addition, we examine the time delays between the left insula -- right insula, left superior temporal -- right superior temporal, left transverse -- right transverse, and superior temporal -- transverse evolution pairs. This experiment reveals if some separate estimated sources, while being localized correctly, are too correlated, i.e., appear and vanish at the same time and the same rate. In such a case, the cross-correlation between paired time series is maximized at time point 0 \unit{\milli\second}. Because this unwanted correlation happens with the last pair, we display the three highest cross-correlation peaks for each method.

The inversion is computed using multiple variations of Kalman filtering: The basic Kalman filter (KF) \cite{GalkaAndreas2004KalmanEEG}, the Standardized Kalman filter (SKF) \cite{Lahtinen2024SKF}, and two variants utilizing the DTI model: one without standardization, DTI-KF, and one with standardization, DTI-SKF. KF and SKF use a simple random walk model introduced in \cite{Lahtinen2024SKF}.

\section{Results}

\subsection{Optimal path test}
\begin{figure}[h!]
    \centering
    \includegraphics[width=0.8\linewidth]{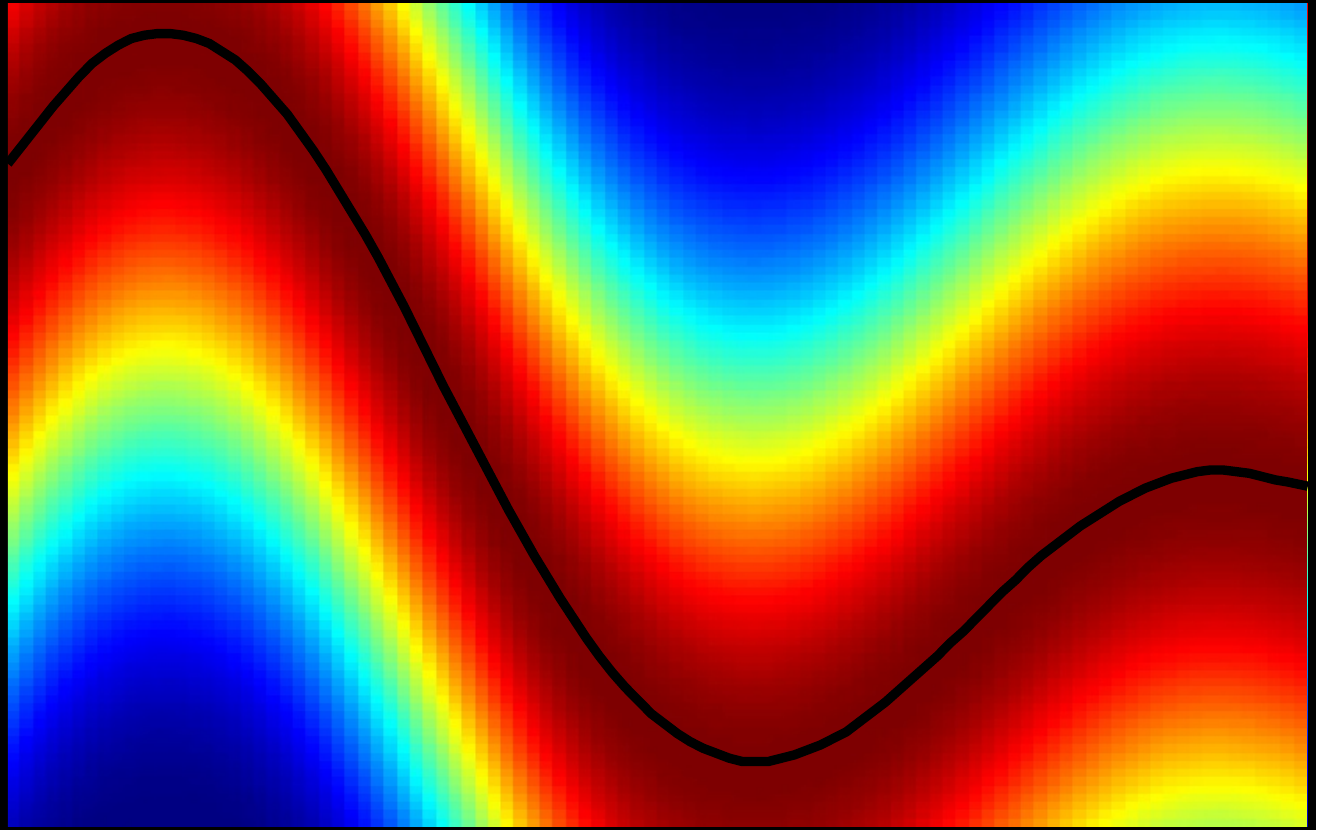}
    \caption{Isotropic diffusion field, where the black curve follows the most diffusive path. Color scale goes from dark red to dark blue, indicating the magnitude of diffusion from left to right.}
    \label{fig:PathExample}
\end{figure}

\begin{figure}
    \centering
    \includegraphics[width=0.75\linewidth]{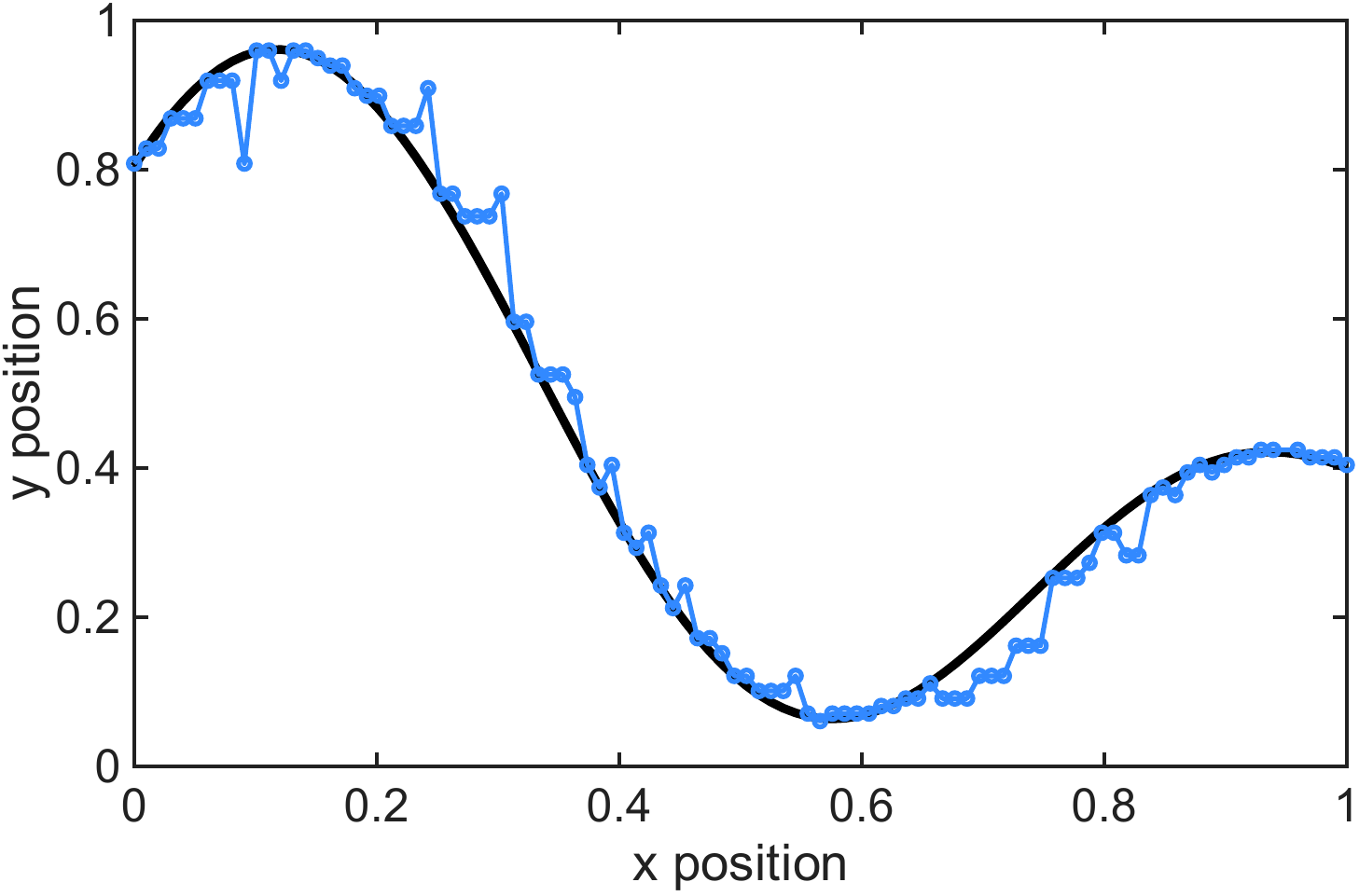}

    \includegraphics[width=0.75\linewidth]{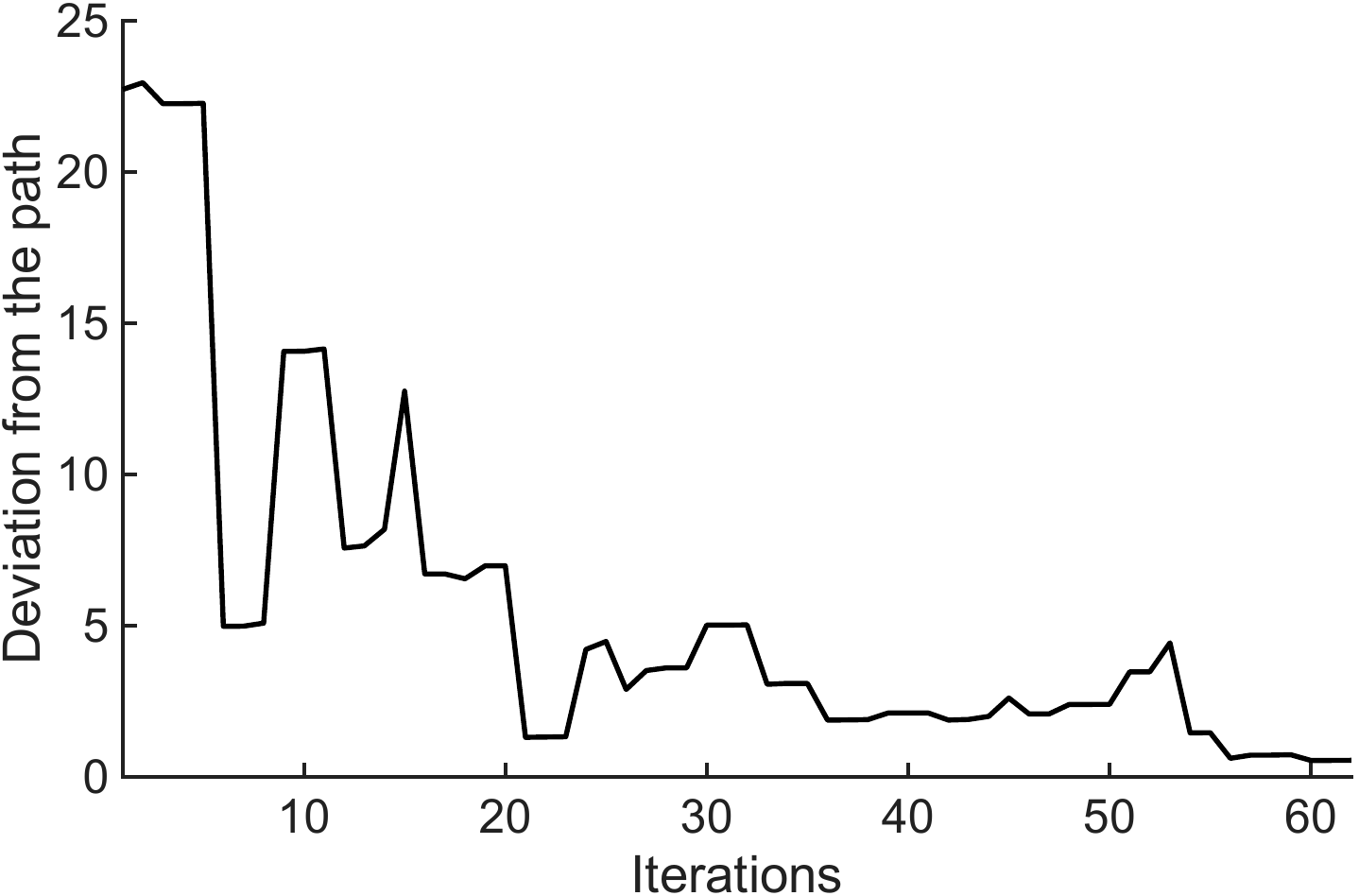}
    \caption{The upper image presents the solution after 100 iterations. The lower image presents the reduction of Euclidean path error up to 65 iterations. }
    \label{fig:PathExampleResult}
\end{figure}

As a first experiment, we demonstrate the accuracy and efficiency of the optimal path-finding algorithm used in this study. In our example, we use the following diffusion field presented in Figure \ref{fig:PathExample}. Based on the results in Figure \ref{fig:PathExampleResult}, although the path is highly non-linear, the linearly segmented tracks and the path show that the error decreases exponentially with iteration steps. 

\newpage

\subsection{Synthesized Somatosensory evoked potentials}
In our experiment with synthesized SEP, we generated the data for four different events, where 14 \unit{\milli\second} activity should appear in the pons, and 16 \unit{\milli\second} have two activity flashes: at the bottom of the brainstem and below the left thalamus, where the afferent volley pathway ends. At 20 and 22 \unit{\milli\second}, there are two highly temporally overlapping sources, one at the left thalamus and another at the somatosensory cortex. The cortical activity is modeled as tangential at 20 \unit{\milli\second} and radial at 22 \unit{\milli\second}.

\begin{figure*}
\def\imageWidth{0.13\linewidth}
\def\RotTextWidth{0.05\linewidth}
    \centering
    \hspace{0.12\linewidth}{\bf KF}\hspace{0.48\linewidth}{\bf SKF}
    \begin{minipage}[b]{\RotTextWidth}
        \rotatebox{90}{\hspace{0.5cm} 14 \unit{\milli\second}}
    \end{minipage}\begin{minipage}[b]{\imageWidth}
        \includegraphics[trim={2cm 0cm 18cm 8.5cm},clip,width=\linewidth]{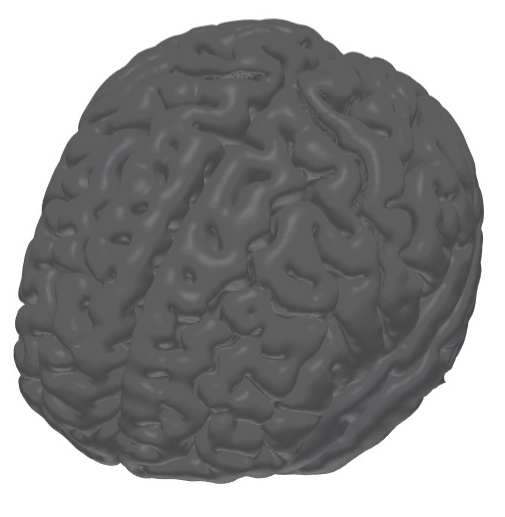}
    \end{minipage}\begin{minipage}[b]{\imageWidth}
        \includegraphics[trim={0cm 0cm 21cm 8cm},clip,width=0.8\linewidth]{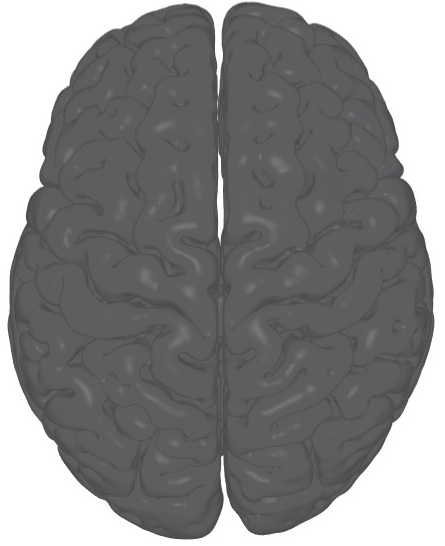}
    \end{minipage}\begin{minipage}[b]{\imageWidth}
        \includegraphics[trim={4cm 2cm 20cm 9.8cm},clip,width=0.8\linewidth]{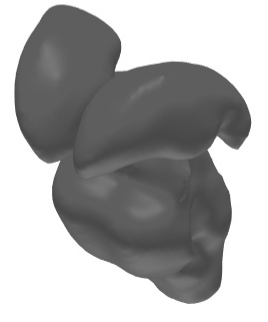}
    \end{minipage}\begin{minipage}[b]{\imageWidth}
        \includegraphics[trim={4.2cm 0.5cm 19.5cm 9cm},clip,width=0.64\linewidth]{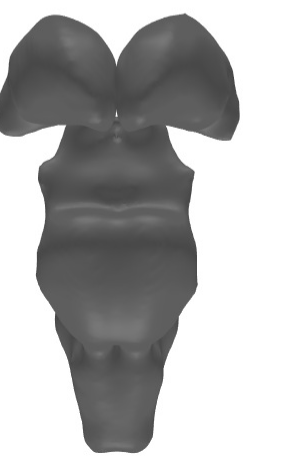}
    \end{minipage}\begin{minipage}[b]{\imageWidth}
        \includegraphics[trim={2cm 0cm 18cm 8.5cm},clip,width=\linewidth]{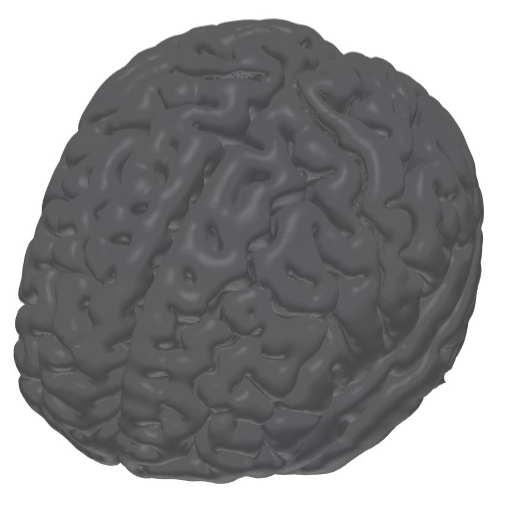}
    \end{minipage}\begin{minipage}[b]{\imageWidth}
        \includegraphics[trim={0cm 0cm 21cm 8cm},clip,width=0.8\linewidth]{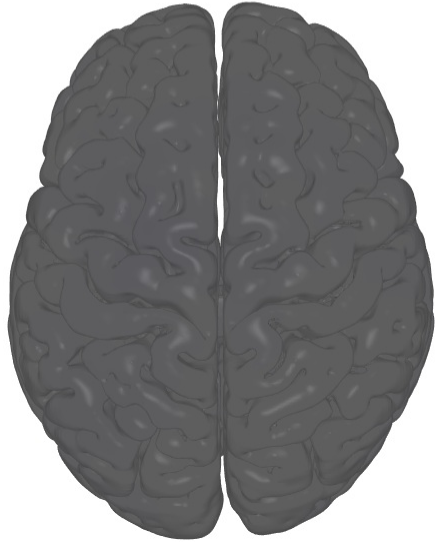}
    \end{minipage}\begin{minipage}[b]{\imageWidth}
        \includegraphics[trim={4cm 2cm 20cm 9.8cm},clip,width=0.8\linewidth]{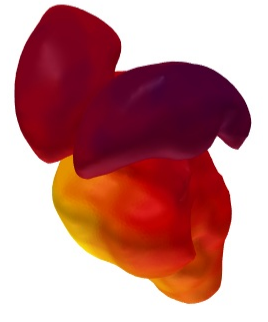}
    \end{minipage}\begin{minipage}[b]{\imageWidth}
        \includegraphics[trim={4.2cm 0.5cm 19.5cm 9cm},clip,width=0.64\linewidth]{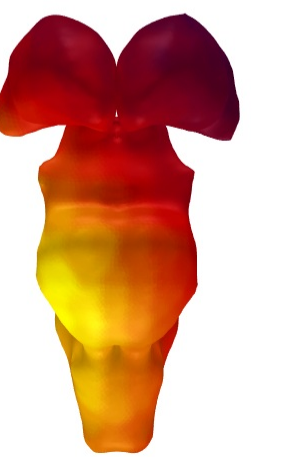}
    \end{minipage}

    \begin{minipage}[b]{\RotTextWidth}
        \rotatebox{90}{\hspace{0.5cm} 16 \unit{\milli\second}}
    \end{minipage}\begin{minipage}[b]{\imageWidth}
        \includegraphics[trim={2cm 0cm 18cm 8.5cm},clip,width=\linewidth]{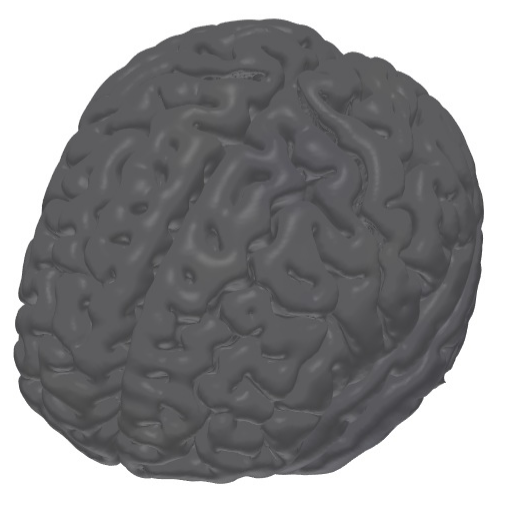}
    \end{minipage}\begin{minipage}[b]{\imageWidth}
        \includegraphics[trim={0cm 0cm 21cm 8cm},clip,width=0.8\linewidth]{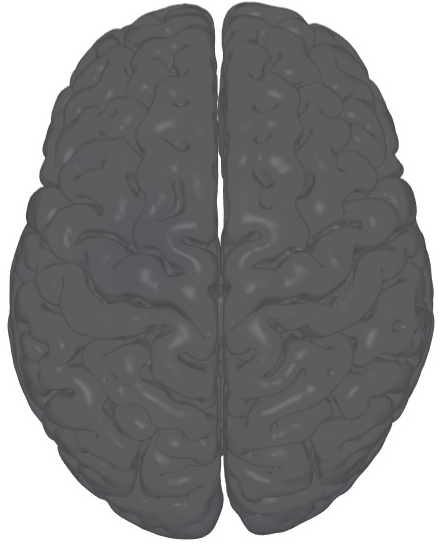}
    \end{minipage}\begin{minipage}[b]{\imageWidth}
        \includegraphics[trim={4cm 2cm 20cm 9.8cm},clip,width=0.8\linewidth]{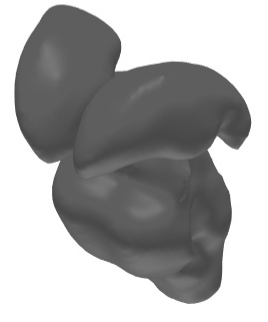}
    \end{minipage}\begin{minipage}[b]{\imageWidth}
        \includegraphics[trim={4.2cm 0.5cm 19.5cm 9cm},clip,width=0.64\linewidth]{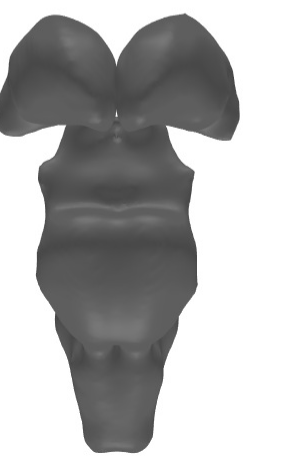}
    \end{minipage}\begin{minipage}[b]{\imageWidth}
        \includegraphics[trim={2cm 0cm 18cm 8.5cm},clip,width=\linewidth]{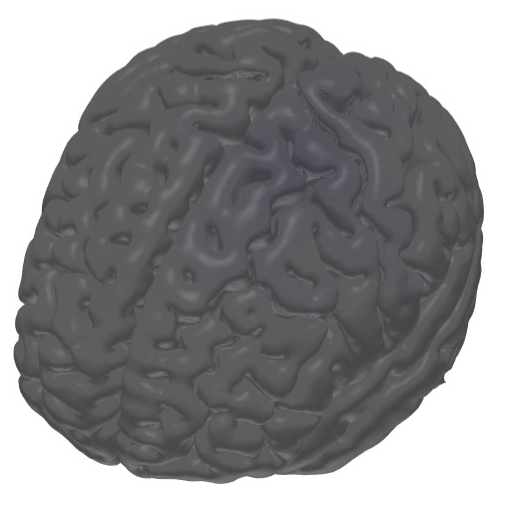}
    \end{minipage}\begin{minipage}[b]{\imageWidth}
        \includegraphics[trim={0cm 0cm 21cm 8cm},clip,width=0.8\linewidth]{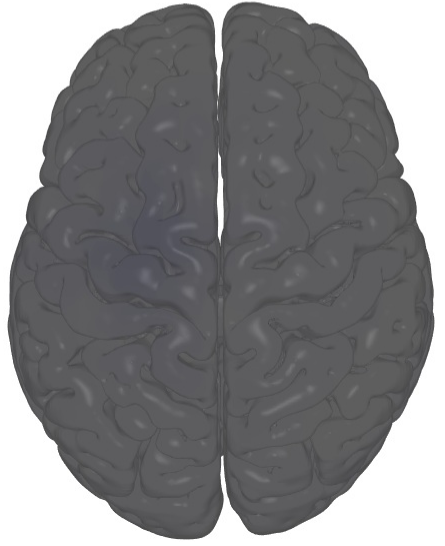}
    \end{minipage}\begin{minipage}[b]{\imageWidth}
        \includegraphics[trim={4cm 2cm 20cm 9.8cm},clip,width=0.8\linewidth]{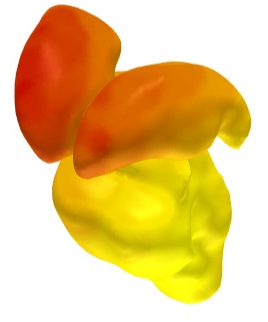}
    \end{minipage}\begin{minipage}[b]{\imageWidth}
        \includegraphics[trim={4.2cm 0.5cm 19.5cm 9cm},clip,width=0.64\linewidth]{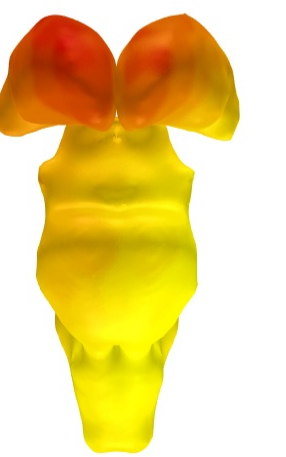}
    \end{minipage}

    \begin{minipage}[b]{\RotTextWidth}
        \rotatebox{90}{\hspace{0.5cm} 20 \unit{\milli\second}}
    \end{minipage}\begin{minipage}[b]{\imageWidth}
        \includegraphics[trim={2cm 0cm 18cm 8.5cm},clip,width=\linewidth]{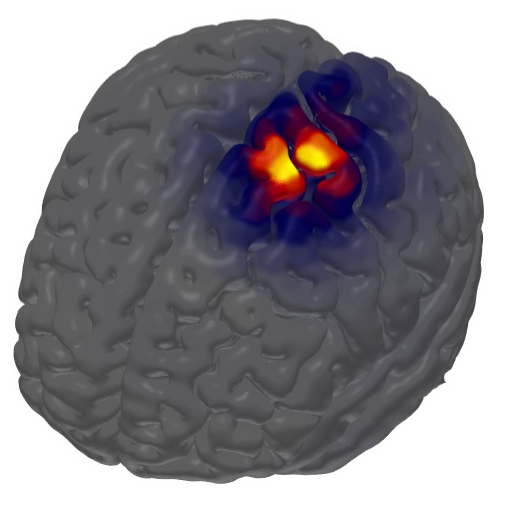}
    \end{minipage}\begin{minipage}[b]{\imageWidth}
        \includegraphics[trim={0cm 0cm 21cm 8cm},clip,width=0.8\linewidth]{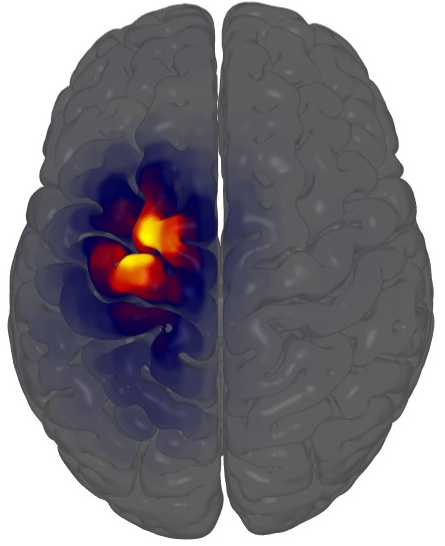}
    \end{minipage}\begin{minipage}[b]{\imageWidth}
        \includegraphics[trim={4cm 2cm 20cm 9.8cm},clip,width=0.8\linewidth]{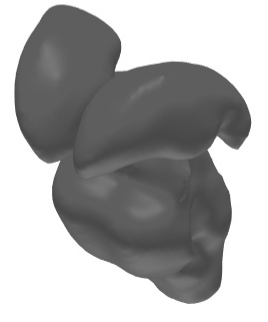}
    \end{minipage}\begin{minipage}[b]{\imageWidth}
        \includegraphics[trim={4.2cm 0.5cm 19.5cm 9cm},clip,width=0.64\linewidth]{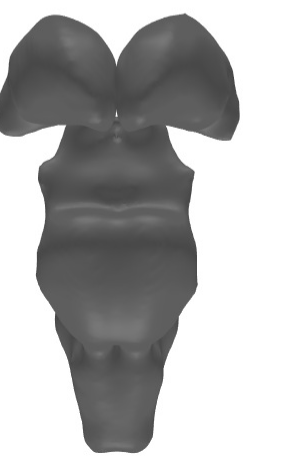}
    \end{minipage}\begin{minipage}[b]{\imageWidth}
        \includegraphics[trim={2cm 0cm 18cm 8.5cm},clip,width=\linewidth]{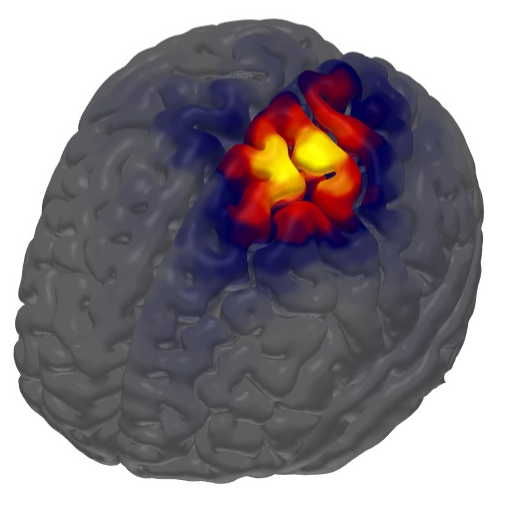}
    \end{minipage}\begin{minipage}[b]{\imageWidth}
        \includegraphics[trim={0cm 0cm 21cm 8cm},clip,width=0.8\linewidth]{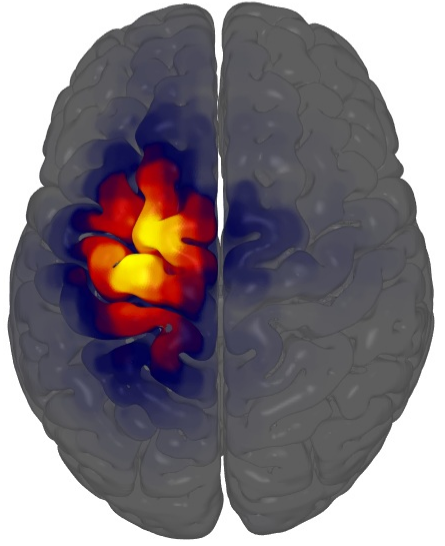}
    \end{minipage}\begin{minipage}[b]{\imageWidth}
        \includegraphics[trim={4cm 2cm 20cm 9.8cm},clip,width=0.8\linewidth]{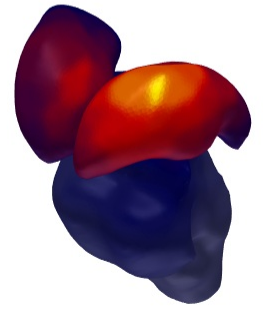}
    \end{minipage}\begin{minipage}[b]{\imageWidth}
        \includegraphics[trim={4.2cm 0.5cm 19.5cm 9cm},clip,width=0.64\linewidth]{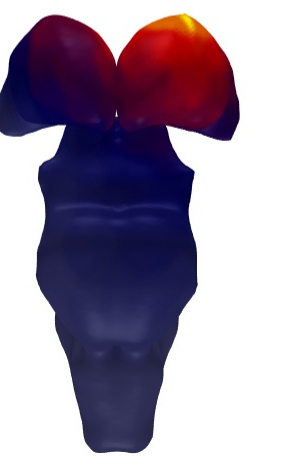}
    \end{minipage}

    \begin{minipage}[b]{\RotTextWidth}
        \rotatebox{90}{\hspace{0.5cm} 22 \unit{\milli\second}}
    \end{minipage}\begin{minipage}[b]{\imageWidth}
        \includegraphics[trim={2cm 0cm 18cm 8.5cm},clip,width=\linewidth]{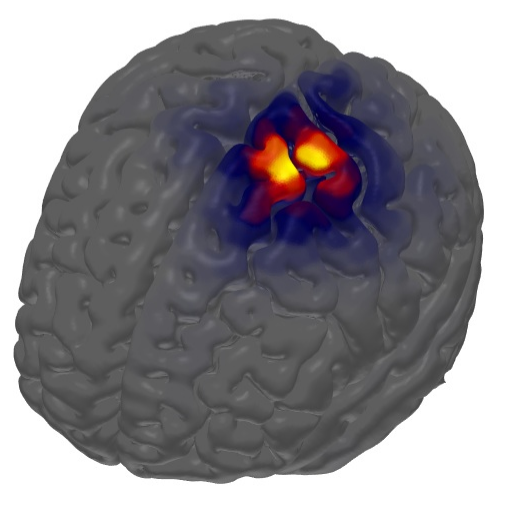}
    \end{minipage}\begin{minipage}[b]{\imageWidth}
        \includegraphics[trim={0cm 0cm 21cm 8cm},clip,width=0.8\linewidth]{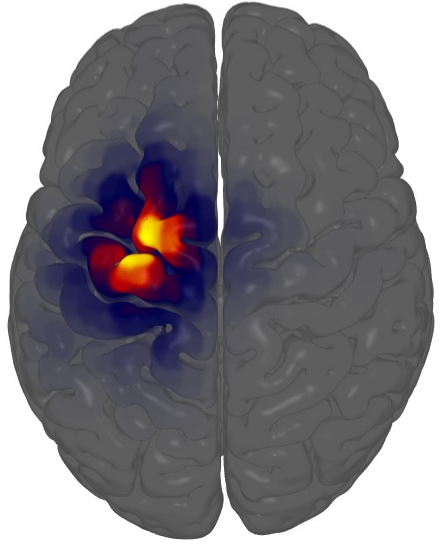}
    \end{minipage}\begin{minipage}[b]{\imageWidth}
        \includegraphics[trim={4cm 2cm 20cm 9.8cm},clip,width=0.8\linewidth]{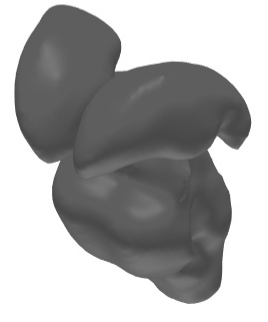}
    \end{minipage}\begin{minipage}[b]{\imageWidth}
        \includegraphics[trim={4.2cm 0.5cm 19.5cm 9cm},clip,width=0.64\linewidth]{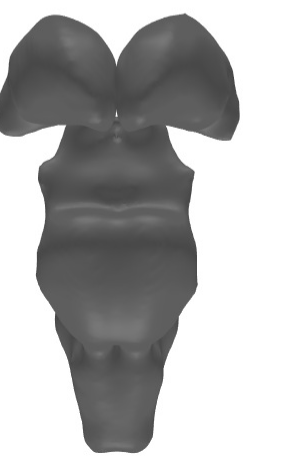}
    \end{minipage}\begin{minipage}[b]{\imageWidth}
        \includegraphics[trim={2cm 0cm 18cm 8.5cm},clip,width=\linewidth]{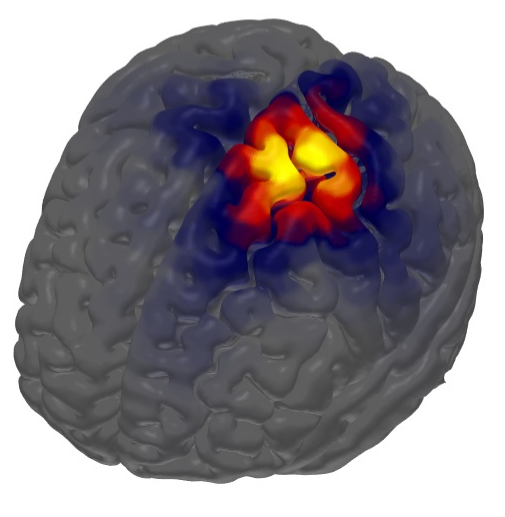}
    \end{minipage}\begin{minipage}[b]{\imageWidth}
        \includegraphics[trim={0cm 0cm 21cm 8cm},clip,width=0.8\linewidth]{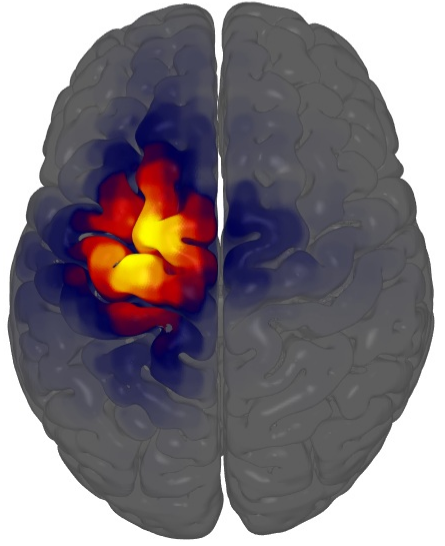}
    \end{minipage}\begin{minipage}[b]{\imageWidth}
        \includegraphics[trim={4cm 2cm 20cm 9.8cm},clip,width=0.8\linewidth]{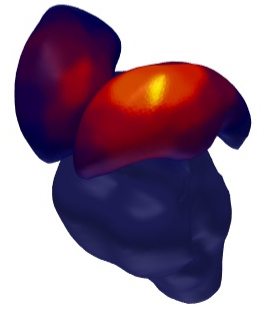}
    \end{minipage}\begin{minipage}[b]{\imageWidth}
        \includegraphics[trim={4.2cm 0.5cm 19.5cm 9cm},clip,width=0.64\linewidth]{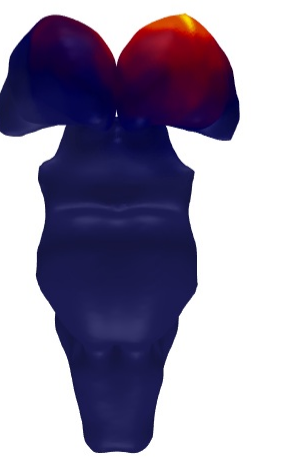}
    \end{minipage}

    \par\noindent\rule{\linewidth}{0.8pt}
    
    \hspace{0.13\linewidth}{\bf DTI-KF}\hspace{0.42\linewidth}{\bf DTI-SKF}
    
    \begin{minipage}[b]{\RotTextWidth}
        \rotatebox{90}{\hspace{0.5cm} 14 \unit{\milli\second}}
    \end{minipage}\begin{minipage}[b]{\imageWidth}
        \includegraphics[trim={2cm 0cm 18cm 8.5cm},clip,width=\linewidth]{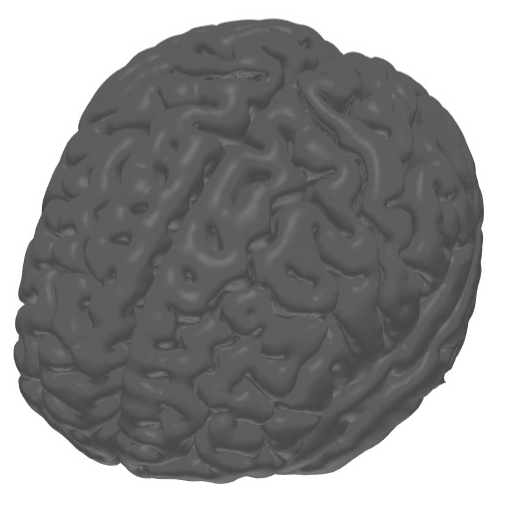}
    \end{minipage}\begin{minipage}[b]{\imageWidth}
        \includegraphics[trim={0cm 0cm 21cm 8cm},clip,width=0.8\linewidth]{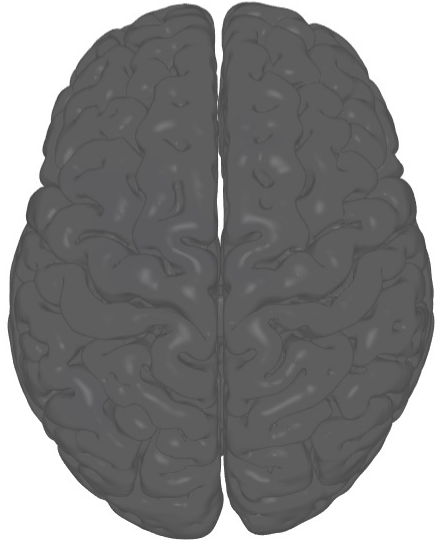}
    \end{minipage}\begin{minipage}[b]{\imageWidth}
        \includegraphics[trim={4cm 2cm 20cm 9.8cm},clip,width=0.8\linewidth]{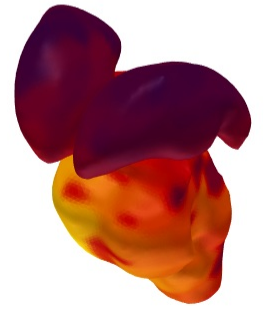}
    \end{minipage}\begin{minipage}[b]{\imageWidth}
        \includegraphics[trim={4.2cm 0.5cm 19.5cm 9cm},clip,width=0.64\linewidth]{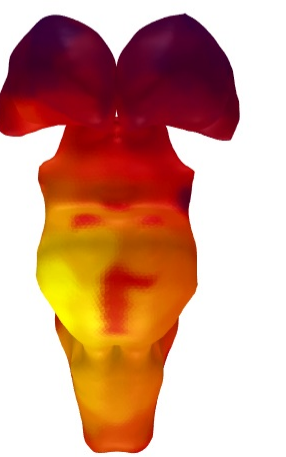}
    \end{minipage}\begin{minipage}[b]{\imageWidth}
        \includegraphics[trim={2cm 0cm 18cm 8.5cm},clip,width=\linewidth]{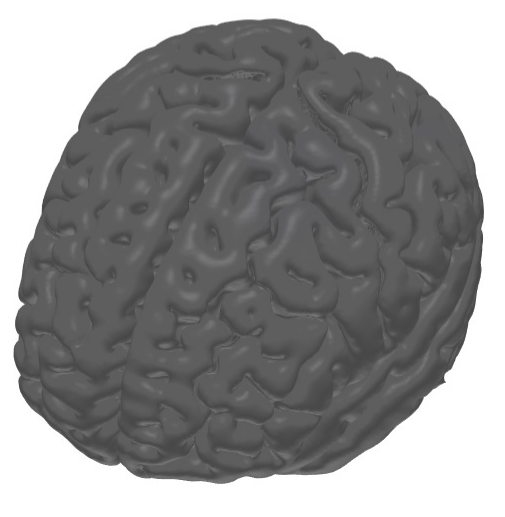}
    \end{minipage}\begin{minipage}[b]{\imageWidth}
        \includegraphics[trim={0cm 0cm 21cm 8cm},clip,width=0.8\linewidth]{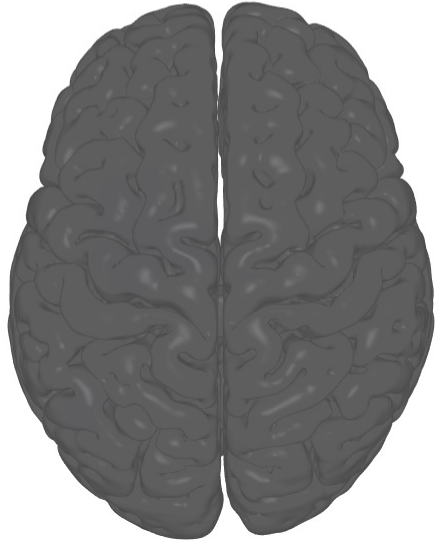}
    \end{minipage}\begin{minipage}[b]{\imageWidth}
        \includegraphics[trim={4cm 2cm 20cm 9.8cm},clip,width=0.8\linewidth]{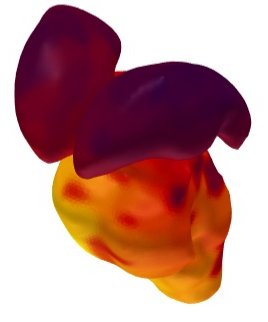}
    \end{minipage}\begin{minipage}[b]{\imageWidth}
        \includegraphics[trim={4.2cm 0.5cm 19.5cm 9cm},clip,width=0.64\linewidth]{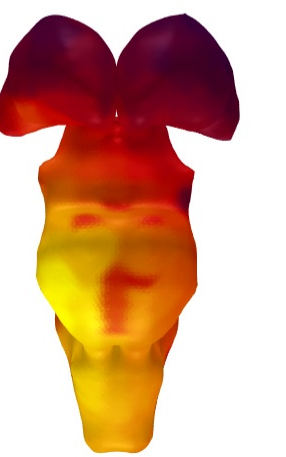}
    \end{minipage}

    \begin{minipage}[b]{\RotTextWidth}
        \rotatebox{90}{\hspace{0.5cm} 16 \unit{\milli\second}}
    \end{minipage}\begin{minipage}[b]{\imageWidth}
        \includegraphics[trim={2cm 0cm 18cm 8.5cm},clip,width=\linewidth]{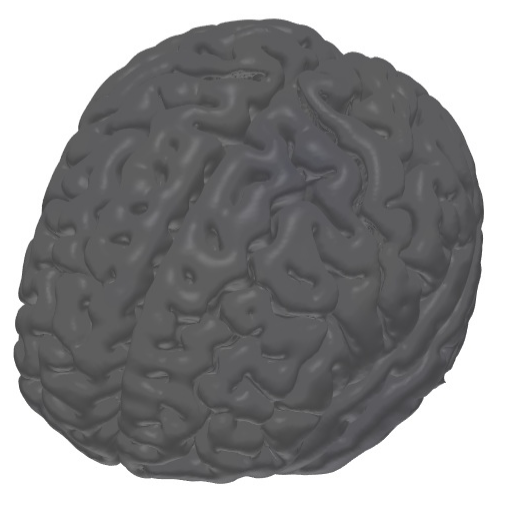}
    \end{minipage}\begin{minipage}[b]{\imageWidth}
        \includegraphics[trim={0cm 0cm 21cm 8cm},clip,width=0.8\linewidth]{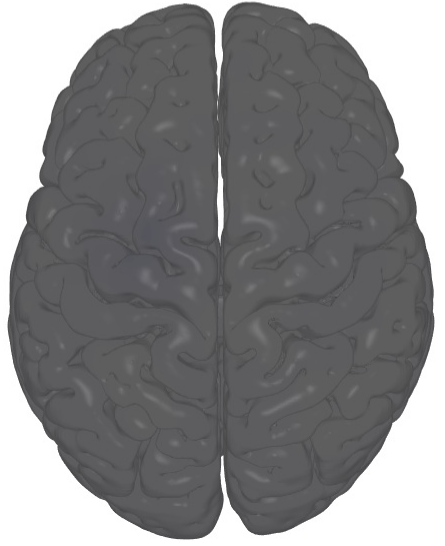}
    \end{minipage}\begin{minipage}[b]{\imageWidth}
        \includegraphics[trim={4cm 2cm 20cm 9.8cm},clip,width=0.8\linewidth]{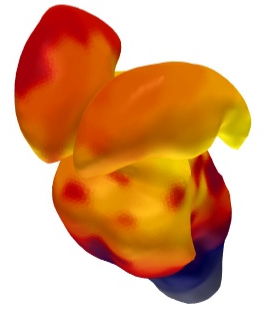}
    \end{minipage}\begin{minipage}[b]{\imageWidth}
        \includegraphics[trim={4.2cm 0.5cm 19.5cm 9cm},clip,width=0.64\linewidth]{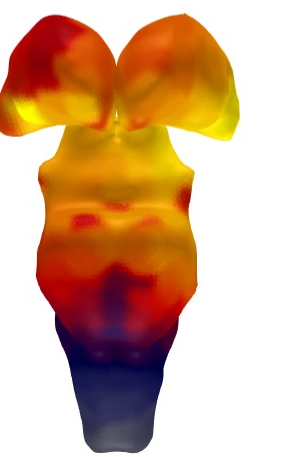}
    \end{minipage}\begin{minipage}[b]{\imageWidth}
        \includegraphics[trim={2cm 0cm 18cm 8.5cm},clip,width=\linewidth]{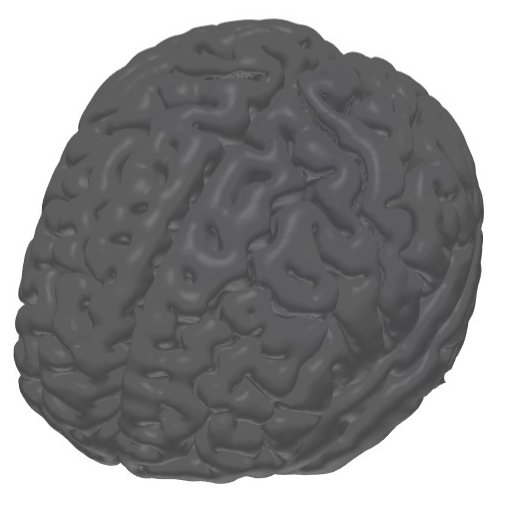}
    \end{minipage}\begin{minipage}[b]{\imageWidth}
        \includegraphics[trim={0cm 0cm 21cm 8cm},clip,width=0.8\linewidth]{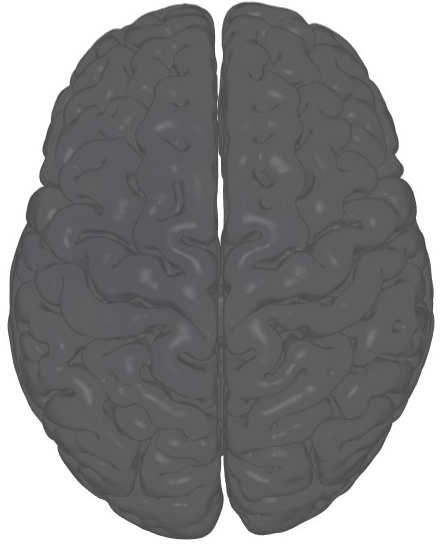}
    \end{minipage}\begin{minipage}[b]{\imageWidth}
        \includegraphics[trim={4cm 2cm 20cm 9.8cm},clip,width=0.8\linewidth]{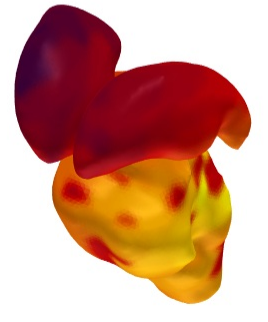}
    \end{minipage}\begin{minipage}[b]{\imageWidth}
        \includegraphics[trim={4.2cm 0.5cm 19.5cm 9cm},clip,width=0.64\linewidth]{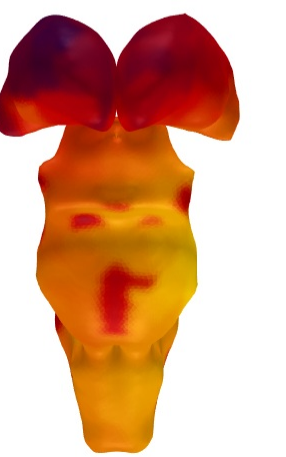}
    \end{minipage}

    \begin{minipage}[b]{\RotTextWidth}
        \rotatebox{90}{\hspace{0.5cm} 20 \unit{\milli\second}}
    \end{minipage}\begin{minipage}[b]{\imageWidth}
        \includegraphics[trim={2cm 0cm 18cm 8.5cm},clip,width=\linewidth]{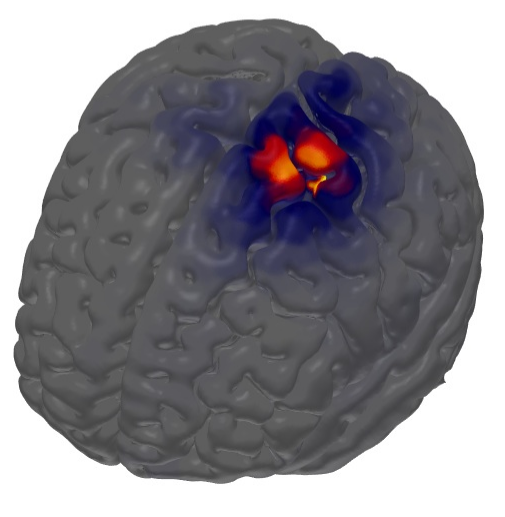}
    \end{minipage}\begin{minipage}[b]{\imageWidth}
        \includegraphics[trim={0cm 0cm 21cm 8cm},clip,width=0.8\linewidth]{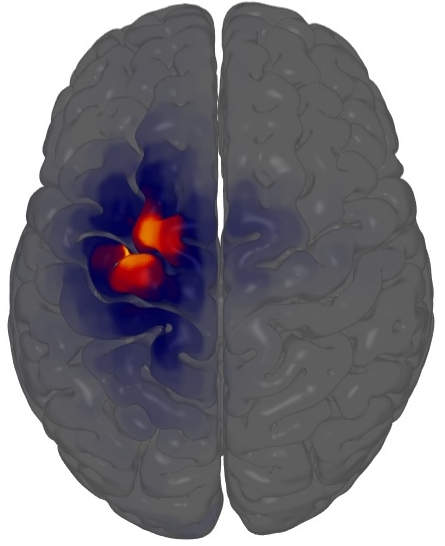}
    \end{minipage}\begin{minipage}[b]{\imageWidth}
        \includegraphics[trim={4cm 2cm 20cm 9.8cm},clip,width=0.8\linewidth]{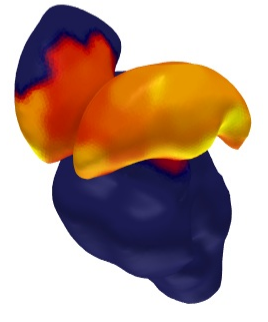}
    \end{minipage}\begin{minipage}[b]{\imageWidth}
        \includegraphics[trim={4.2cm 0.5cm 19.5cm 9cm},clip,width=0.64\linewidth]{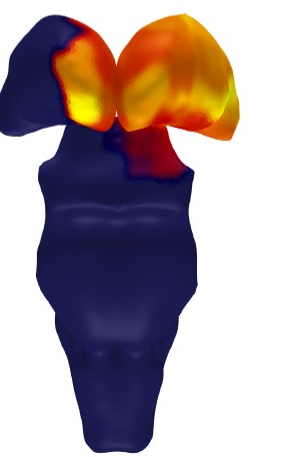}
    \end{minipage}\begin{minipage}[b]{\imageWidth}
        \includegraphics[trim={2cm 0cm 18cm 8.5cm},clip,width=\linewidth]{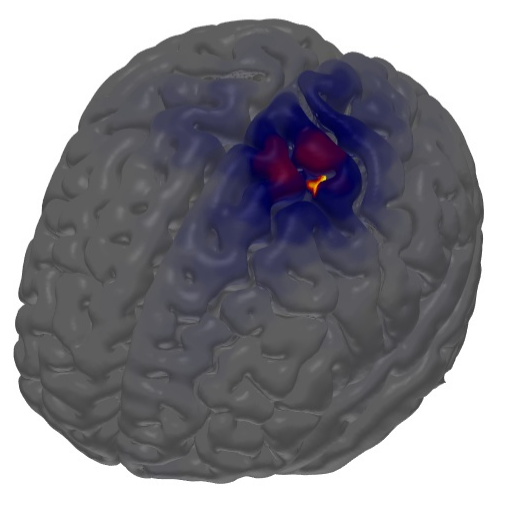}
    \end{minipage}\begin{minipage}[b]{\imageWidth}
        \includegraphics[trim={0cm 0cm 21cm 8cm},clip,width=0.8\linewidth]{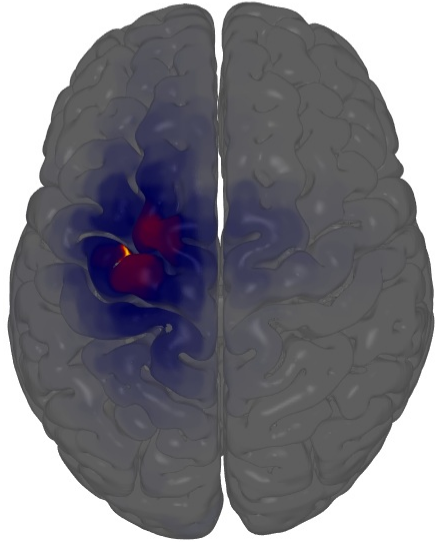}
    \end{minipage}\begin{minipage}[b]{\imageWidth}
        \includegraphics[trim={4cm 2cm 20cm 9.8cm},clip,width=0.8\linewidth]{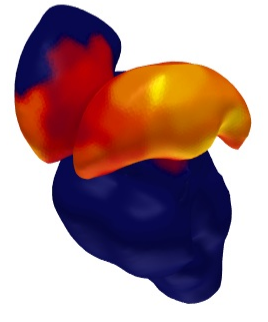}
    \end{minipage}\begin{minipage}[b]{\imageWidth}
        \includegraphics[trim={4.2cm 0.5cm 19.5cm 9cm},clip,width=0.64\linewidth]{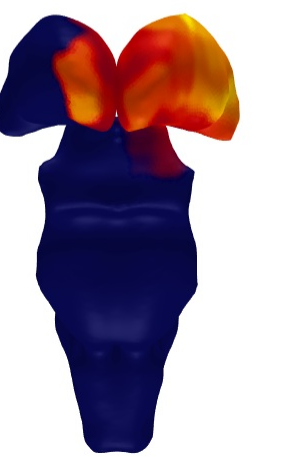}
    \end{minipage}

    \begin{minipage}[b]{\RotTextWidth}
        \rotatebox{90}{\hspace{0.5cm} 22 \unit{\milli\second}}
    \end{minipage}\begin{minipage}[b]{\imageWidth}
        \includegraphics[trim={2cm 0cm 18cm 8.5cm},clip,width=\linewidth]{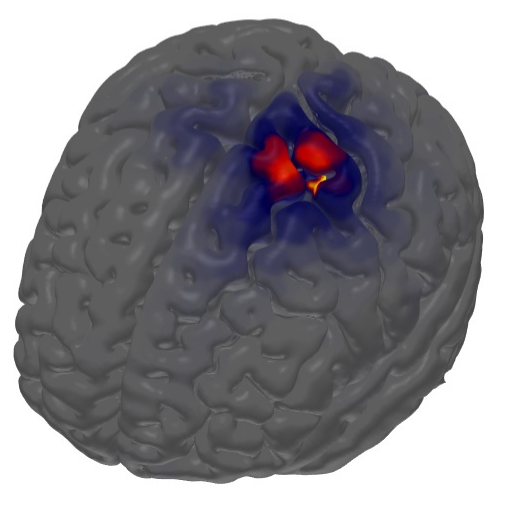}
    \end{minipage}\begin{minipage}[b]{\imageWidth}
        \includegraphics[trim={0cm 0cm 21cm 8cm},clip,width=0.8\linewidth]{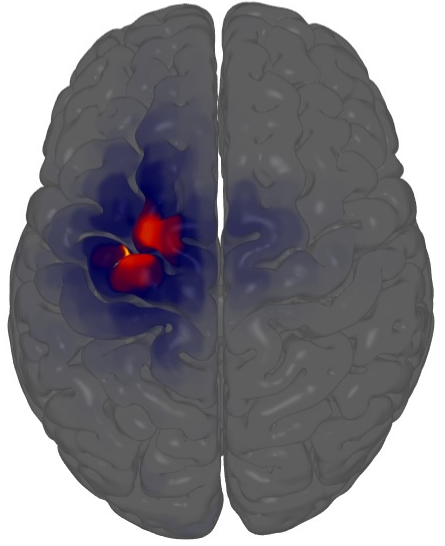}
    \end{minipage}\begin{minipage}[b]{\imageWidth}
        \includegraphics[trim={4cm 2cm 20cm 9.8cm},clip,width=0.8\linewidth]{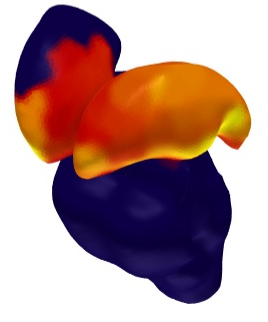}
    \end{minipage}\begin{minipage}[b]{\imageWidth}
        \includegraphics[trim={4.2cm 0.5cm 19.5cm 9cm},clip,width=0.64\linewidth]{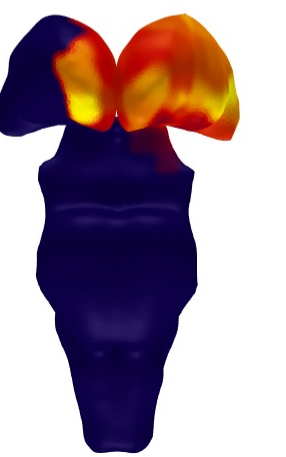}
    \end{minipage}\begin{minipage}[b]{\imageWidth}
        \includegraphics[trim={2cm 0cm 18cm 8.5cm},clip,width=\linewidth]{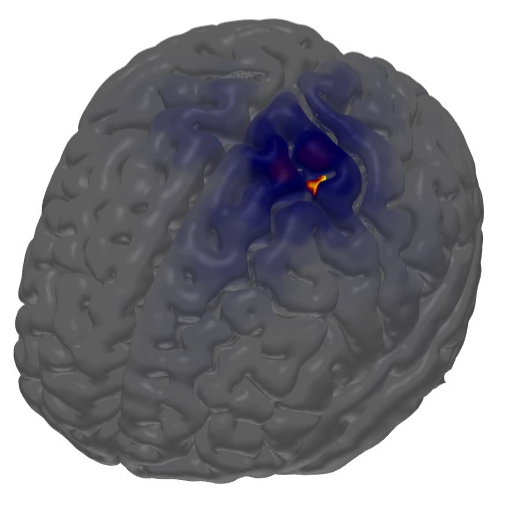}
    \end{minipage}\begin{minipage}[b]{\imageWidth}
        \includegraphics[trim={0cm 0cm 21cm 8cm},clip,width=0.8\linewidth]{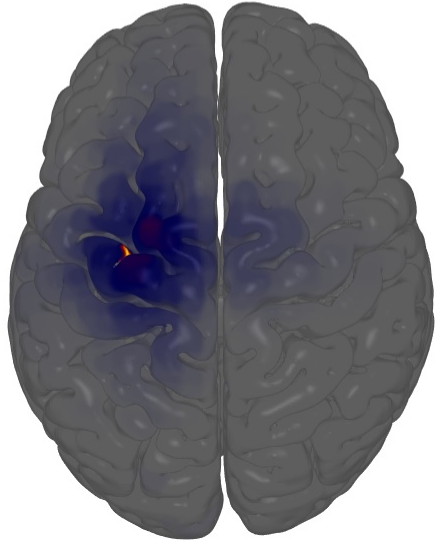}
    \end{minipage}\begin{minipage}[b]{\imageWidth}
        \includegraphics[trim={4cm 2cm 20cm 9.8cm},clip,width=0.8\linewidth]{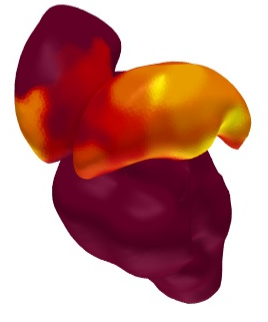}
    \end{minipage}\begin{minipage}[b]{\imageWidth}
        \includegraphics[trim={4.2cm 0.5cm 19.5cm 9cm},clip,width=0.64\linewidth]{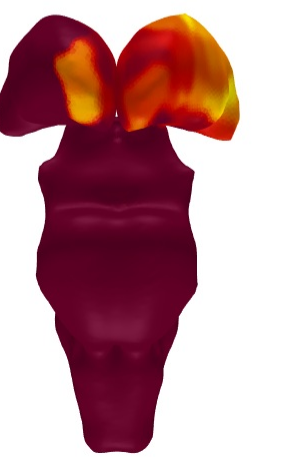}
    \end{minipage}
    
    \caption{Source estimation projections on the brain model for four events appearing in SEP simulation. The upper block shows the result for the Kalman filter and the Standardized Kalman filter using a random walk evolution model, and the lower block shows the corresponding methods that use a DTI-based transition model.}
    \label{fig:SEPbrain}
\end{figure*}

Based on the visual inspection of Figure \ref{fig:SEPbrain}, KF can detect only the cortical activity. There is also no notable difference between estimations of 20 and 22 \unit{\milli\second}, indicating that the method cannot differentiate nearby sources from each other even if their orientation is significantly diverse. The difference from SKF is that, due to standardization, the deep activity can be localized with satisfactory accuracy. Cortical activity is more spread than with KF, and it suffers from similarity issues between clearly distinct cortical components. The later thalamic components are estimated to be too anterior, a distinct phenomenon across all methods analyzed next.

In contrast to KF, DTI-KF can detect deep sources due to the advanced evolution model. The estimation of 14 \unit{\milli\second} exhibits close similarity to the SKF estimate of the said component. 16 \unit{\milli\second} component is missing the pair at the bottom of the brainstem. The cortical estimates at 20 and 22 \unit{\milli\second} are more focal than with KF, and the estimates differ distinctly. However, the highest values are located in the position of the 22 \unit{\milli\second} component, showing that the rising radial component overshadows the tangential 20 \unit{\milli\second} component at the cortex. The thalamic components are more spread out than with SKF. The tail of the highest estimation magnitudes is shifted toward the posterior lateral part of the thalamus, where the true source is located. DTI-SKF provides nearly the same estimation of the 14 \unit{\milli\second} component as DTI-KF. The highest 16 \unit{\milli\second} estimates stay on the left side of the model, unlike SKF. The cortical source estimates are extremely focal and nearly accurate. The thalamic estimates are less spread than DTI-KF but more than SKF. The tailing of the estimation indicates more clearly posterior-lateral thalamic contribution than DTI-KF.

\begin{figure*}[b!]
    \centering
    \begin{minipage}{0.2\linewidth}
    \centering
    {\bf KF}
        \includegraphics[width=\linewidth]{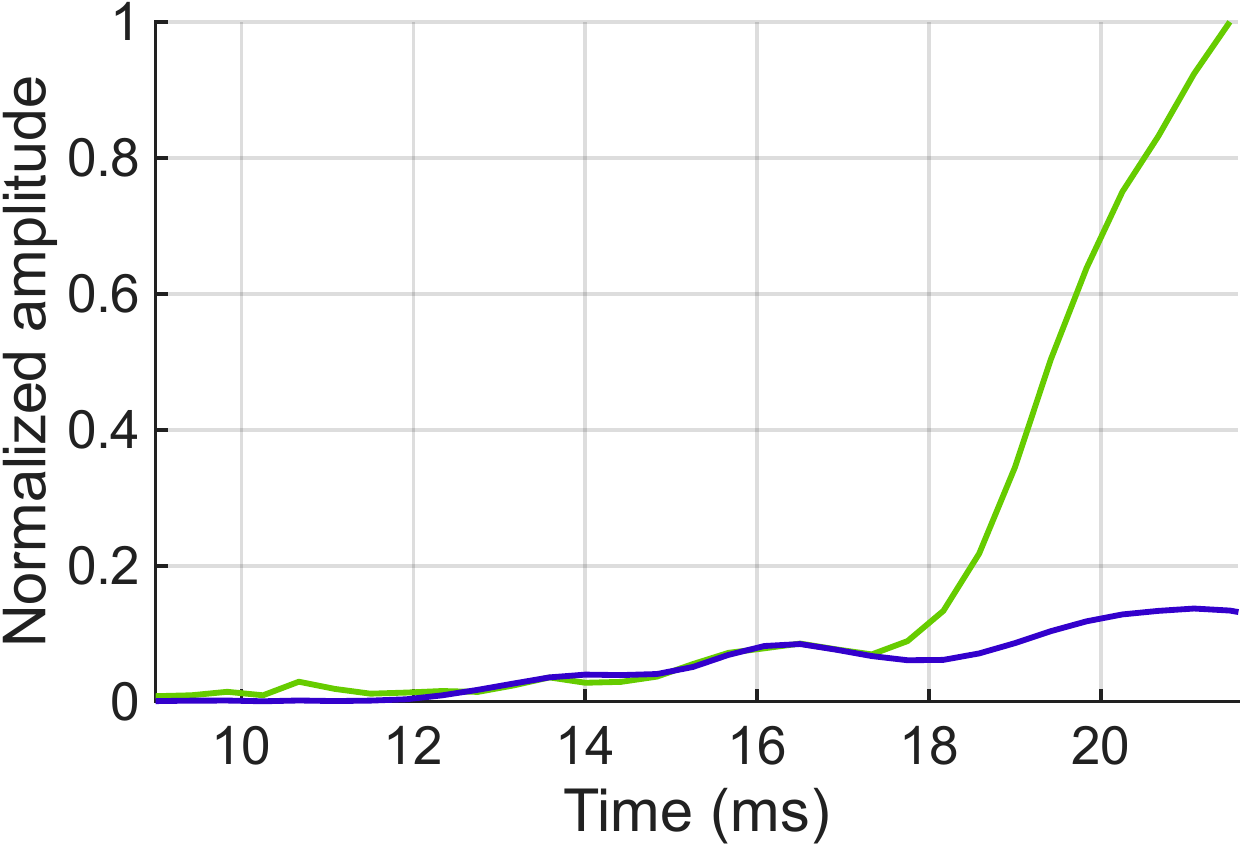}
    \end{minipage}\hspace{0.5cm}\begin{minipage}{0.2\linewidth}
    \centering
    {\bf SKF}
        \includegraphics[width=\linewidth]{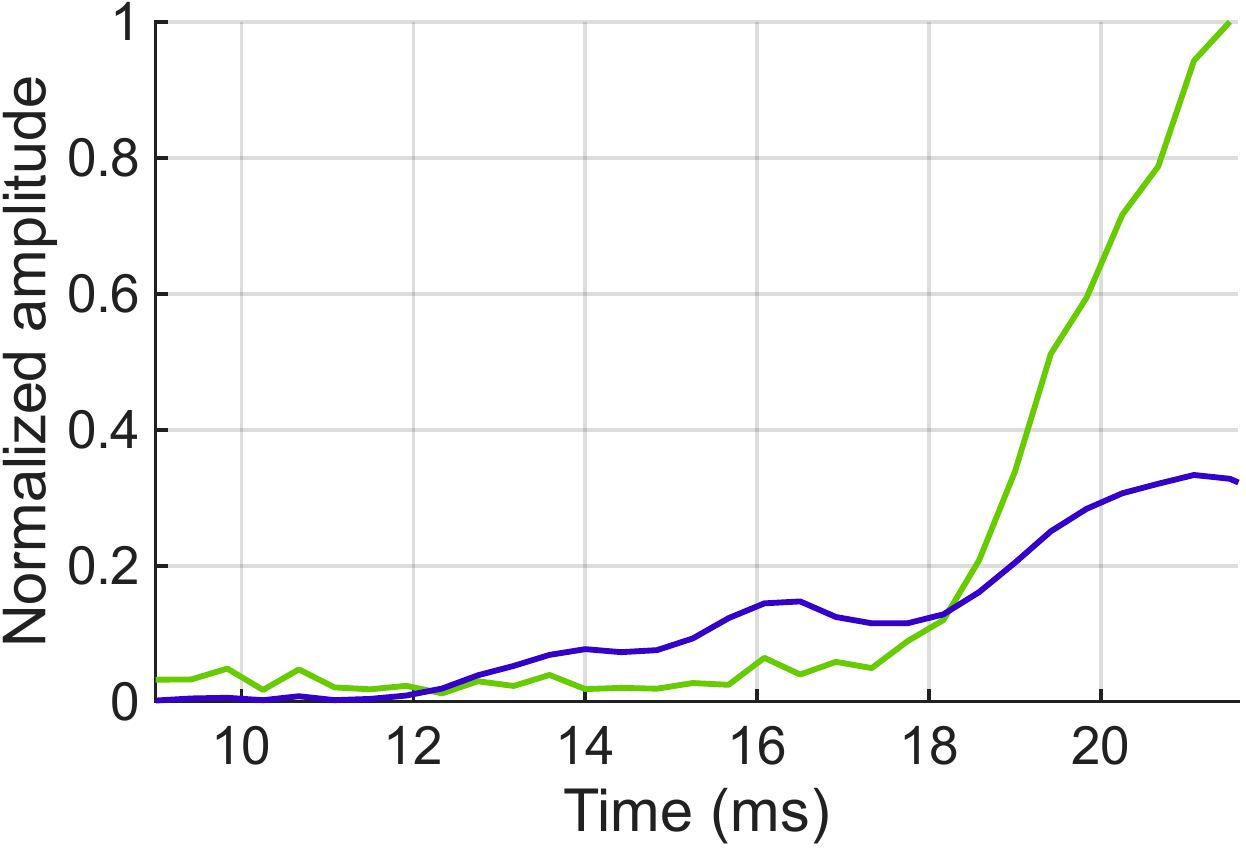}
    \end{minipage}\hspace{0.5cm}\begin{minipage}{0.2\linewidth}
    \centering
    {\bf DTI-KF}
        \includegraphics[width=\linewidth]{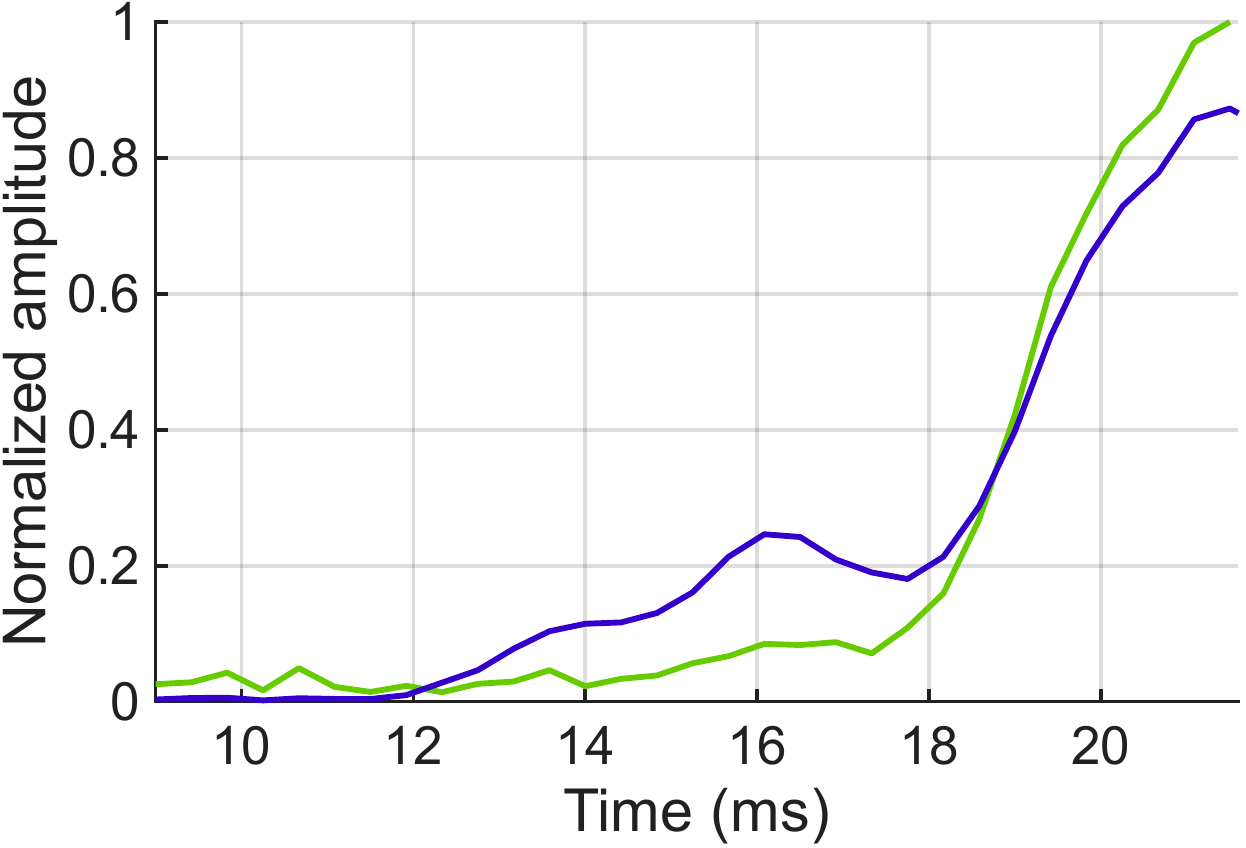}
        \end{minipage}\hspace{0.5cm}\begin{minipage}{0.2\linewidth}
        \centering
    {\bf DTI-SKF}
        \includegraphics[width=\linewidth]{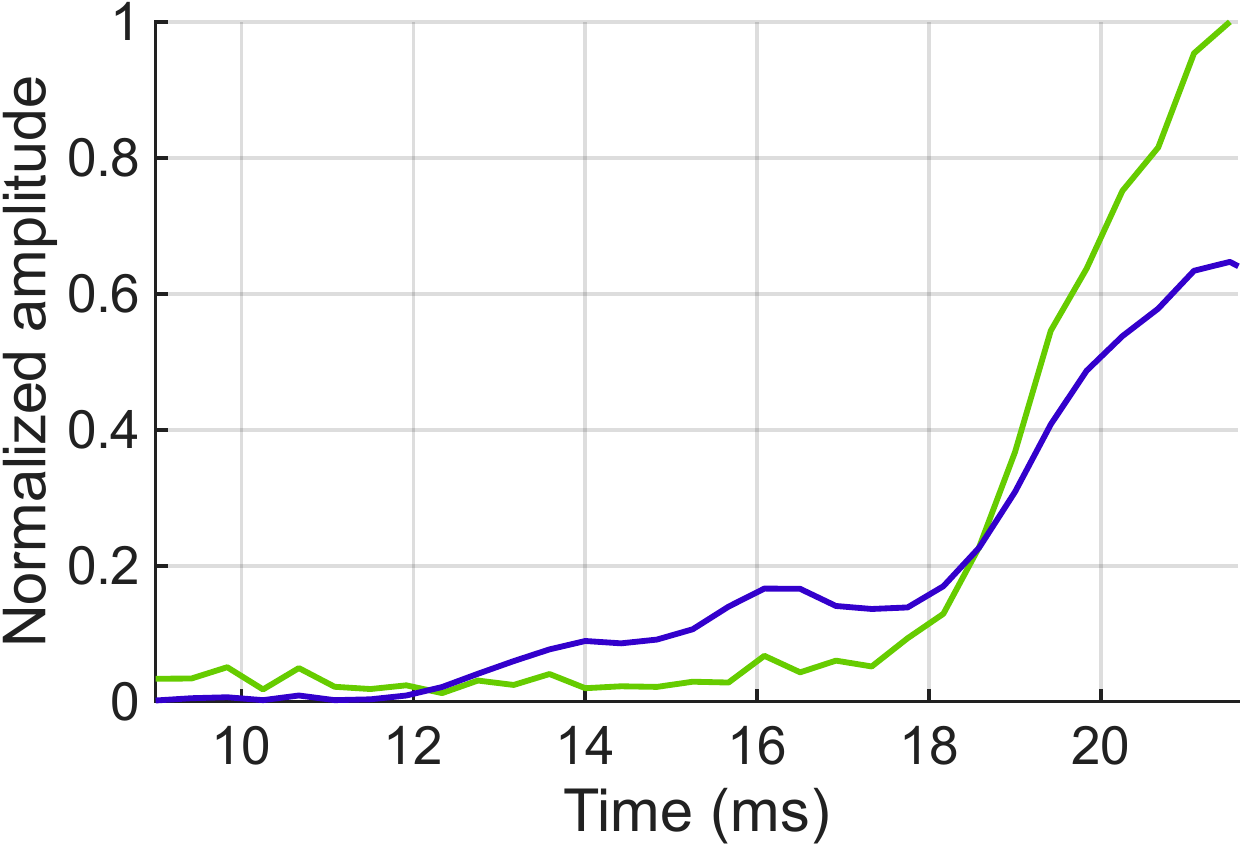}
        \end{minipage}

    \begin{minipage}{0.2\linewidth}
    \centering
        \includegraphics[width=\linewidth]{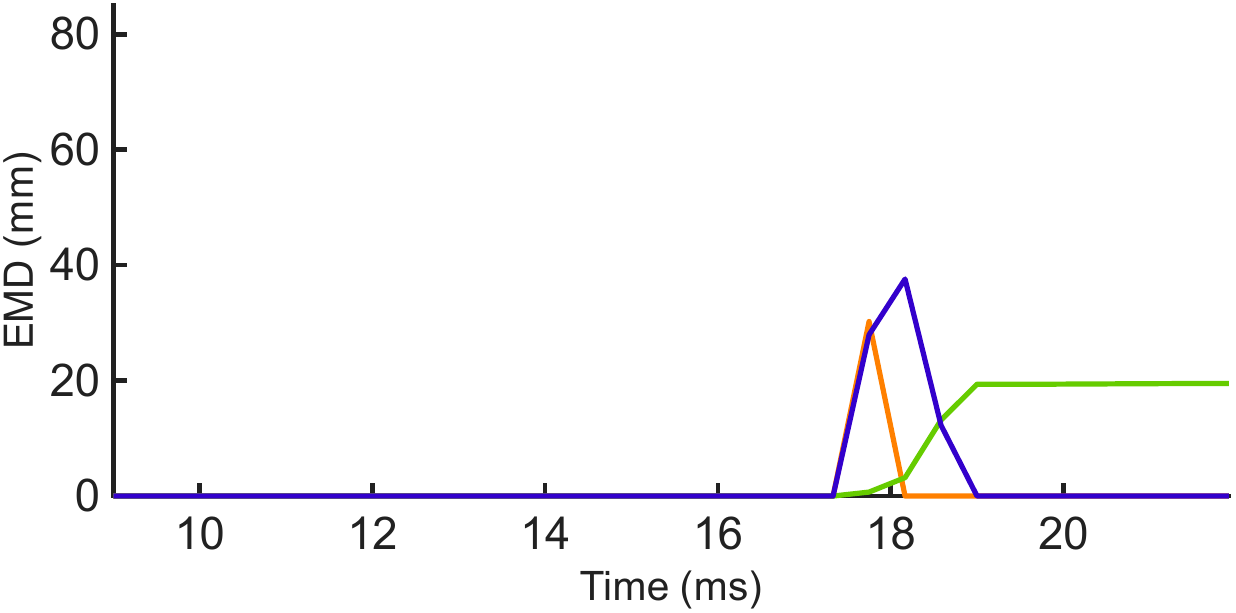}
    \end{minipage}\hspace{0.5cm}\begin{minipage}{0.2\linewidth}
    \centering
        \includegraphics[width=\linewidth]{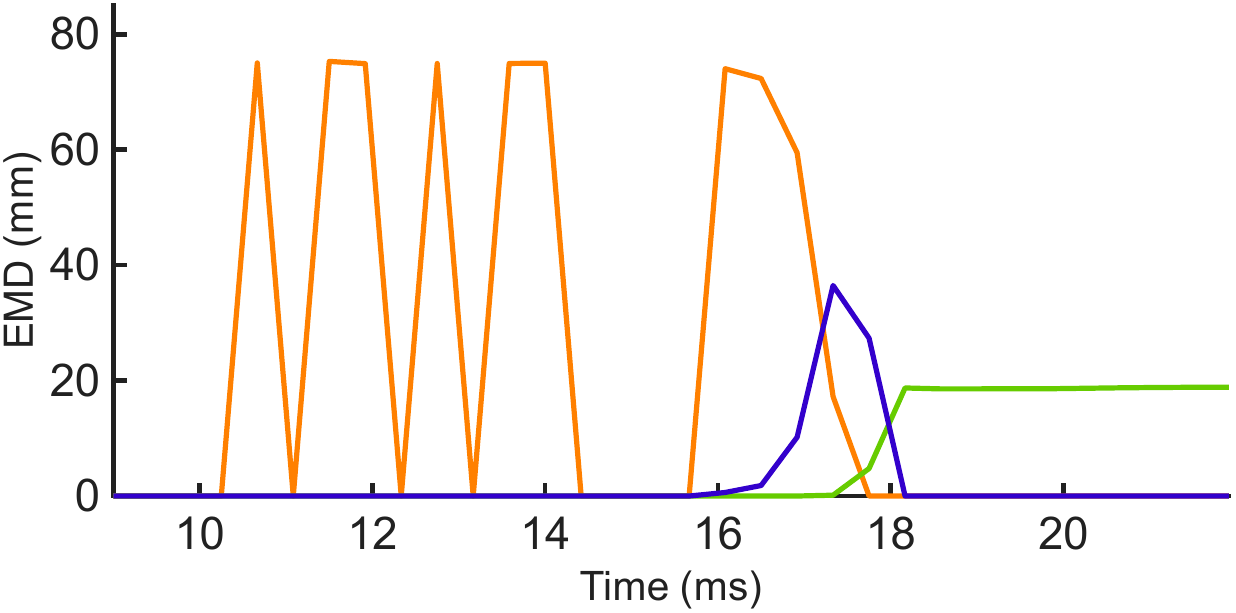}
    \end{minipage}\hspace{0.5cm}\begin{minipage}{0.2\linewidth}
    \centering
        \includegraphics[width=\linewidth]{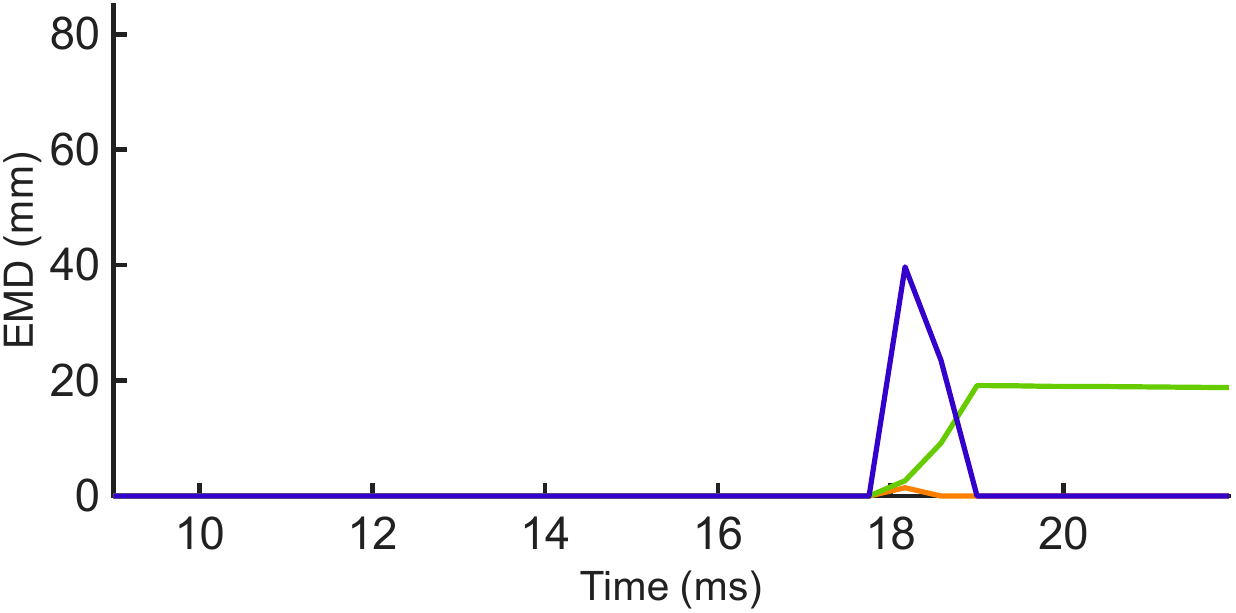}
        \end{minipage}\hspace{0.5cm}\begin{minipage}{0.2\linewidth}
        \centering
        \includegraphics[width=\linewidth]{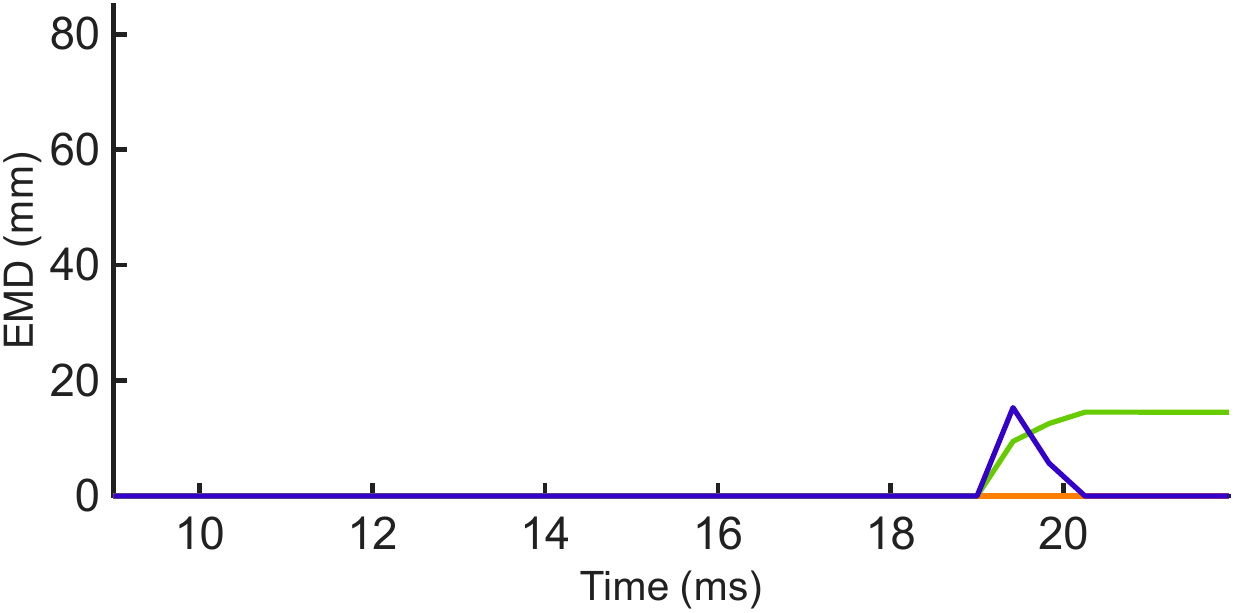}
        \end{minipage}

        \begin{minipage}{0.2\linewidth}
    \centering
        \includegraphics[width=\linewidth]{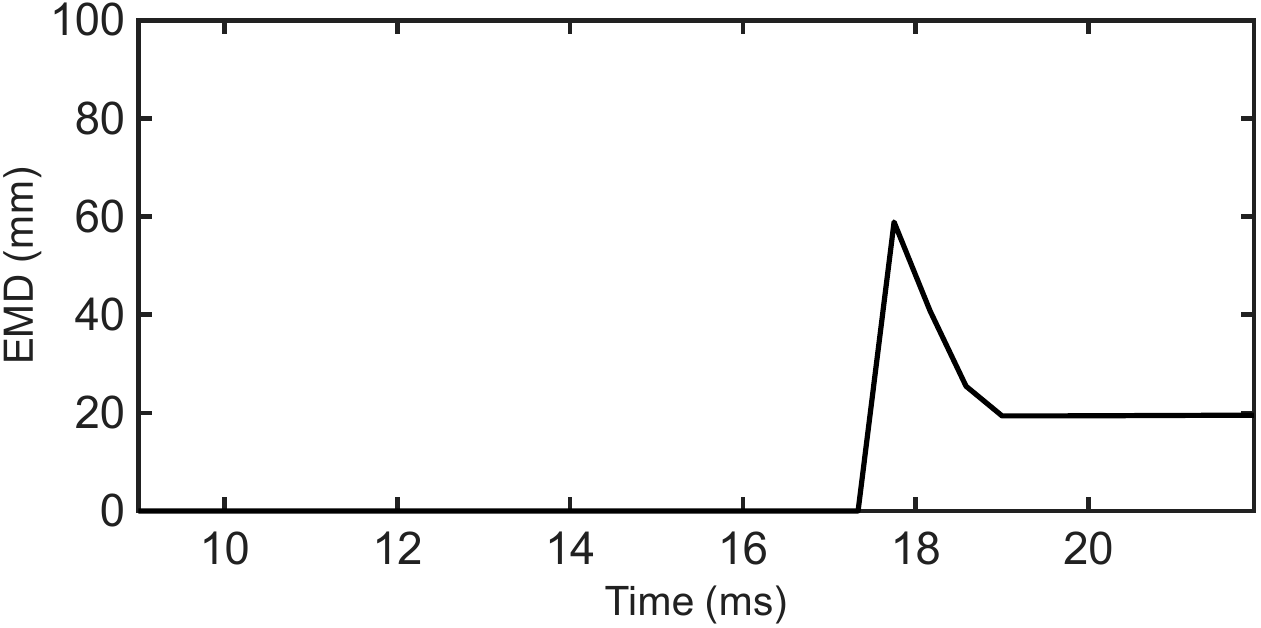}
    \end{minipage}\hspace{0.5cm}\begin{minipage}{0.2\linewidth}
    \centering
        \includegraphics[width=\linewidth]{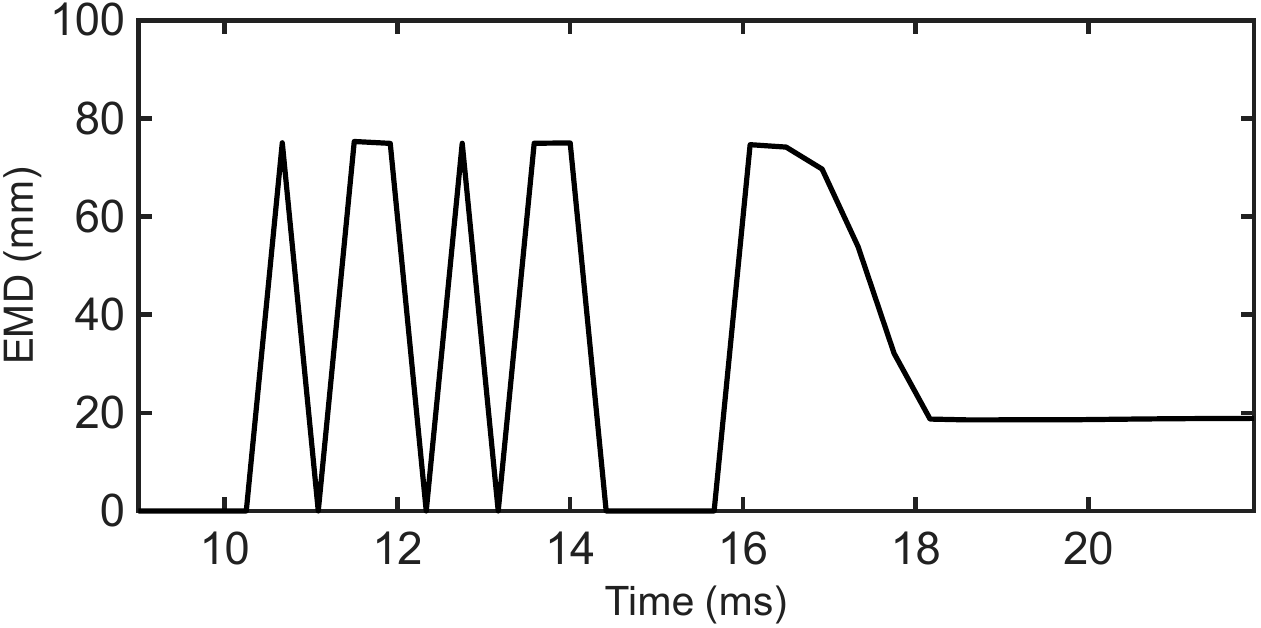}
    \end{minipage}\hspace{0.5cm}\begin{minipage}{0.2\linewidth}
    \centering
        \includegraphics[width=\linewidth]{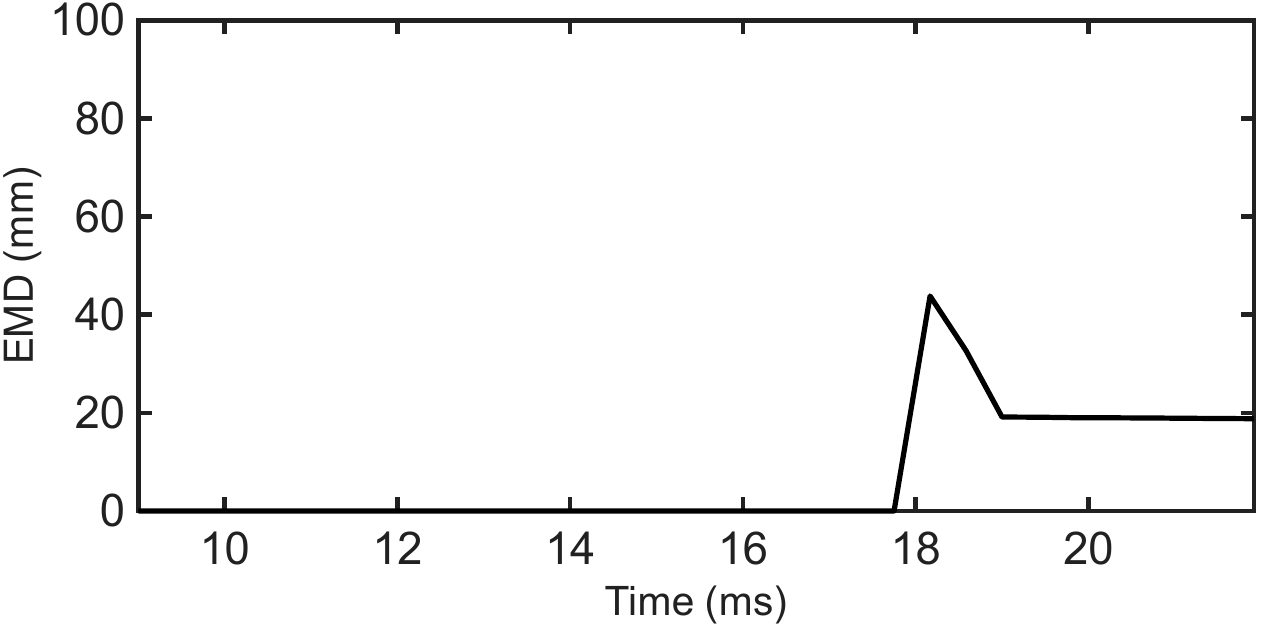}
        \end{minipage}\hspace{0.5cm}\begin{minipage}{0.2\linewidth}
        \centering
        \includegraphics[width=\linewidth]{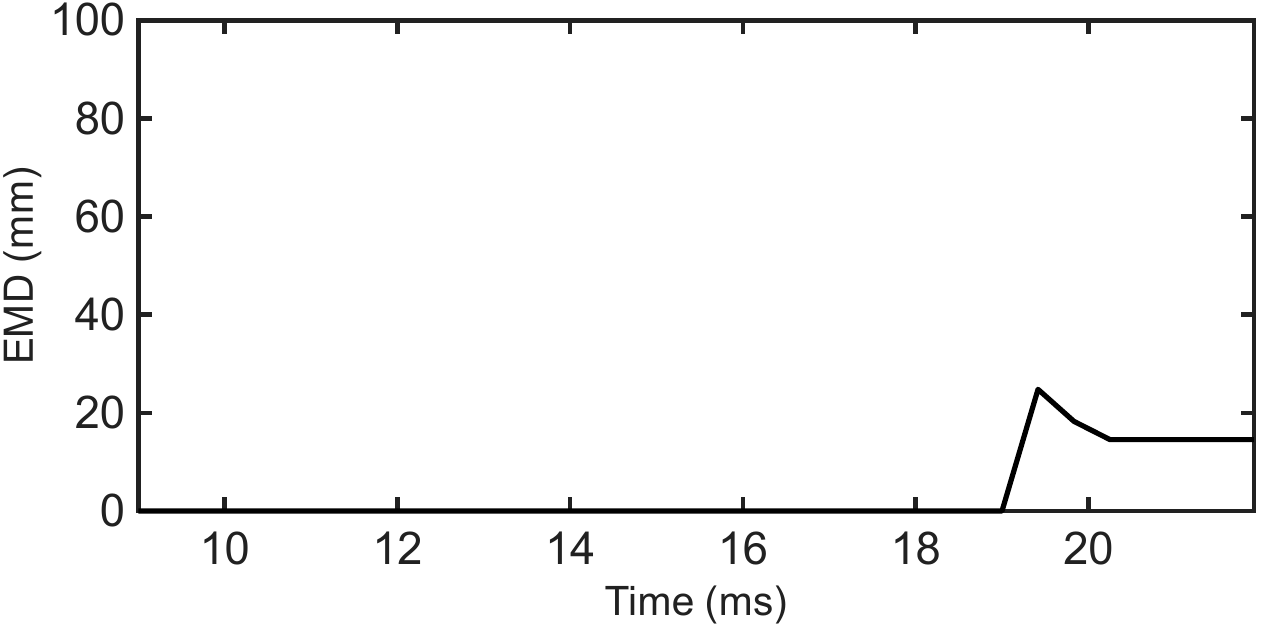}
        \end{minipage}
    \caption{Upper row shows the estimation magnitude time series in the somatosensory cortex (green) and thalamus (blue). The second row displays the Earth Mover's Distance towards the brainstem (orange), thalamus (blue), and somatosensory cortex (green). The last row presents the total Earth Mover's Distance of these three anatomical regions.}
    \label{fig:SEPcurvenEMD}
\end{figure*}

By looking at the time courses of the brain region-specific estimations at the top of Figure \ref{fig:SEPcurvenEMD}, we have a notable ratio difference between cortical and deep estimates between DTI-modelled Kalman filters and their random-walk counterparts. DTI-KF yields the most balanced estimates; otherwise, the curves are nearly the same. Two other image rows show the EMD results measured from the estimations. The middle rows show separately the virtual minimum workload needed to push the reconstruction on the brainstem (orange), left thalamus (blue), and somatosensory cortex. The capacity of each region is dictated by the true time evolution. Moreover, if somatosensory is already "filled", the remaining estimation is transported to deep structures. So early spikes in deep EMD indicate unwanted activity anywhere in the source space, not necessarily deep activity. The total EMD, shown in the bottom row, is more important and indicates the total estimation error. For example, with KF, we see that the deep activity is projected onto the cortex. As is well known, it spreads across the cortex, so activity clusters do not emerge in the visualization. The cortical error appears steady at around 20 \unit{\milli\meter}. With SKF, we see extremely strong spikes in brainstem error because the deep estimation spreads across the deepest structures. The estimations are pushed in the thalamus and somatosensory cortex 2 \unit{\milli\second} earlier than with KF. With DTI-KF, the leakage does not start before 18 \unit{\milli\second}. DTI-KF produces the highest EMD spike in the thalamus across all methods. EMD of the brainstem is nearly zero throughout the duration of the test data. DTI-SKF shows the smallest EMD across all examined regions.

\clearpage

\clearpage
\clearpage

\subsection{Synthesized auditory evoked potential}
In the synthesized AEP data, we have three nearby activities at the insula, superior temporal gyrus, and transverse temporal gyrus on both sides of the brain. The activity alternates with the corresponding gyri on both sides, and the activity starts at the left insula.

\begin{figure*}
\def\BoxWidth{0.25\linewidth}
\def\ImgWidth{0.35\linewidth}
\def\MidImgWidth{0.245\linewidth}
\def\RotTextWidth{0.02\linewidth}
    \centering
    \begin{minipage}[b]{\RotTextWidth}
        \rotatebox{90}{\hspace{-0.8cm}\small Insula L}
    \end{minipage}\begin{minipage}{\BoxWidth}
    \begin{center}
        {\bf KF}
    \begin{minipage}{\ImgWidth}
        \includegraphics[trim={0cm 0cm 19cm 9cm},clip,width=\linewidth]{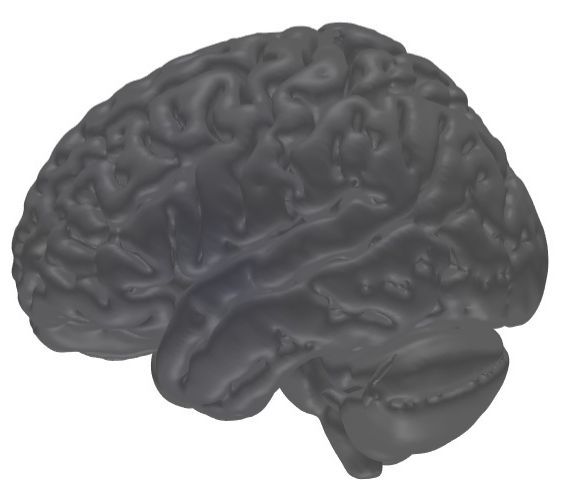}    
    \end{minipage}\begin{minipage}{\MidImgWidth}
    \includegraphics[trim={0cm 0cm 21cm 8cm},clip,width=\linewidth]{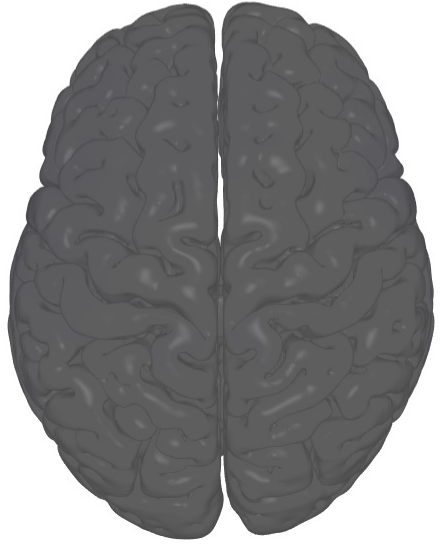}
    \end{minipage}\begin{minipage}{\ImgWidth}
    \includegraphics[trim={0cm 0cm 19cm 9cm},clip,width=\linewidth]{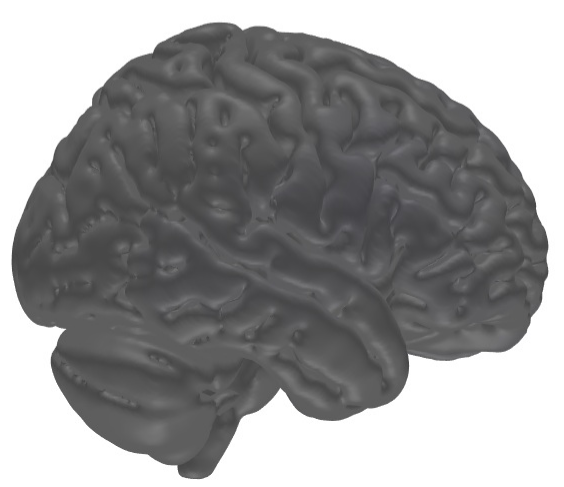}\end{minipage}
    \end{center}
    \end{minipage}\begin{minipage}{\BoxWidth}
    \begin{center}
        {\bf SKF}
    \begin{minipage}{\ImgWidth}
        \includegraphics[trim={0cm 0cm 19cm 9cm},clip,width=\linewidth]{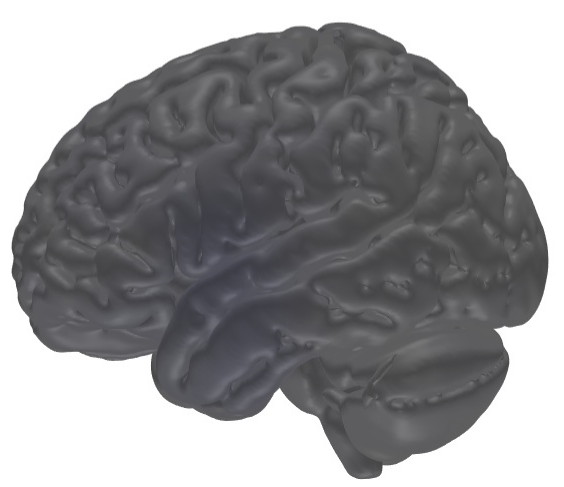}    
    \end{minipage}\begin{minipage}{\MidImgWidth}
    \includegraphics[trim={0cm 0cm 21cm 8cm},clip,width=\linewidth]{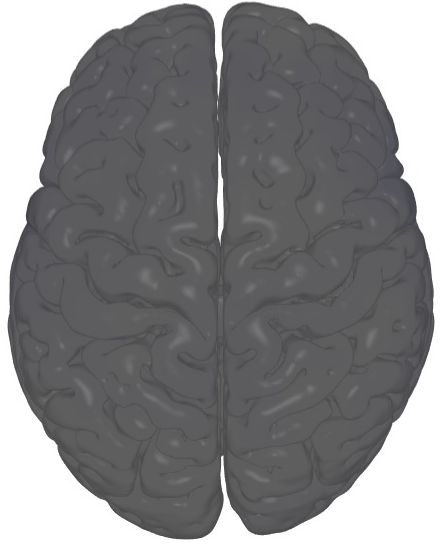}
    \end{minipage}\begin{minipage}{\ImgWidth}
    \includegraphics[trim={0cm 0cm 19cm 9cm},clip,width=\linewidth]{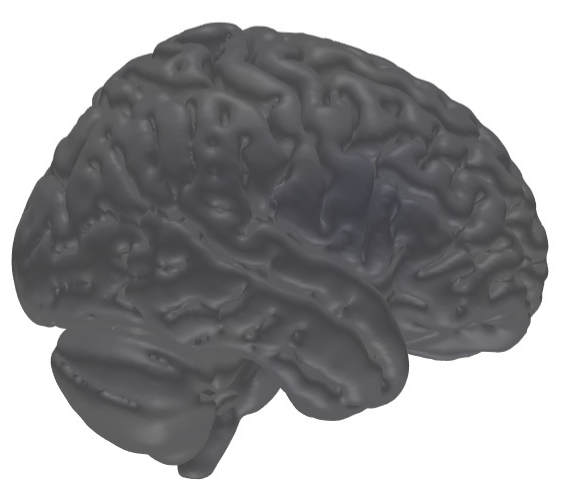}\end{minipage}
    \end{center}
    \end{minipage}\begin{minipage}{\BoxWidth}
    \begin{center}
        {\bf DTI-KF}
    \begin{minipage}{\ImgWidth}
        \includegraphics[trim={0cm 0cm 19cm 9cm},clip,width=\linewidth]{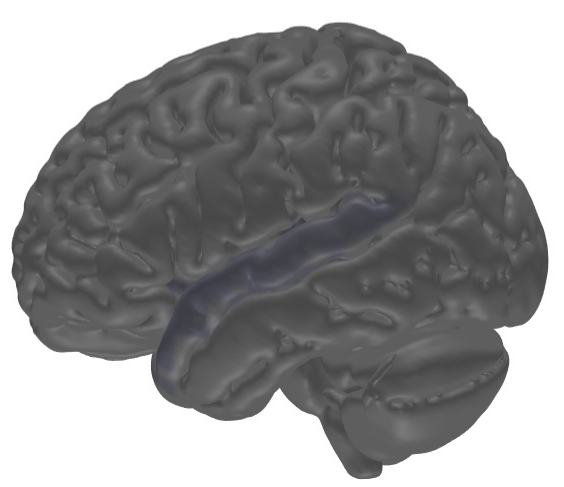}    
    \end{minipage}\begin{minipage}{\MidImgWidth}
    \includegraphics[trim={0cm 0cm 21cm 8cm},clip,width=\linewidth]{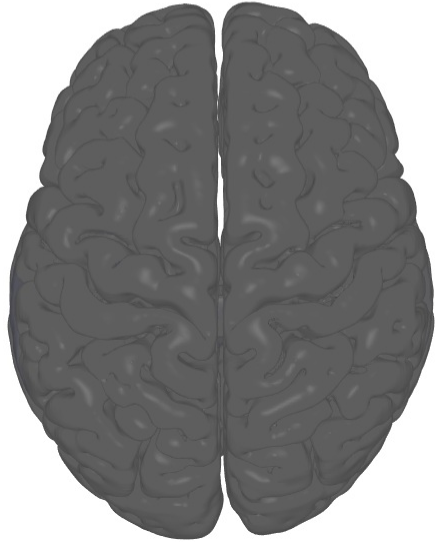}
    \end{minipage}\begin{minipage}{\ImgWidth}
    \includegraphics[trim={0cm 0cm 19cm 9cm},clip,width=\linewidth]{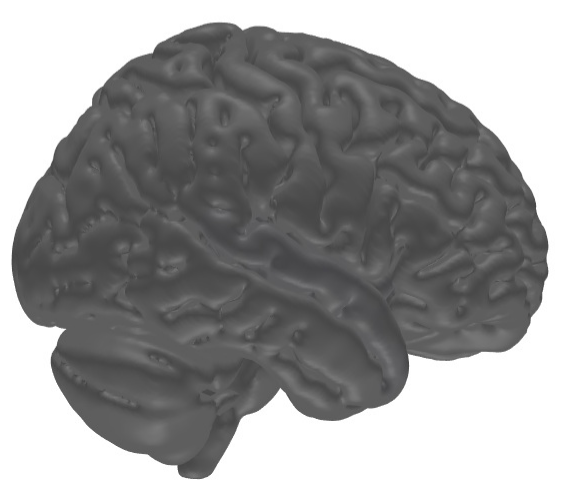}\end{minipage}
    \end{center}
    \end{minipage}\begin{minipage}{\BoxWidth}
    \begin{center}
        {\bf DTI-SKF}
    \begin{minipage}{\ImgWidth}
        \includegraphics[trim={0cm 0cm 19cm 9cm},clip,width=\linewidth]{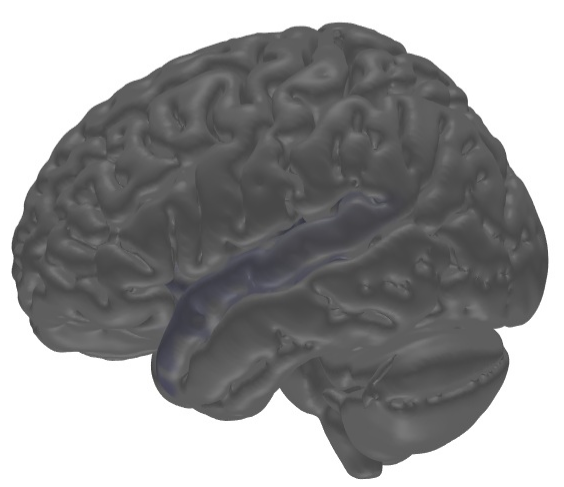}    
    \end{minipage}\begin{minipage}{\MidImgWidth}
    \includegraphics[trim={0cm 0cm 21cm 8cm},clip,width=\linewidth]{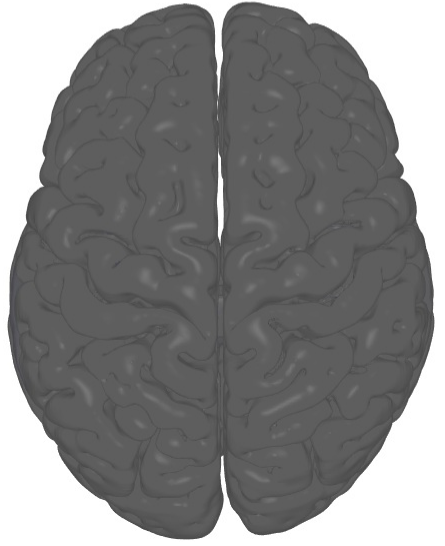}
    \end{minipage}\begin{minipage}{\ImgWidth}
    \includegraphics[trim={0cm 0cm 19cm 9cm},clip,width=\linewidth]{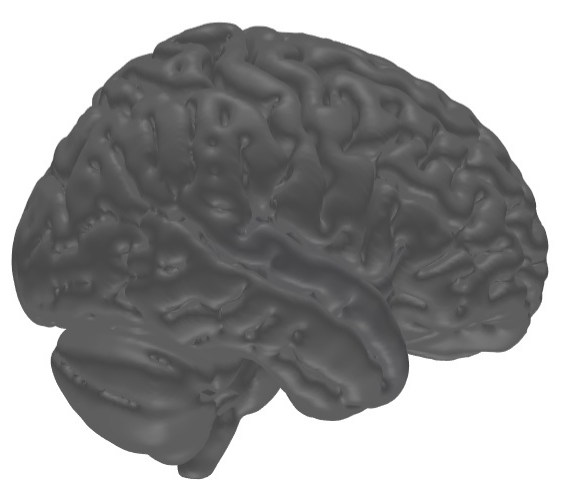}\end{minipage}
    \end{center}
    \end{minipage}

    \begin{minipage}[b]{\RotTextWidth}
        \rotatebox{90}{\hspace{-0.5cm}\small S.T L}
    \end{minipage}\begin{minipage}{\BoxWidth}
    \begin{center}
        {\bf KF}
    \begin{minipage}{\ImgWidth}
        \includegraphics[trim={0cm 0cm 19cm 9cm},clip,width=\linewidth]{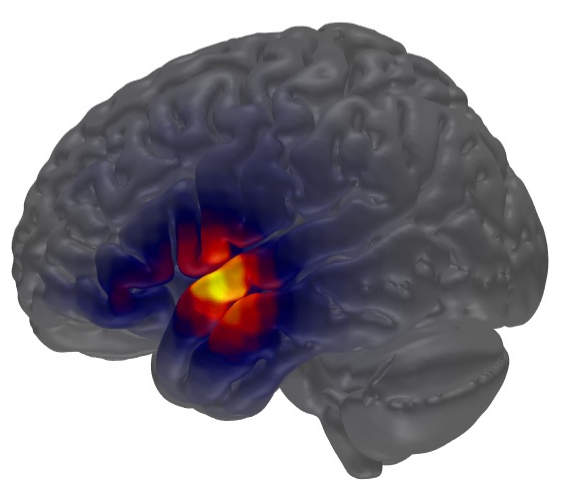}    
    \end{minipage}\begin{minipage}{\MidImgWidth}
    \includegraphics[trim={0cm 0cm 21cm 8cm},clip,width=\linewidth]{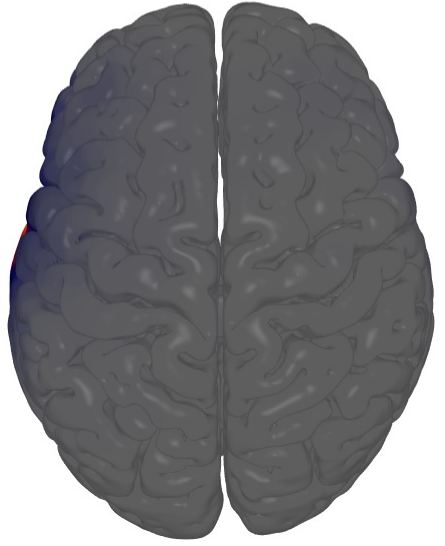}
    \end{minipage}\begin{minipage}{\ImgWidth}
    \includegraphics[trim={0cm 0cm 19cm 9cm},clip,width=\linewidth]{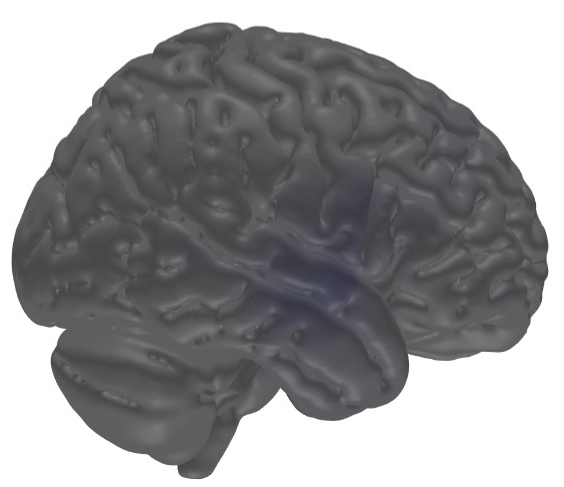}\end{minipage}
    \end{center}
    \end{minipage}\begin{minipage}{\BoxWidth}
    \begin{center}
        {\bf SKF}
    \begin{minipage}{\ImgWidth}
        \includegraphics[trim={0cm 0cm 19cm 9cm},clip,width=\linewidth]{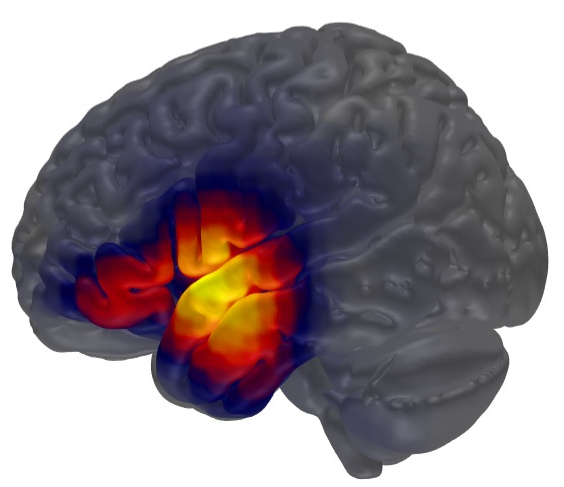}    
    \end{minipage}\begin{minipage}{\MidImgWidth}
    \includegraphics[trim={0cm 0cm 21cm 8cm},clip,width=\linewidth]{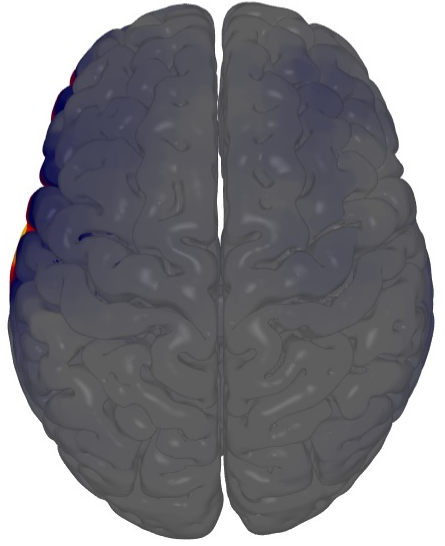}
    \end{minipage}\begin{minipage}{\ImgWidth}
    \includegraphics[trim={0cm 0cm 19cm 9cm},clip,width=\linewidth]{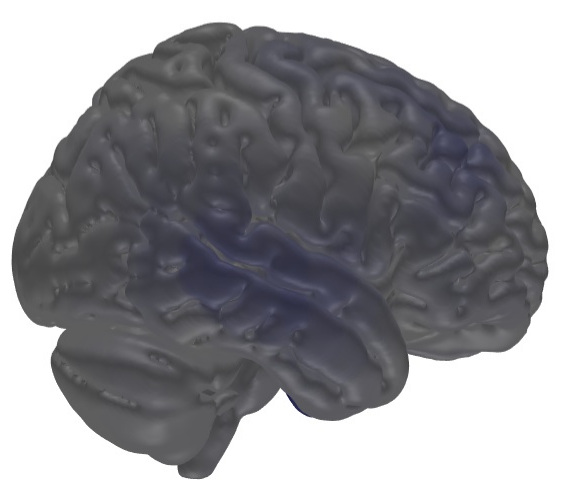}\end{minipage}
    \end{center}
    \end{minipage}\begin{minipage}{\BoxWidth}
    \begin{center}
        {\bf DTI-KF}
    \begin{minipage}{\ImgWidth}
        \includegraphics[trim={0cm 0cm 19cm 9cm},clip,width=\linewidth]{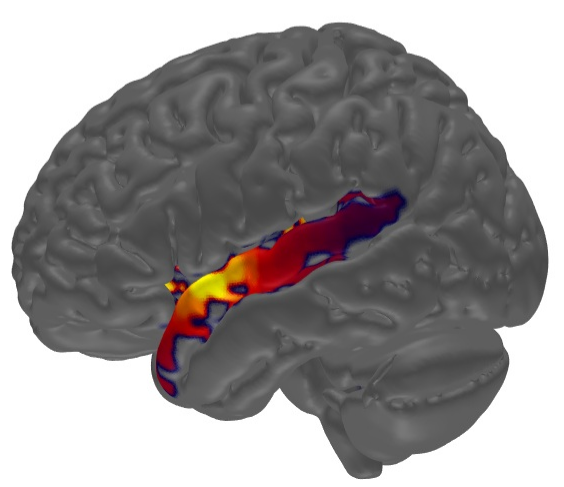}    
    \end{minipage}\begin{minipage}{\MidImgWidth}
    \includegraphics[trim={0cm 0cm 21cm 8cm},clip,width=\linewidth]{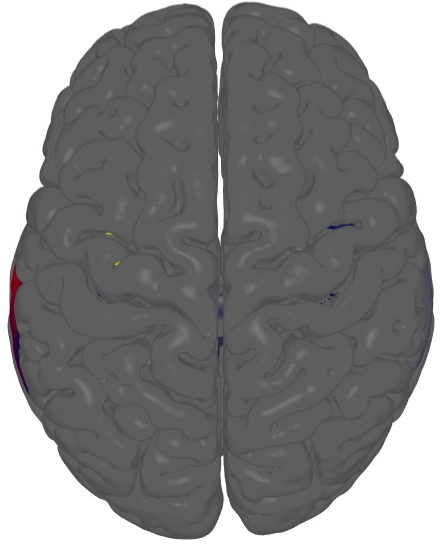}
    \end{minipage}\begin{minipage}{\ImgWidth}
    \includegraphics[trim={0cm 0cm 19cm 9cm},clip,width=\linewidth]{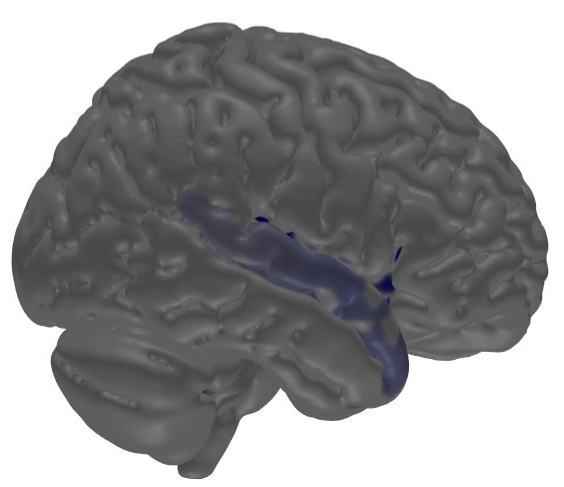}\end{minipage}
    \end{center}
    \end{minipage}\begin{minipage}{\BoxWidth}
    \begin{center}
        {\bf DTI-SKF}
    \begin{minipage}{\ImgWidth}
        \includegraphics[trim={0cm 0cm 19cm 9cm},clip,width=\linewidth]{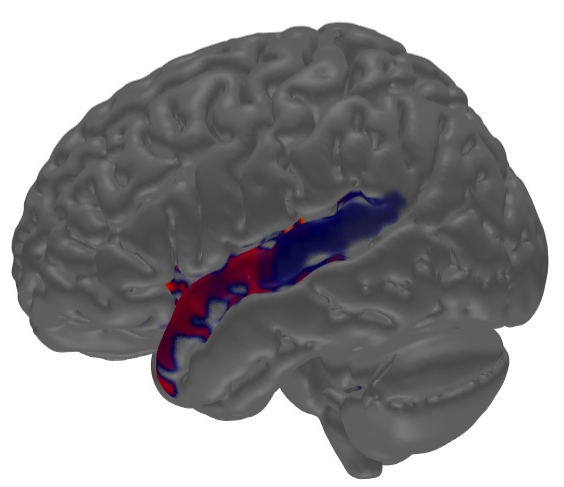}    
    \end{minipage}\begin{minipage}{\MidImgWidth}
    \includegraphics[trim={0cm 0cm 21cm 8cm},clip,width=\linewidth]{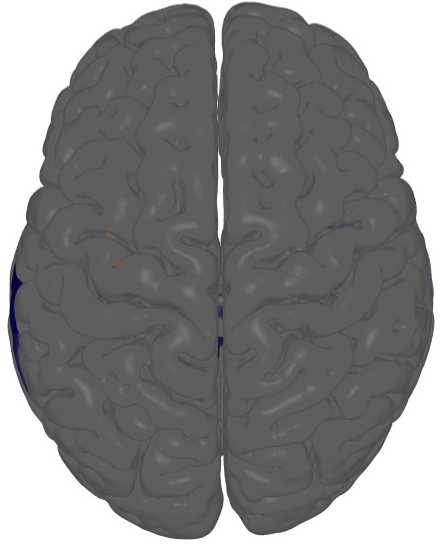}
    \end{minipage}\begin{minipage}{\ImgWidth}
    \includegraphics[trim={0cm 0cm 19cm 9cm},clip,width=\linewidth]{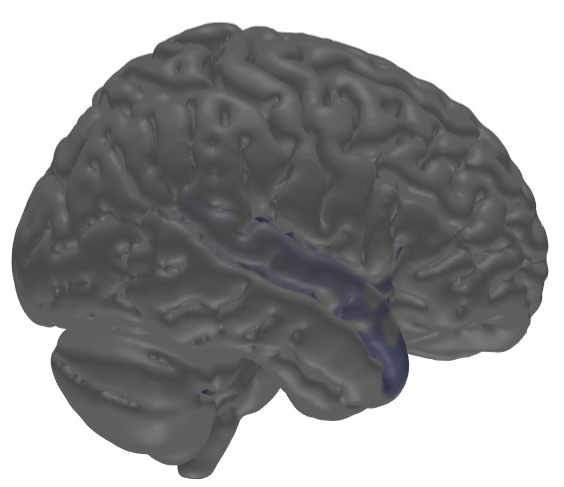}\end{minipage}
    \end{center}
    \end{minipage}

    \begin{minipage}[b]{\RotTextWidth}
        \rotatebox{90}{\hspace{-0.5cm}\small T.T L}
    \end{minipage}\begin{minipage}{\BoxWidth}
    \begin{center}
        {\bf KF}
    \begin{minipage}{\ImgWidth}
        \includegraphics[trim={0cm 0cm 19cm 9cm},clip,width=\linewidth]{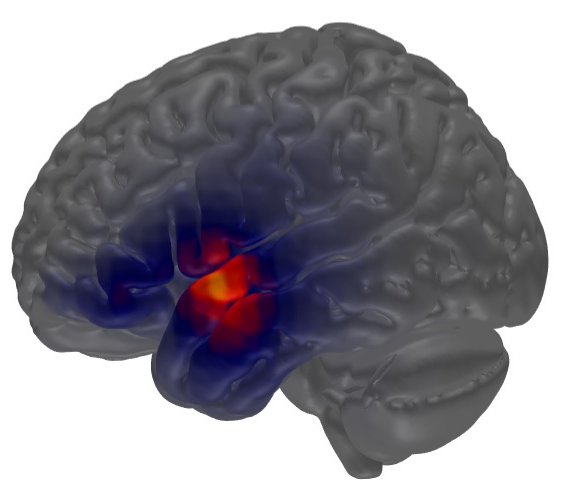}    
    \end{minipage}\begin{minipage}{\MidImgWidth}
    \includegraphics[trim={0cm 0cm 21cm 8cm},clip,width=\linewidth]{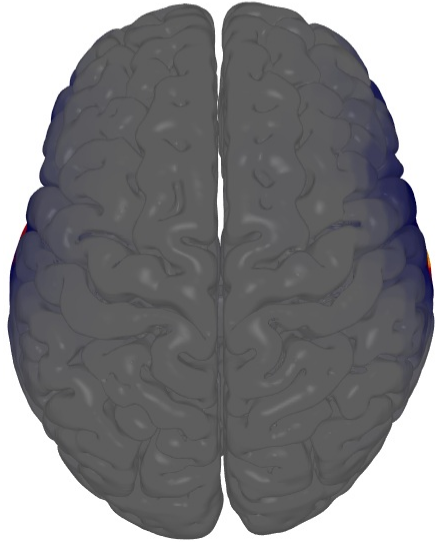}
    \end{minipage}\begin{minipage}{\ImgWidth}
    \includegraphics[trim={0cm 0cm 19cm 9cm},clip,width=\linewidth]{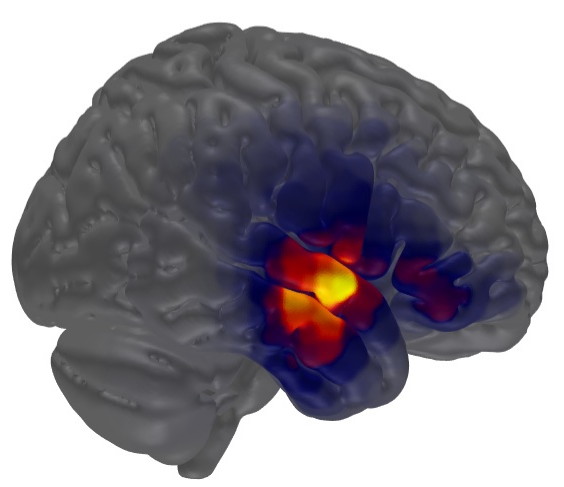}\end{minipage}
    \end{center}
    \end{minipage}\begin{minipage}{\BoxWidth}
    \begin{center}
        {\bf SKF}
    \begin{minipage}{\ImgWidth}
        \includegraphics[trim={0cm 0cm 19cm 9cm},clip,width=\linewidth]{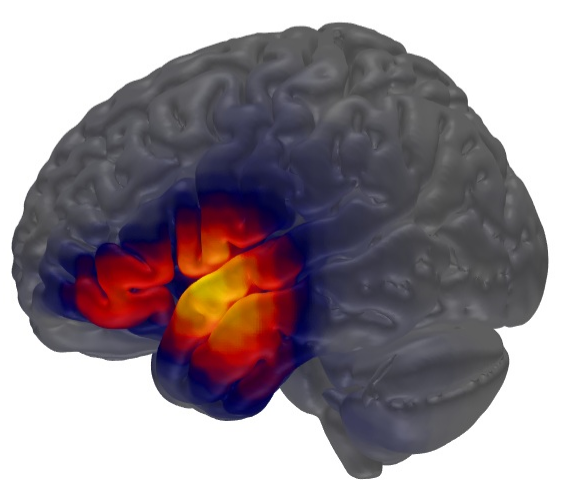}    
    \end{minipage}\begin{minipage}{\MidImgWidth}
    \includegraphics[trim={0cm 0cm 21cm 8cm},clip,width=\linewidth]{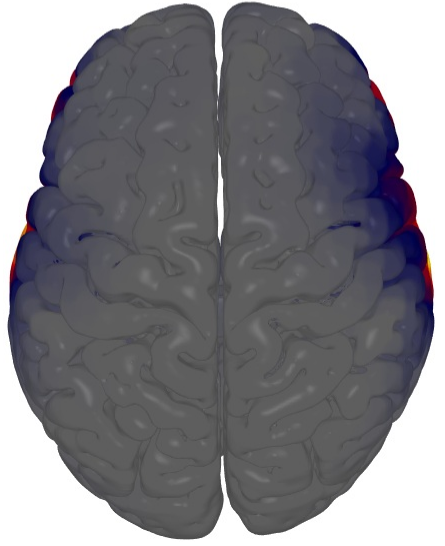}
    \end{minipage}\begin{minipage}{\ImgWidth}
    \includegraphics[trim={0cm 0cm 19cm 9cm},clip,width=\linewidth]{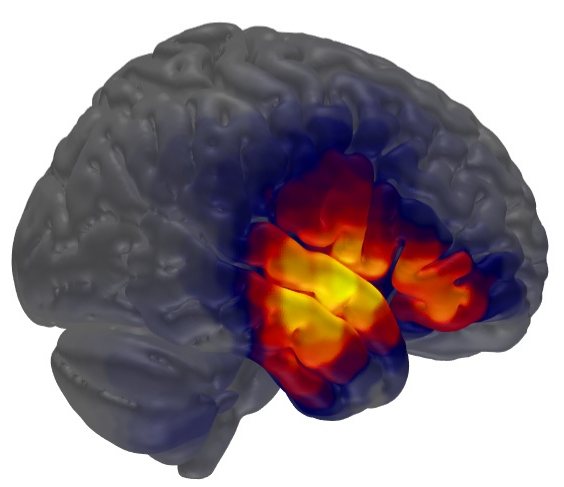}\end{minipage}
    \end{center}
    \end{minipage}\begin{minipage}{\BoxWidth}
    \begin{center}
        {\bf DTI-KF}
    \begin{minipage}{\ImgWidth}
        \includegraphics[trim={0cm 0cm 19cm 9cm},clip,width=\linewidth]{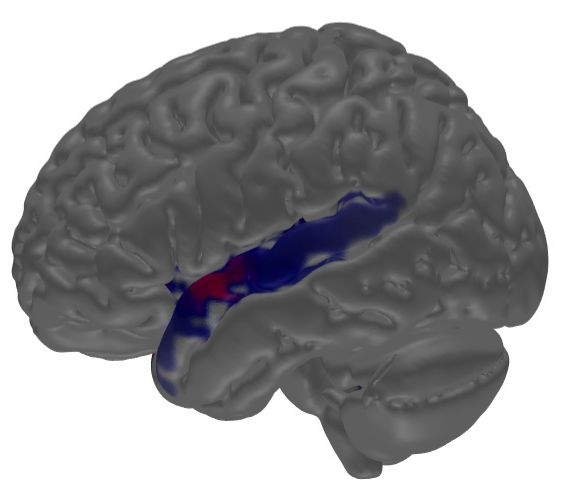}    
    \end{minipage}\begin{minipage}{\MidImgWidth}
    \includegraphics[trim={0cm 0cm 21cm 8cm},clip,width=\linewidth]{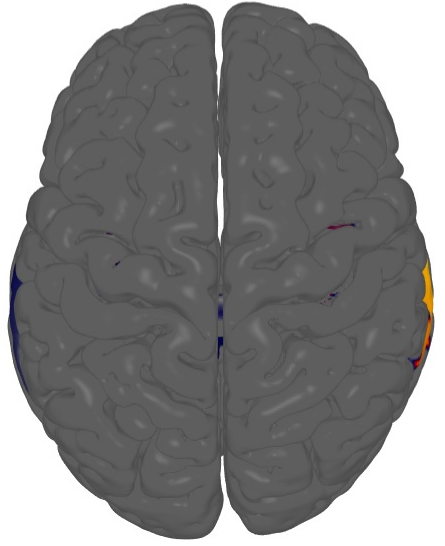}
    \end{minipage}\begin{minipage}{\ImgWidth}
    \includegraphics[trim={0cm 0cm 19cm 9cm},clip,width=\linewidth]{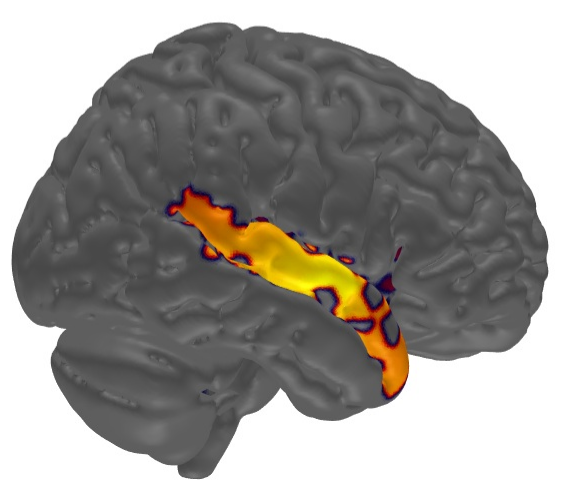}\end{minipage}
    \end{center}
    \end{minipage}\begin{minipage}{\BoxWidth}
    \begin{center}
        {\bf DTI-SKF}
    \begin{minipage}{\ImgWidth}
        \includegraphics[trim={0cm 0cm 19cm 9cm},clip,width=\linewidth]{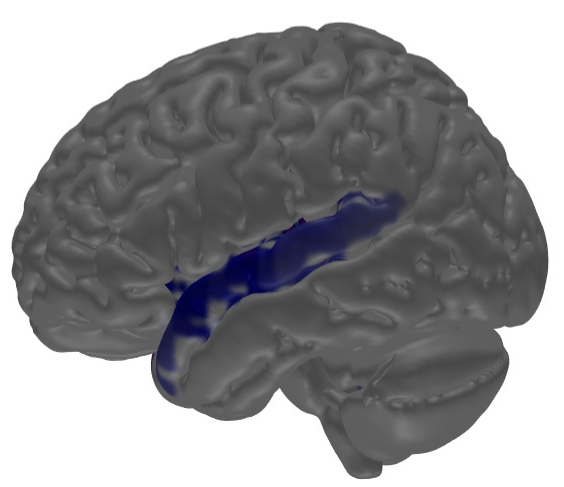}    
    \end{minipage}\begin{minipage}{\MidImgWidth}
    \includegraphics[trim={0cm 0cm 21cm 8cm},clip,width=\linewidth]{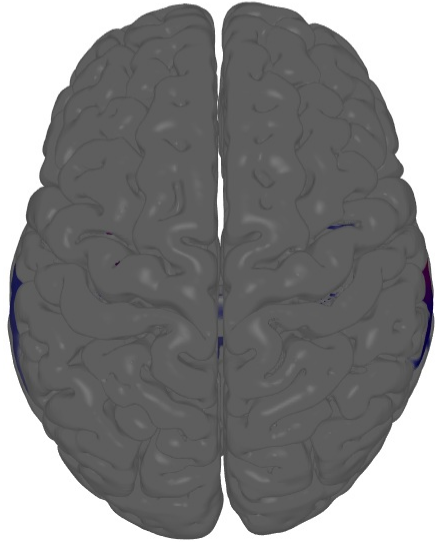}
    \end{minipage}\begin{minipage}{\ImgWidth}
    \includegraphics[trim={0cm 0cm 19cm 9cm},clip,width=\linewidth]{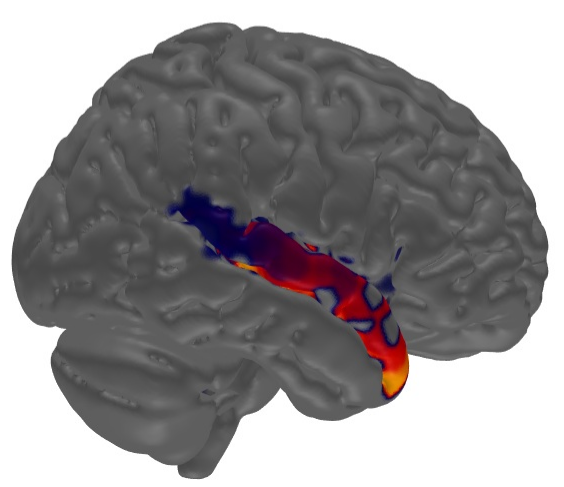}\end{minipage}
    \end{center}
    \end{minipage}

    \par\noindent\rule{\linewidth}{0.8pt}

    \begin{minipage}[b]{\RotTextWidth}
        \rotatebox{90}{\hspace{-0.8cm}\small Insula R}
    \end{minipage}\begin{minipage}{\BoxWidth}
    \begin{center}
        {\bf KF}
    \begin{minipage}{\ImgWidth}
        \includegraphics[trim={0cm 0cm 19cm 9cm},clip,width=\linewidth]{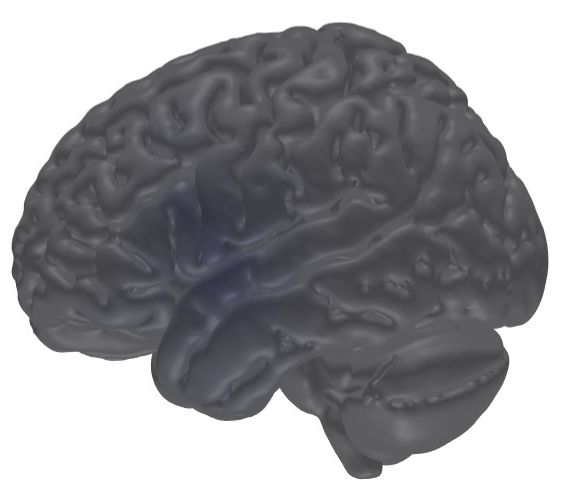}    
    \end{minipage}\begin{minipage}{\MidImgWidth}
    \includegraphics[trim={0cm 0cm 21cm 8cm},clip,width=\linewidth]{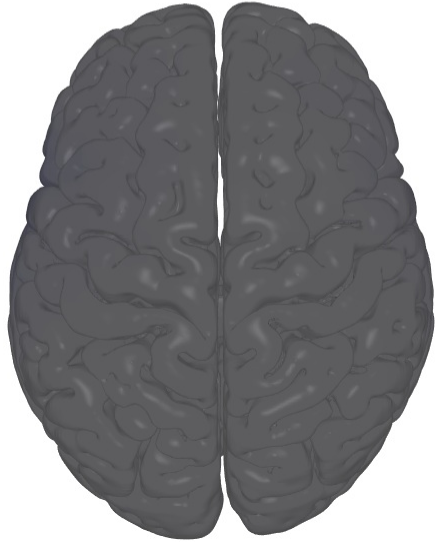}
    \end{minipage}\begin{minipage}{\ImgWidth}
    \includegraphics[trim={0cm 0cm 19cm 9cm},clip,width=\linewidth]{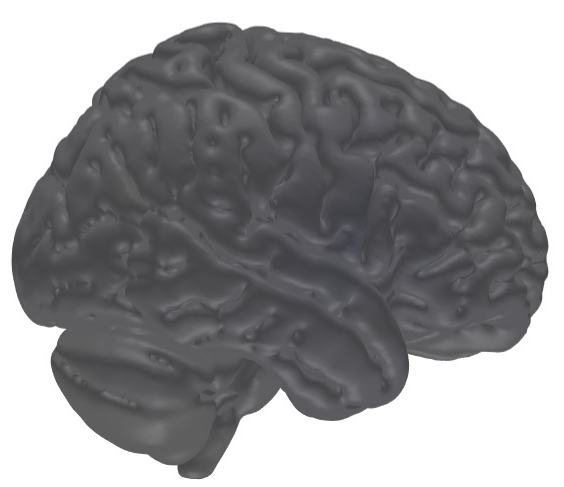}\end{minipage}
    \end{center}
    \end{minipage}\begin{minipage}{\BoxWidth}
    \begin{center}
        {\bf SKF}
    \begin{minipage}{\ImgWidth}
        \includegraphics[trim={0cm 0cm 19cm 9cm},clip,width=\linewidth]{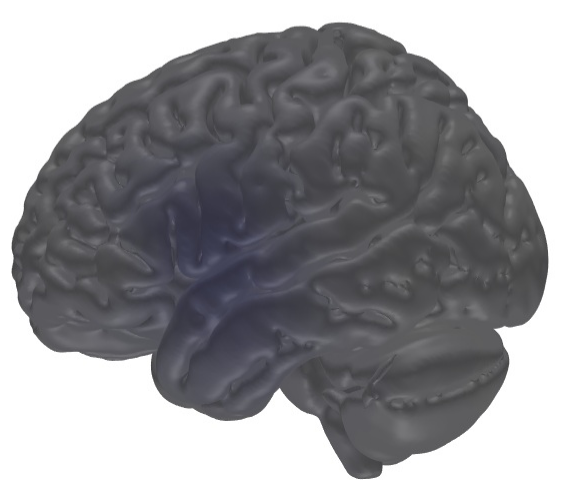}    
    \end{minipage}\begin{minipage}{\MidImgWidth}
    \includegraphics[trim={0cm 0cm 21cm 8cm},clip,width=\linewidth]{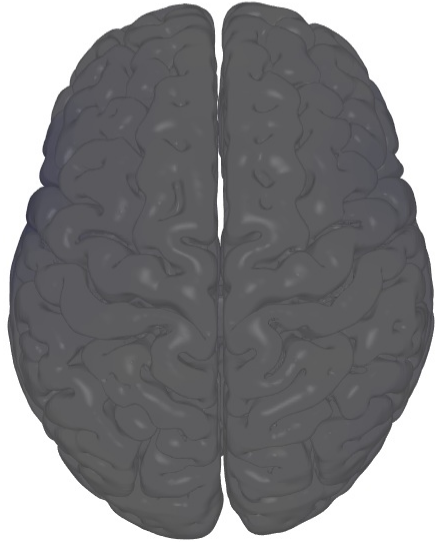}
    \end{minipage}\begin{minipage}{\ImgWidth}
    \includegraphics[trim={0cm 0cm 19cm 9cm},clip,width=\linewidth]{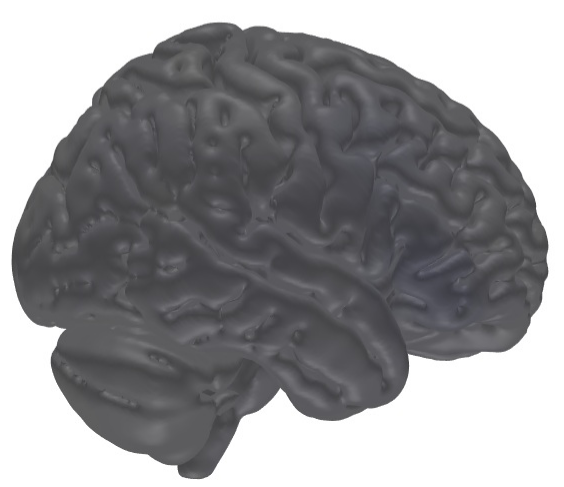}\end{minipage}
    \end{center}
    \end{minipage}\begin{minipage}{\BoxWidth}
    \begin{center}
        {\bf DTI-KF}
    \begin{minipage}{\ImgWidth}
        \includegraphics[trim={0cm 0cm 19cm 9cm},clip,width=\linewidth]{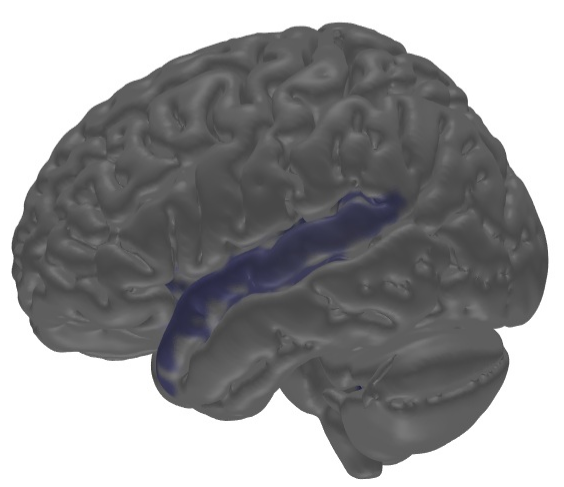}    
    \end{minipage}\begin{minipage}{\MidImgWidth}
    \includegraphics[trim={0cm 0cm 21cm 8cm},clip,width=\linewidth]{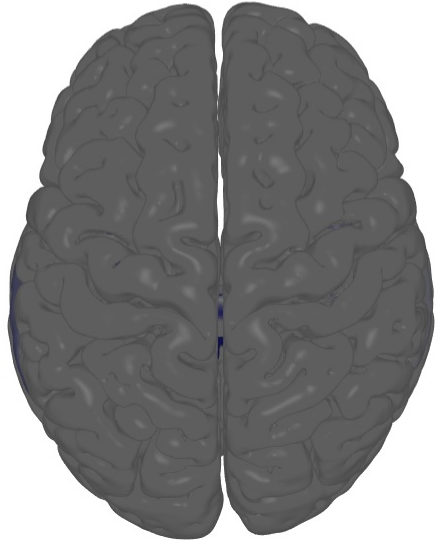}
    \end{minipage}\begin{minipage}{\ImgWidth}
    \includegraphics[trim={0cm 0cm 19cm 9cm},clip,width=\linewidth]{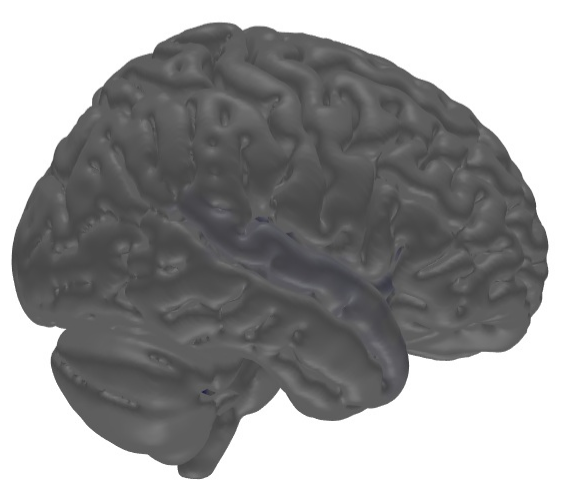}\end{minipage}
    \end{center}
    \end{minipage}\begin{minipage}{\BoxWidth}
    \begin{center}
        {\bf DTI-SKF}
    \begin{minipage}{\ImgWidth}
        \includegraphics[trim={0cm 0cm 19cm 9cm},clip,width=\linewidth]{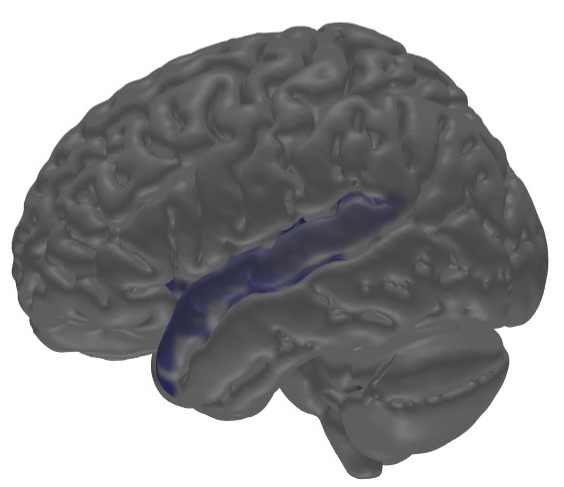}    
    \end{minipage}\begin{minipage}{\MidImgWidth}
    \includegraphics[trim={0cm 0cm 21cm 8cm},clip,width=\linewidth]{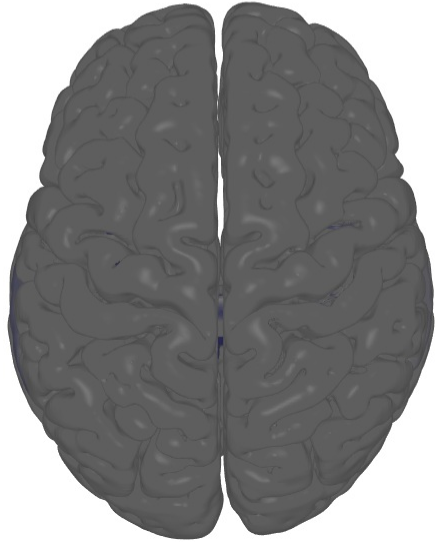}
    \end{minipage}\begin{minipage}{\ImgWidth}
    \includegraphics[trim={0cm 0cm 19cm 9cm},clip,width=\linewidth]{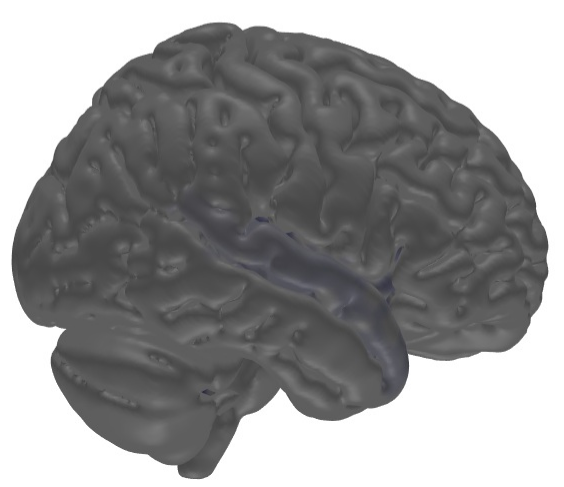}\end{minipage}
    \end{center}
    \end{minipage}

    \begin{minipage}[b]{\RotTextWidth}
        \rotatebox{90}{\hspace{-0.5cm}\small S.T R}
    \end{minipage}\begin{minipage}{\BoxWidth}
    \begin{center}
        {\bf KF}
    \begin{minipage}{\ImgWidth}
        \includegraphics[trim={0cm 0cm 19cm 9cm},clip,width=\linewidth]{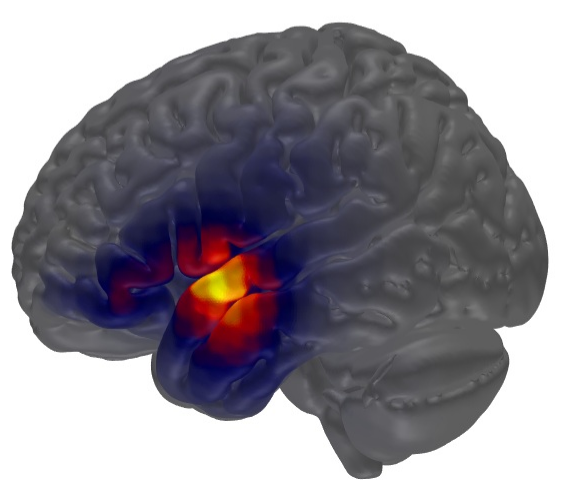}    
    \end{minipage}\begin{minipage}{\MidImgWidth}
    \includegraphics[trim={0cm 0cm 21cm 8cm},clip,width=\linewidth]{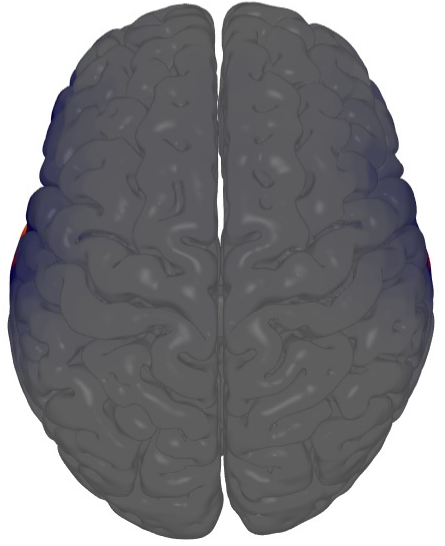}
    \end{minipage}\begin{minipage}{\ImgWidth}
    \includegraphics[trim={0cm 0cm 19cm 9cm},clip,width=\linewidth]{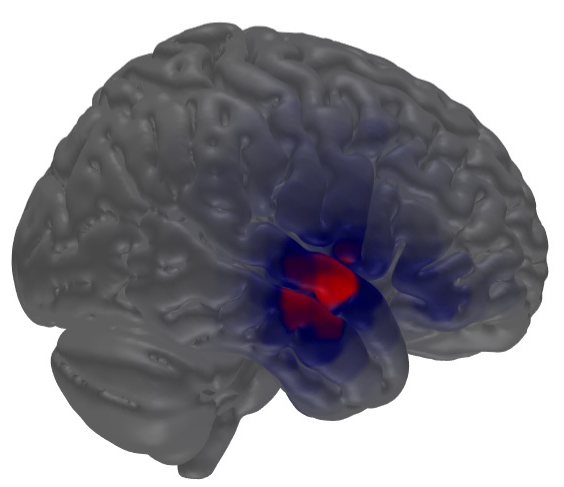}\end{minipage}
    \end{center}
    \end{minipage}\begin{minipage}{\BoxWidth}
    \begin{center}
        {\bf SKF}
    \begin{minipage}{\ImgWidth}
        \includegraphics[trim={0cm 0cm 19cm 9cm},clip,width=\linewidth]{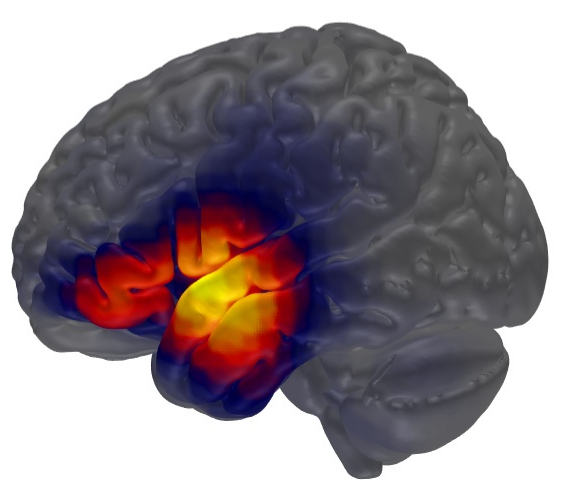}    
    \end{minipage}\begin{minipage}{\MidImgWidth}
    \includegraphics[trim={0cm 0cm 21cm 8cm},clip,width=\linewidth]{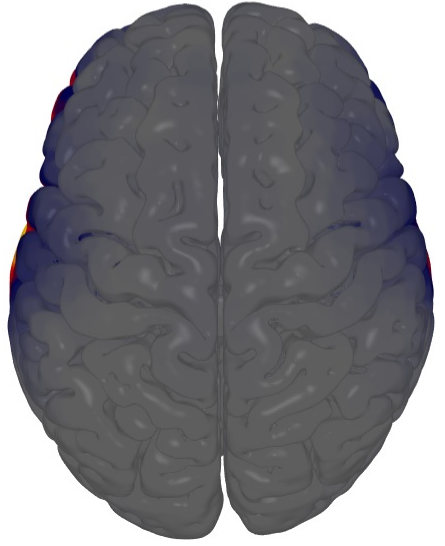}
    \end{minipage}\begin{minipage}{\ImgWidth}
    \includegraphics[trim={0cm 0cm 19cm 9cm},clip,width=\linewidth]{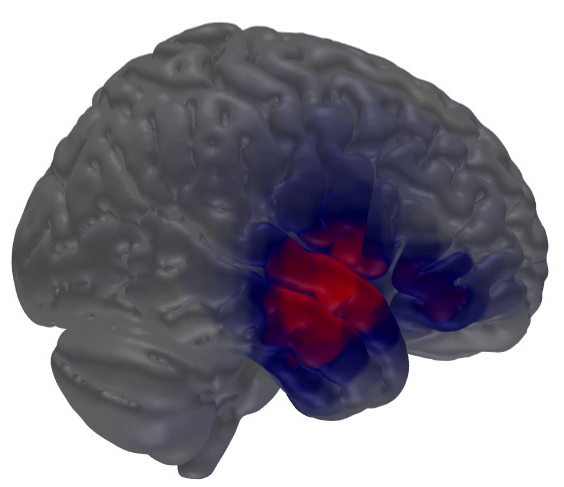}\end{minipage}
    \end{center}
    \end{minipage}\begin{minipage}{\BoxWidth}
    \begin{center}
        {\bf DTI-KF}
    \begin{minipage}{\ImgWidth}
        \includegraphics[trim={0cm 0cm 19cm 9cm},clip,width=\linewidth]{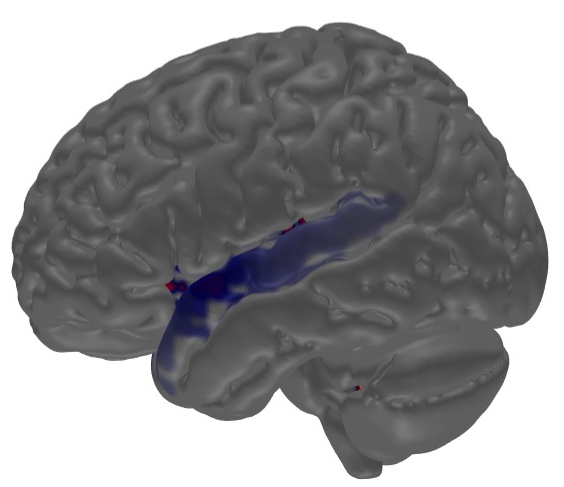}    
    \end{minipage}\begin{minipage}{\MidImgWidth}
    \includegraphics[trim={0cm 0cm 21cm 8cm},clip,width=\linewidth]{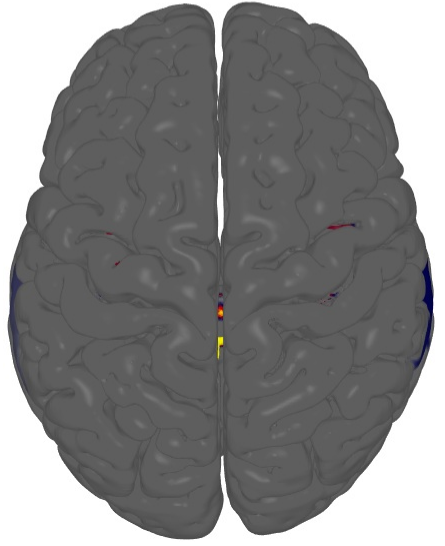}
    \end{minipage}\begin{minipage}{\ImgWidth}
    \includegraphics[trim={0cm 0cm 19cm 9cm},clip,width=\linewidth]{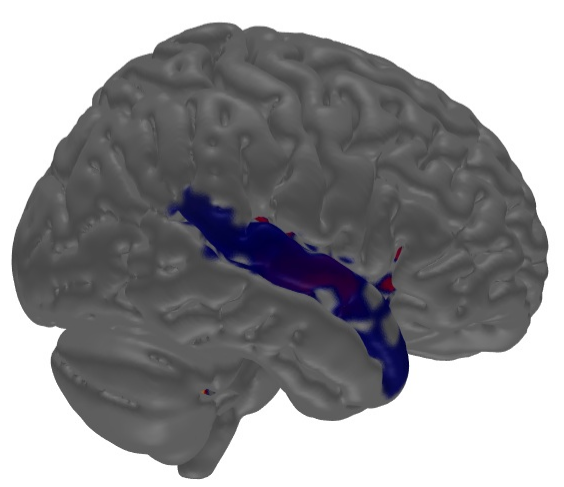}\end{minipage}
    \end{center}
    \end{minipage}\begin{minipage}{\BoxWidth}
    \begin{center}
        {\bf DTI-SKF}
    \begin{minipage}{\ImgWidth}
        \includegraphics[trim={0cm 0cm 19cm 9cm},clip,width=\linewidth]{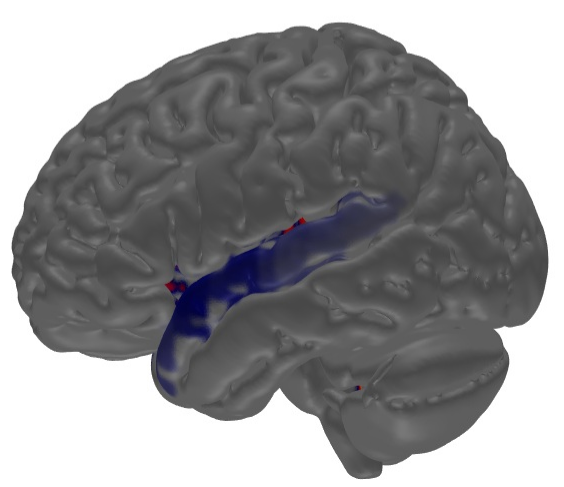}    
    \end{minipage}\begin{minipage}{\MidImgWidth}
    \includegraphics[trim={0cm 0cm 21cm 8cm},clip,width=\linewidth]{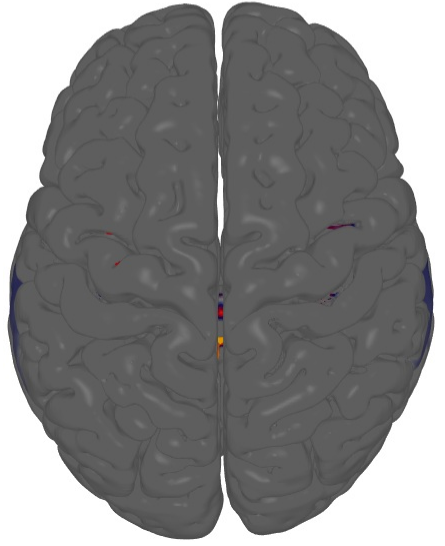}
    \end{minipage}\begin{minipage}{\ImgWidth}
    \includegraphics[trim={0cm 0cm 19cm 9cm},clip,width=\linewidth]{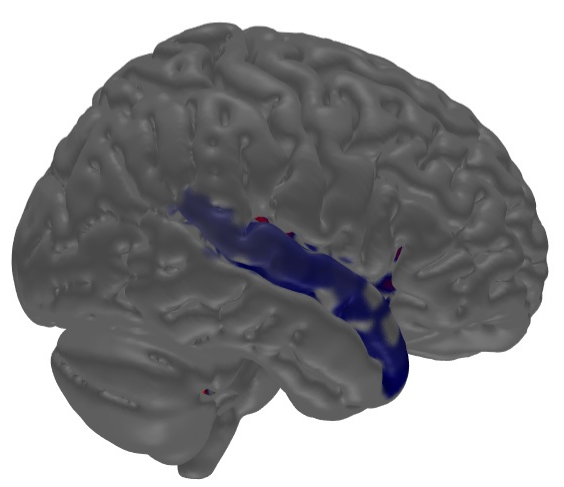}\end{minipage}
    \end{center}
    \end{minipage}

    \begin{minipage}[b]{\RotTextWidth}
        \rotatebox{90}{\hspace{-0.5cm}\small T.T. R}
    \end{minipage}\begin{minipage}{\BoxWidth}
    \begin{center}
        {\bf KF}
    \begin{minipage}{\ImgWidth}
        \includegraphics[trim={0cm 0cm 19cm 9cm},clip,width=\linewidth]{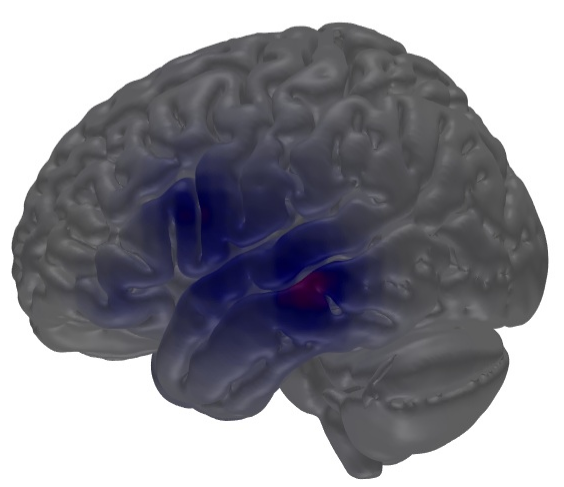}    
    \end{minipage}\begin{minipage}{\MidImgWidth}
    \includegraphics[trim={0cm 0cm 21cm 8cm},clip,width=\linewidth]{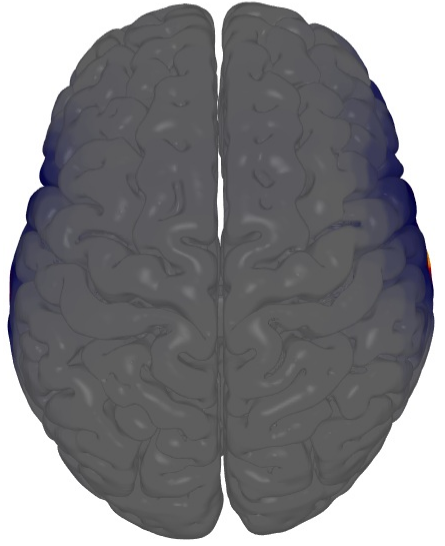}
    \end{minipage}\begin{minipage}{\ImgWidth}
    \includegraphics[trim={0cm 0cm 19cm 9cm},clip,width=\linewidth]{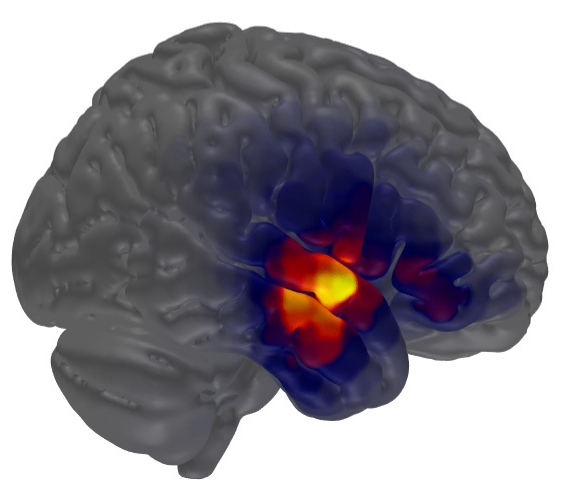}\end{minipage}
    \end{center}
    \end{minipage}\begin{minipage}{\BoxWidth}
    \begin{center}
        {\bf SKF}
    \begin{minipage}{\ImgWidth}
        \includegraphics[trim={0cm 0cm 19cm 9cm},clip,width=\linewidth]{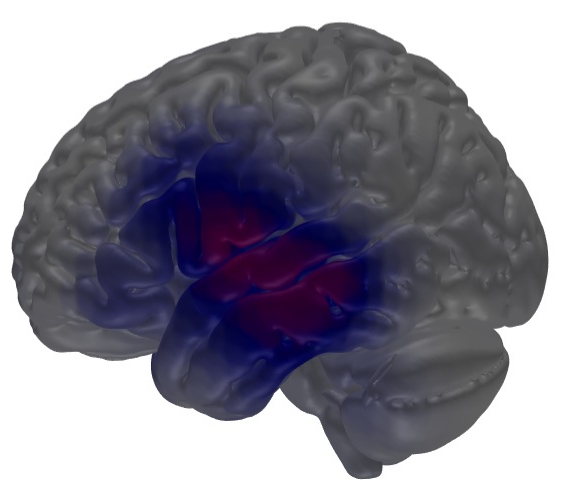}    
    \end{minipage}\begin{minipage}{\MidImgWidth}
    \includegraphics[trim={0cm 0cm 21cm 8cm},clip,width=\linewidth]{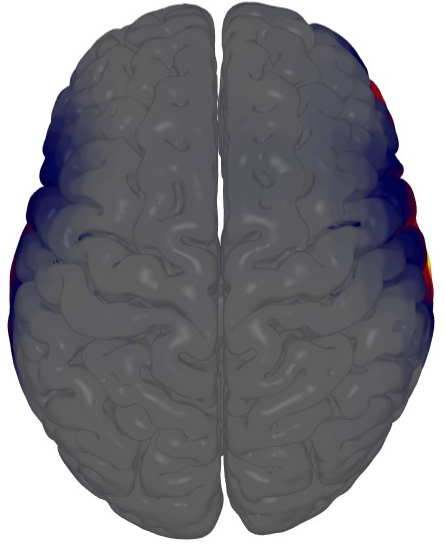}
    \end{minipage}\begin{minipage}{\ImgWidth}
    \includegraphics[trim={0cm 0cm 19cm 9cm},clip,width=\linewidth]{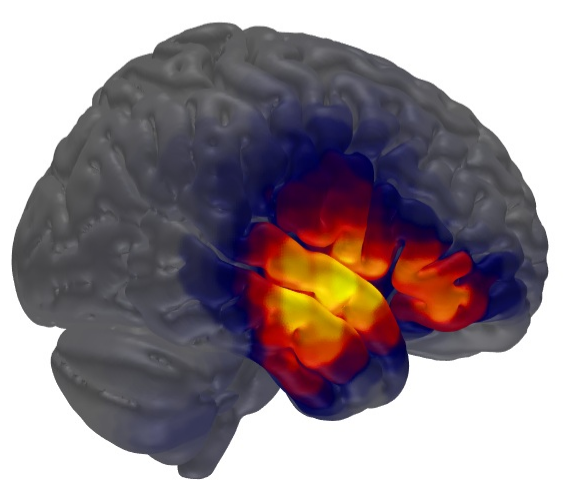}\end{minipage}
    \end{center}
    \end{minipage}\begin{minipage}{\BoxWidth}
    \begin{center}
        {\bf DTI-KF}
    \begin{minipage}{\ImgWidth}
        \includegraphics[trim={0cm 0cm 19cm 9cm},clip,width=\linewidth]{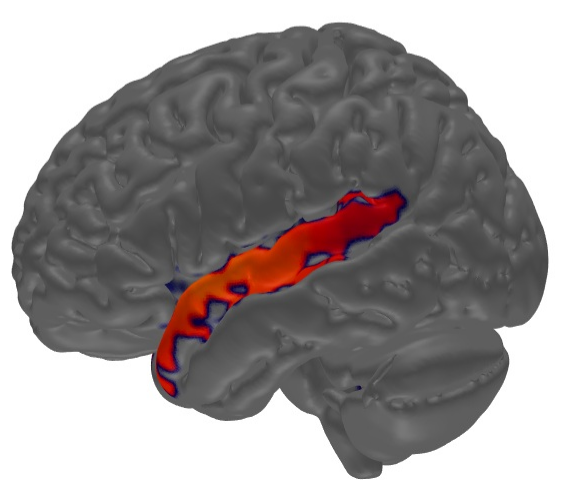}    
    \end{minipage}\begin{minipage}{\MidImgWidth}
    \includegraphics[trim={0cm 0cm 21cm 8cm},clip,width=\linewidth]{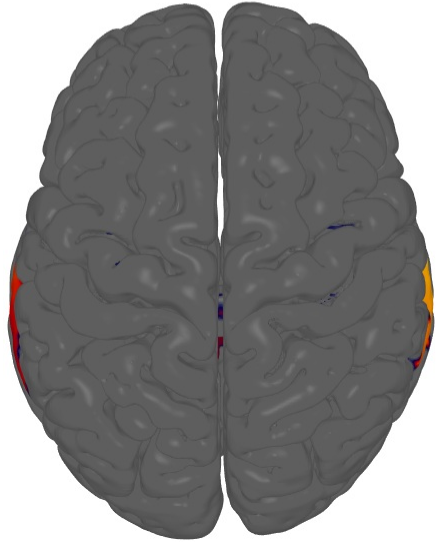}
    \end{minipage}\begin{minipage}{\ImgWidth}
    \includegraphics[trim={0cm 0cm 19cm 9cm},clip,width=\linewidth]{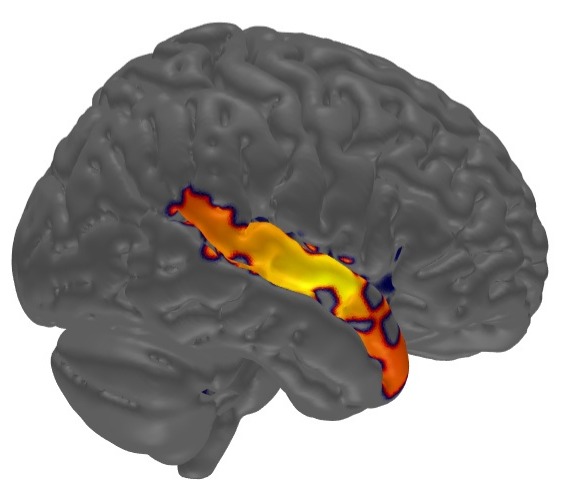}\end{minipage}
    \end{center}
    \end{minipage}\begin{minipage}{\BoxWidth}
    \begin{center}
        {\bf DTI-SKF}
    \begin{minipage}{\ImgWidth}
        \includegraphics[trim={0cm 0cm 19cm 9cm},clip,width=\linewidth]{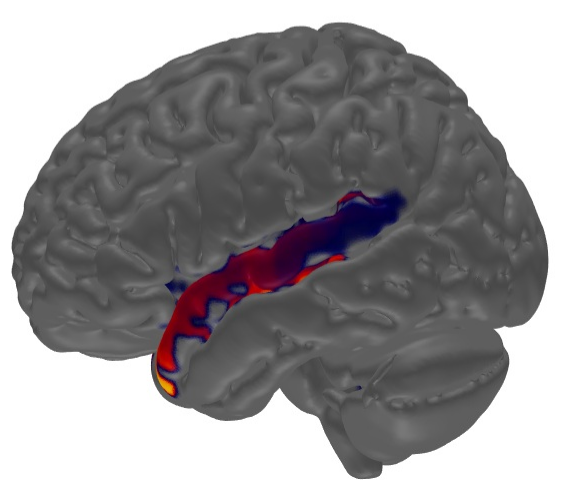}    
    \end{minipage}\begin{minipage}{\MidImgWidth}
    \includegraphics[trim={0cm 0cm 21cm 8cm},clip,width=\linewidth]{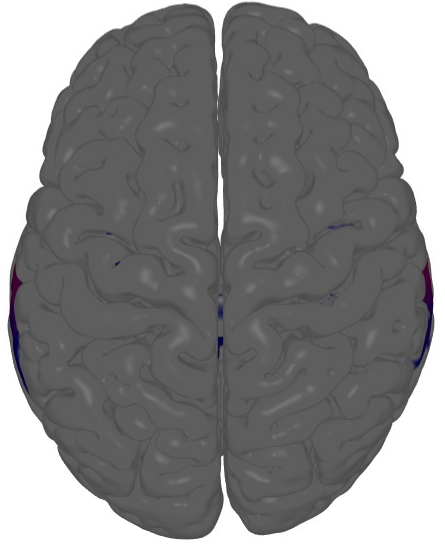}
    \end{minipage}\begin{minipage}{\ImgWidth}
    \includegraphics[trim={0cm 0cm 19cm 9cm},clip,width=\linewidth]{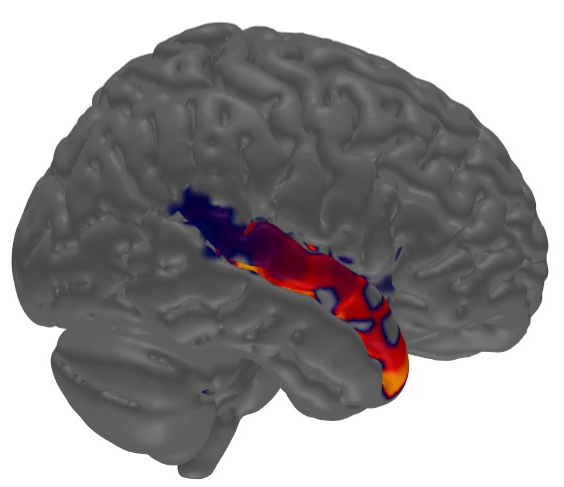}\end{minipage}
    \end{center}
    \end{minipage}
    
    \caption{Source estimation projections on the brain model for four events appearing in the AEP simulation for each anatomical region and method; the brain is shown from the left, top, and right, respectively. The upper block shows the result when the brain region on the left is active, indicated by L after the anatomical region abbreviation. In the lower block, R after the region abbreviation denotes the right hemisphere. S.T. means superior temporal gyrus; T.T on the other hand, is transverse temporal gyrus. 
    The first two rows present  the Kalman filter and the Standardized Kalman filter using a random walk evolution model, and the latter two show the corresponding methods that use a DTI-based transition model.}
    \label{fig:AUDITORYbrain}
\end{figure*}

In Figure \ref{fig:AUDITORYbrain}, at the beginning, we see faint activity in the left insula across all methods. The activity progresses more slowly than it should, causing the left insula to show more pronounced activity once the right insula peaks. After that, we see activity peaking correctly in the left, except that SKF slightly mislocalizes the activity. When the activity is supposed to move right, the peak stays left with random walk evolution-modeled methods. The correct fading left and peaking right happens with DTI-evolution-modeled methods. In the third phase, when activity should be located in the left transverse temporal gyrus, all methods show stronger activity in the right hemisphere. All methods detect the last component on the correct side; however, DTI-KF shows the most accurate localization.

\begin{figure*}
    \centering
    \begin{minipage}{0.49\linewidth}
    \centering
    {\bf KF}
        \begin{minipage}{0.5\linewidth}
        \centering
        \small{left}
        \includegraphics[width=0.98\linewidth]{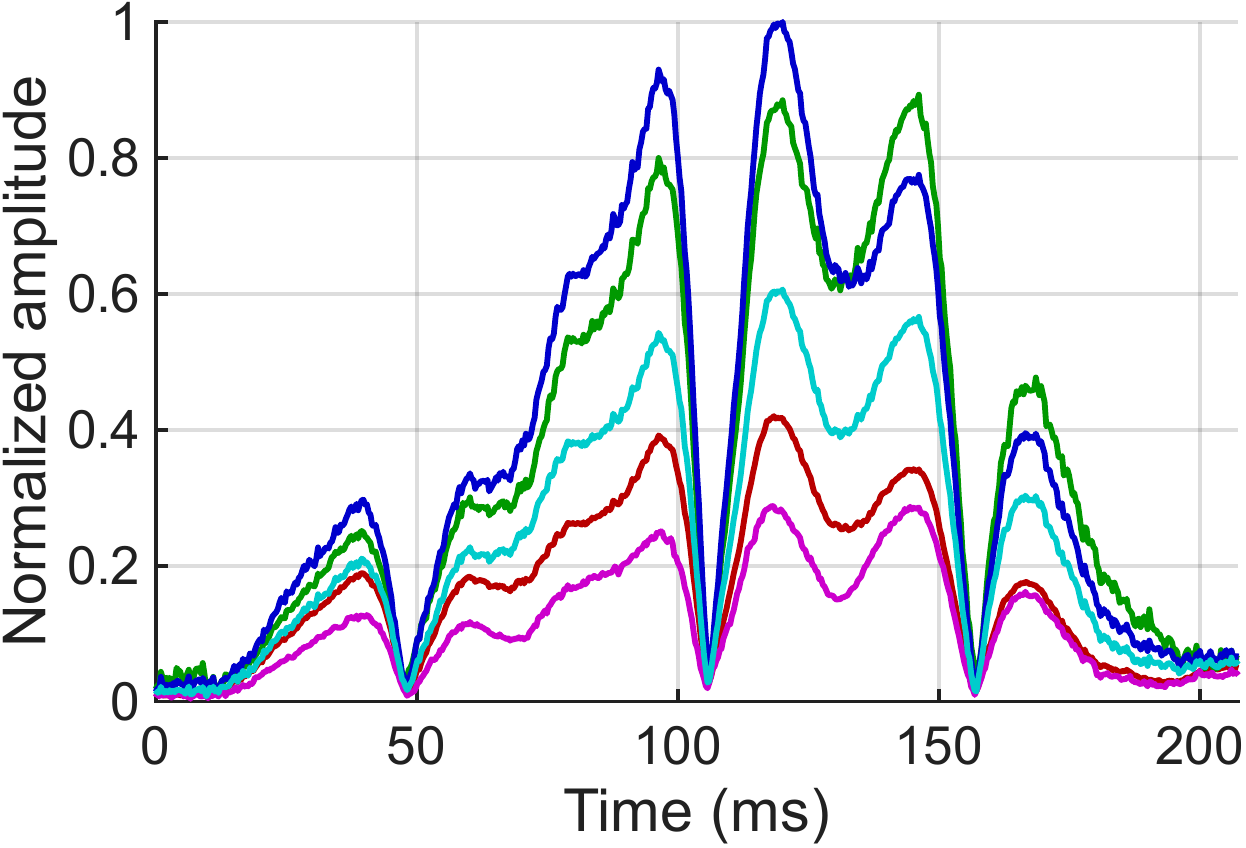}
    \end{minipage}\begin{minipage}{0.5\linewidth}
    \centering
        \small{right}
        \includegraphics[width=0.98\linewidth]{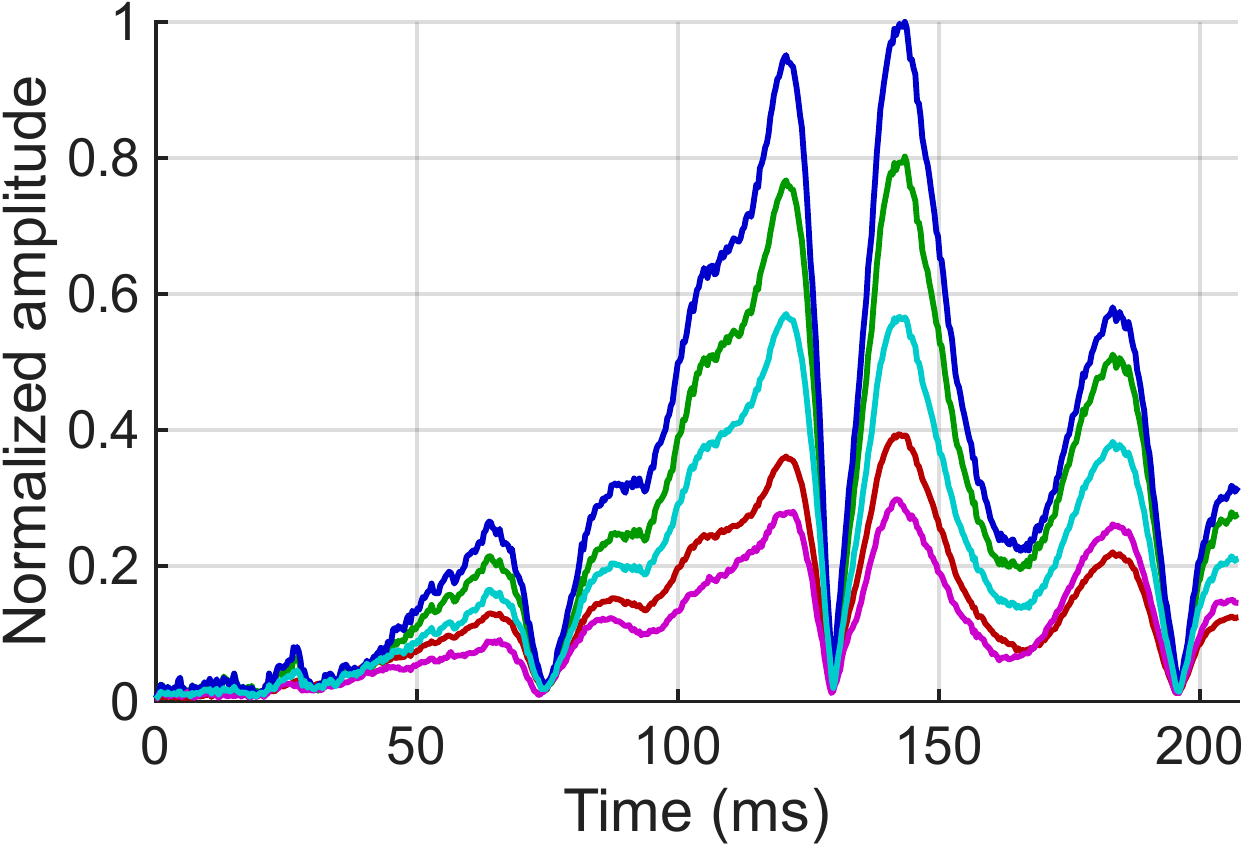}
    \end{minipage}\begin{minipage}{0.49\linewidth}
        
    \end{minipage}
    \end{minipage}\hspace{0.2cm}\begin{minipage}{0.49\linewidth}
    \centering
    {\bf SKF}
        \begin{minipage}{0.5\linewidth}
        \centering
        \small{left}
        \includegraphics[width=0.98\linewidth]{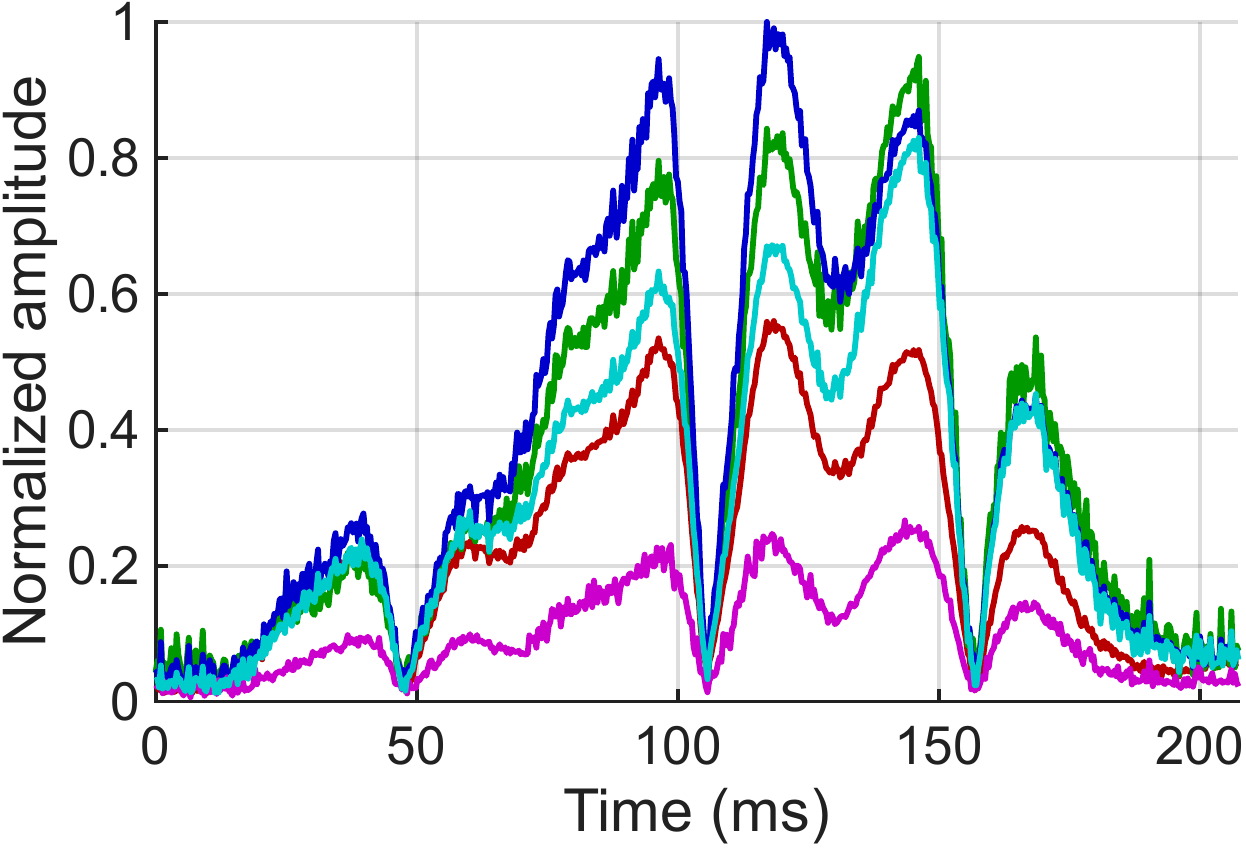}
    \end{minipage}\begin{minipage}{0.5\linewidth}
    \centering
        \small{right}
        \includegraphics[width=0.98\linewidth]{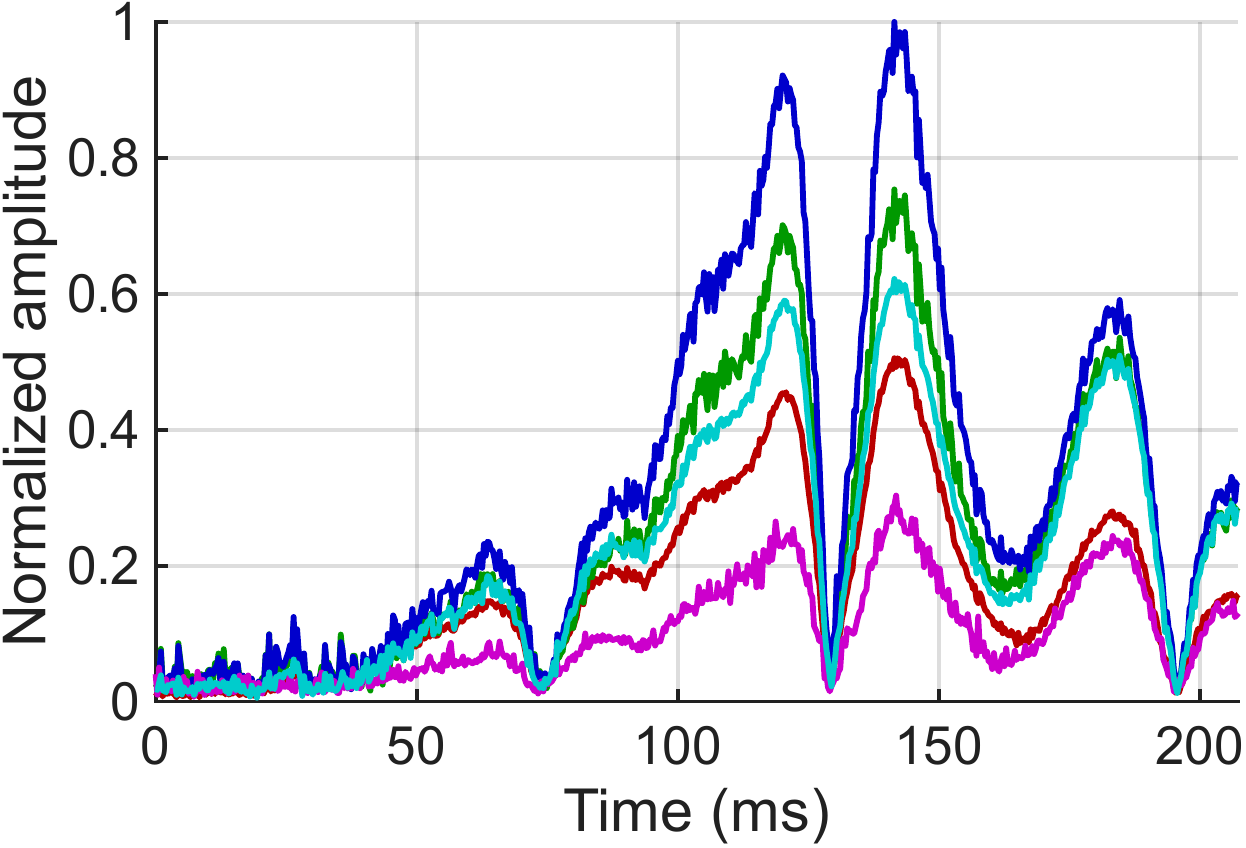}
    \end{minipage}\begin{minipage}{0.49\linewidth}
        
    \end{minipage}
    \end{minipage}

    \begin{minipage}{0.49\linewidth}
    \centering
    {\bf DTI-KF}
        \begin{minipage}{0.5\linewidth}
        \centering
        \small{left}
        \includegraphics[width=0.98\linewidth]{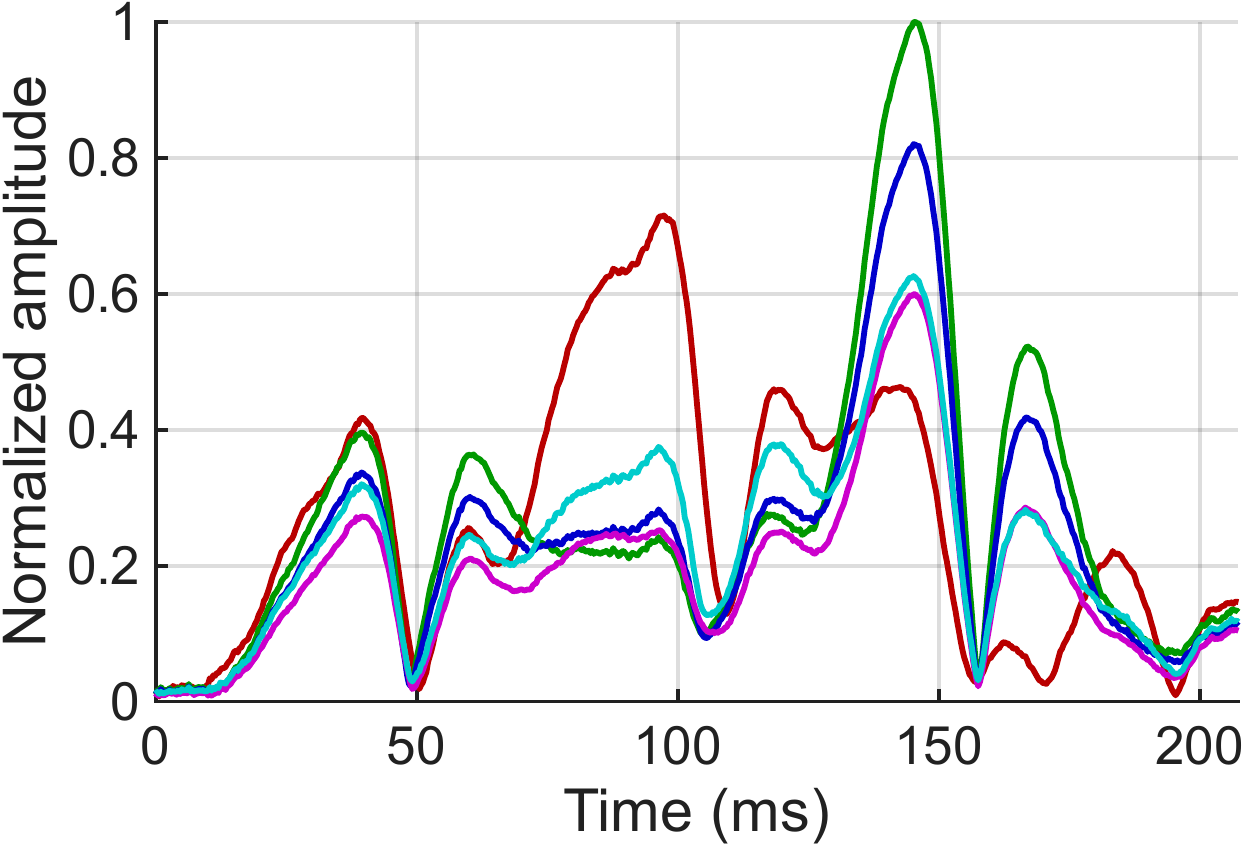}
    \end{minipage}\begin{minipage}{0.5\linewidth}
    \centering
        \small{right}
        \includegraphics[width=0.98\linewidth]{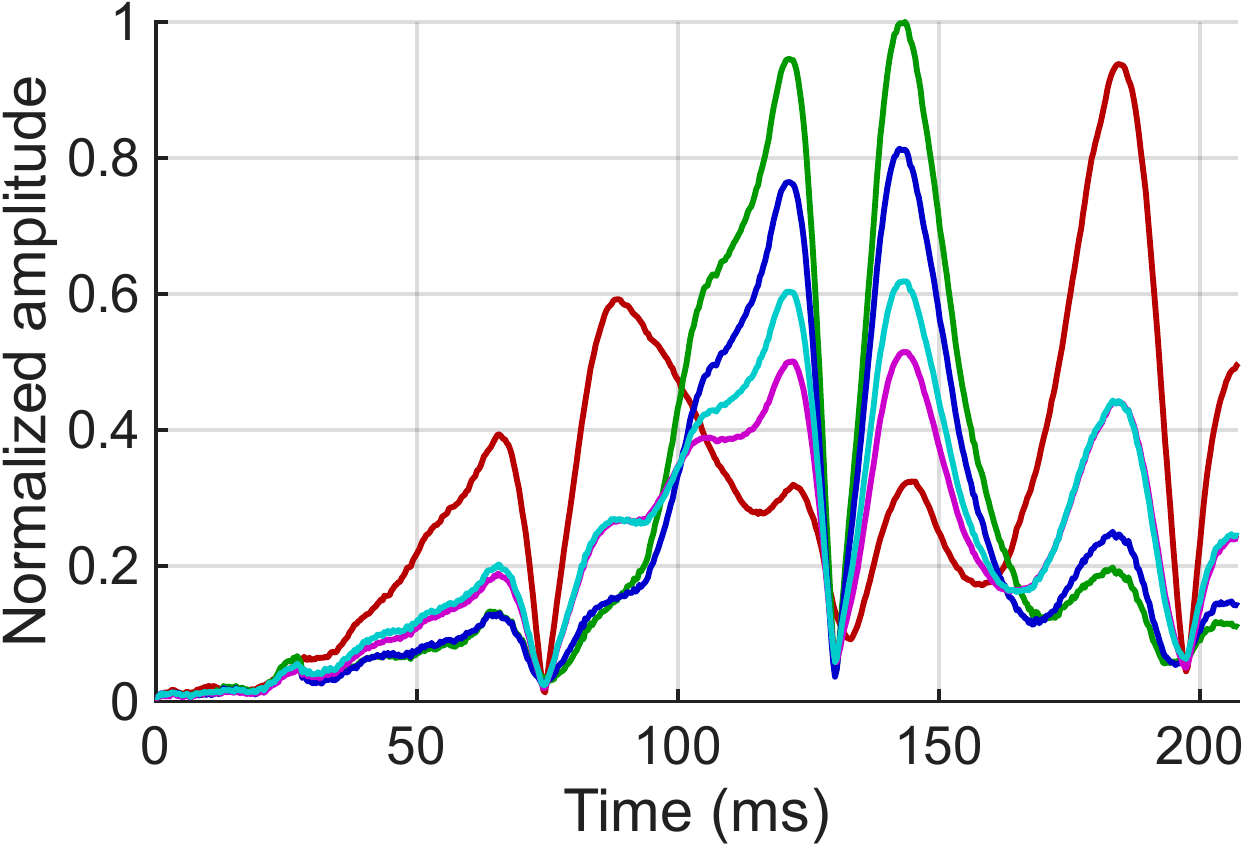}
    \end{minipage}\begin{minipage}{0.49\linewidth}
        
    \end{minipage}
    \end{minipage}\hspace{0.2cm}\begin{minipage}{0.49\linewidth}
    \centering
    {\bf DTI-SKF}
        \begin{minipage}{0.5\linewidth}
        \centering
        \small{left}
        \includegraphics[width=0.98\linewidth]{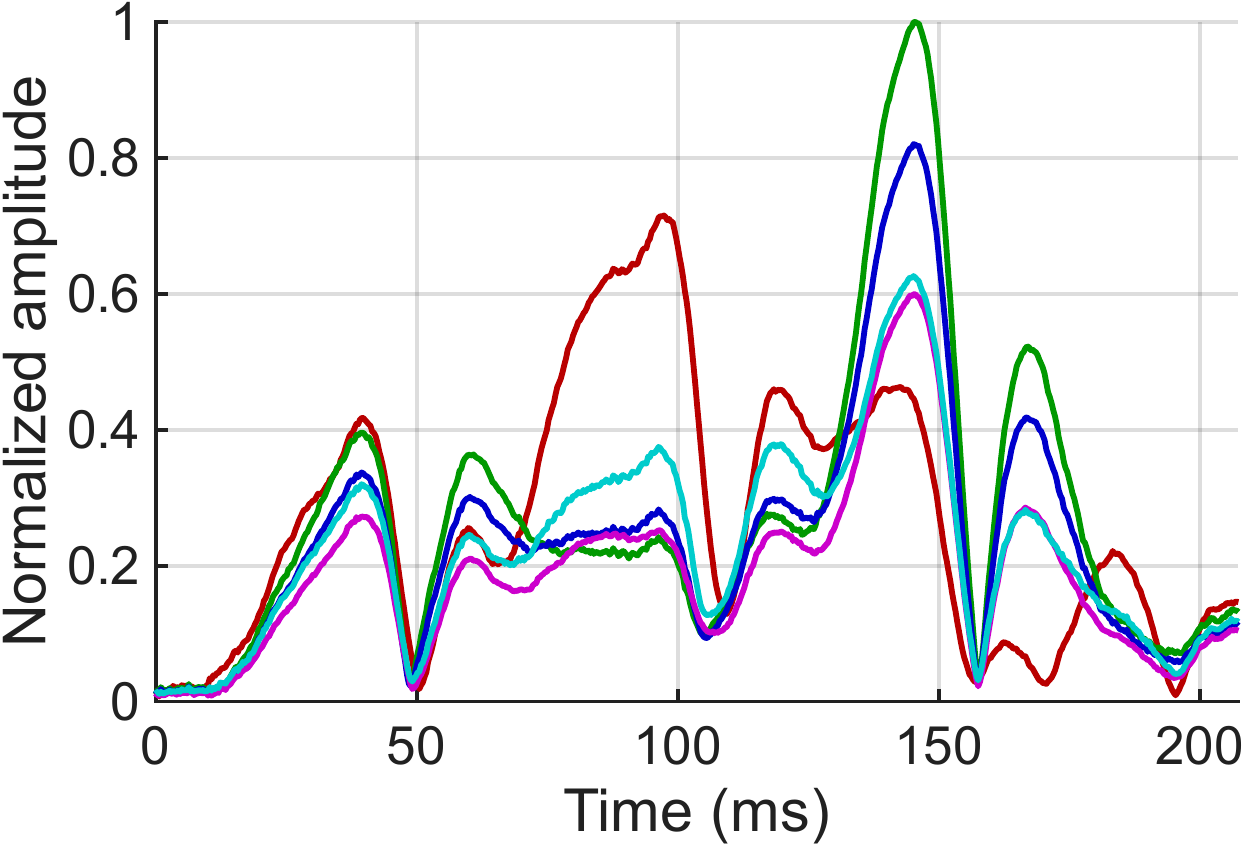}
    \end{minipage}\begin{minipage}{0.5\linewidth}
    \centering
        \small{right}
        \includegraphics[width=0.98\linewidth]{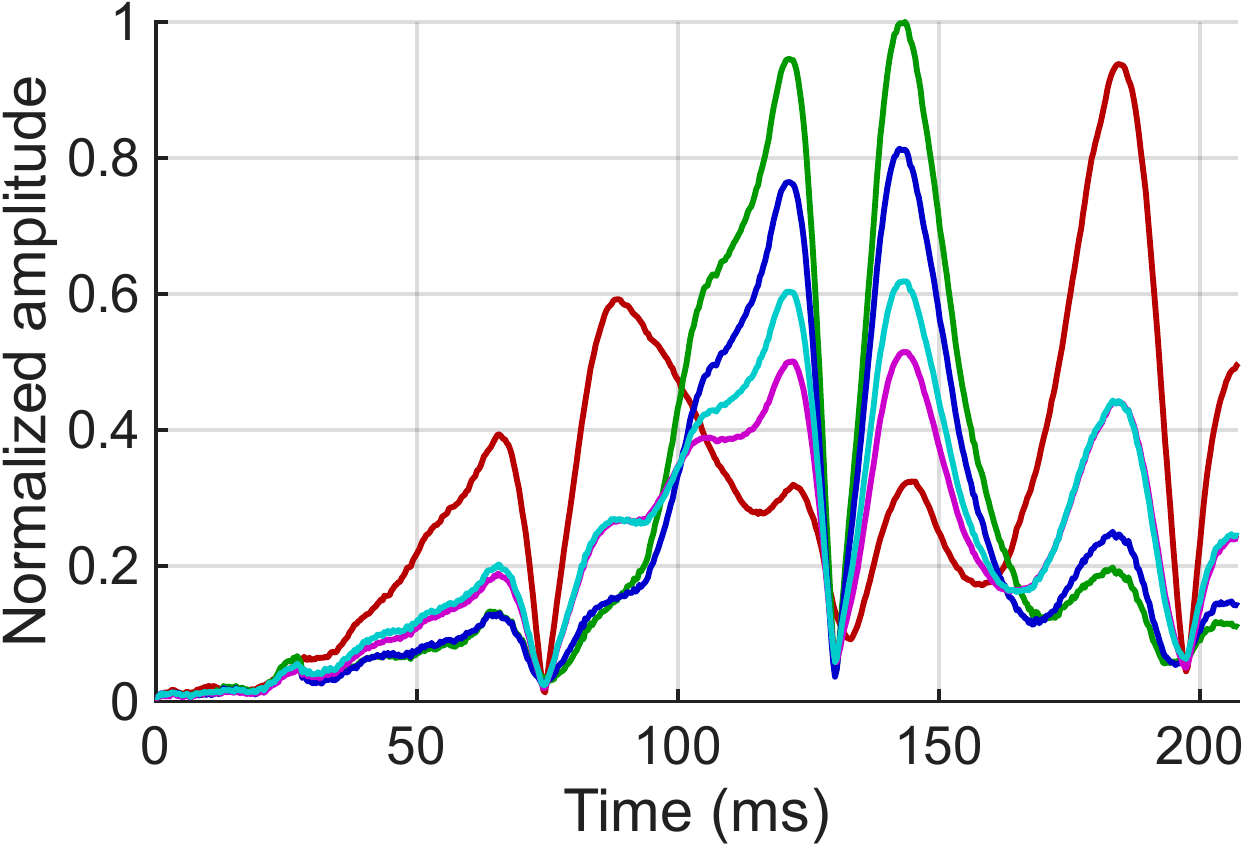}
    \end{minipage}\begin{minipage}{0.49\linewidth}
        
    \end{minipage}
    \end{minipage}

\begin{minipage}{0.49\linewidth}
    \centering
    {\bf Ground truth}
        \begin{minipage}{0.5\linewidth}
        \centering
        \small{left}
        \includegraphics[width=0.98\linewidth]{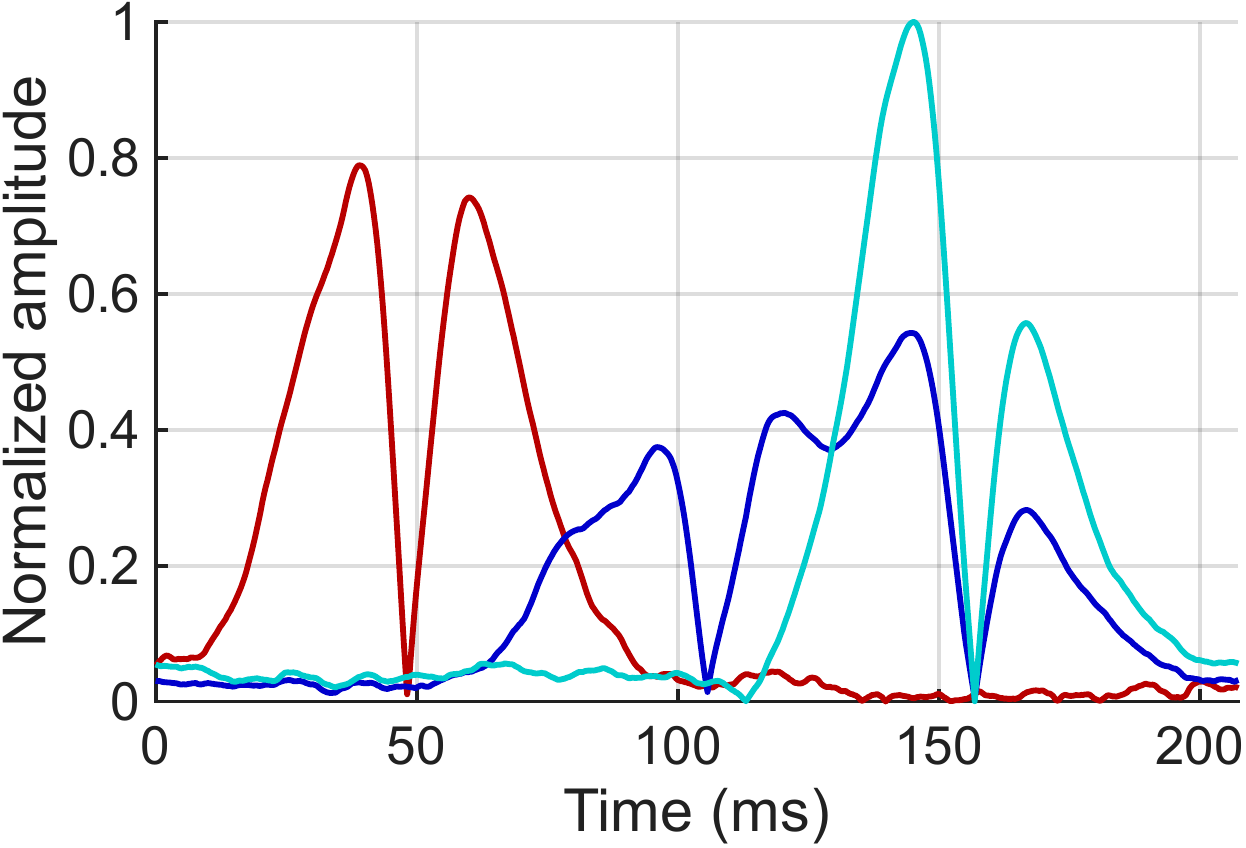}
    \end{minipage}\begin{minipage}{0.5\linewidth}
    \centering
        \small{right}
        \includegraphics[width=0.98\linewidth]{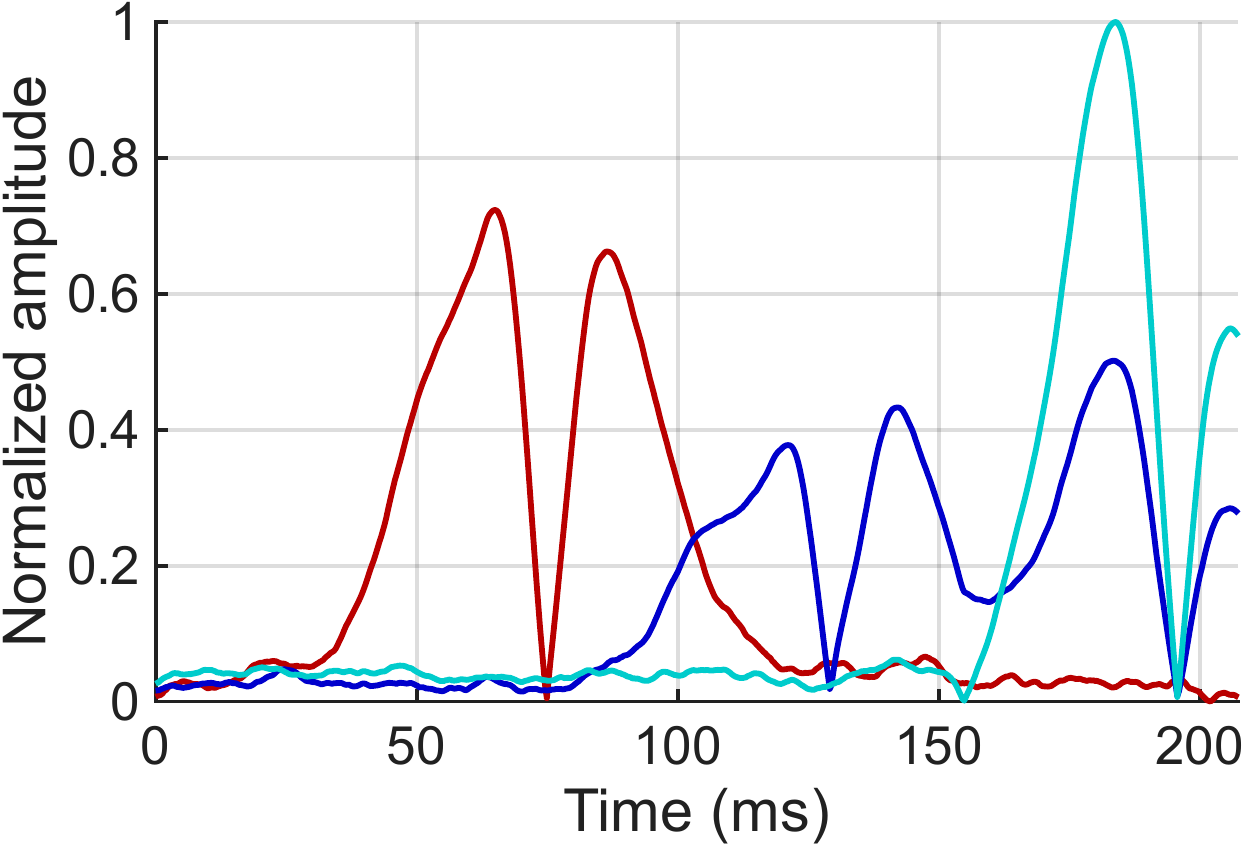}
    \end{minipage}\begin{minipage}{0.49\linewidth}
        
    \end{minipage}
    \end{minipage}
    
    \caption{Normalized estimation magnitude time series shown for the left and right hemispheres separately. The bottom pair shows the magnitude of the true time series. Different anatomical regions are depicted by different colors: insula (dark red), middle temporal gyrus (dark green), superior temporal gyrus (blue), superior marginal gyrus (pink), transverse temporal (teal). }
    \label{fig:AUDITORYtimeseries}
\end{figure*}

\begin{figure*}
\def\imageWidth{0.235\linewidth}
\def\inverseMethodNameWidth{0.02\linewidth}
\def\sideWidth{0.08\linewidth}
    \centering
    \begin{minipage}{\inverseMethodNameWidth}
        \rotatebox{90}{\bf KF}
    \end{minipage}\begin{minipage}{\imageWidth}
        \centering
        \begin{minipage}{\sideWidth}
            \rotatebox{90}{left}
        \end{minipage}\begin{minipage}{\linewidth}
        \centering
        \small{Insula}
            \includegraphics[width=\linewidth]{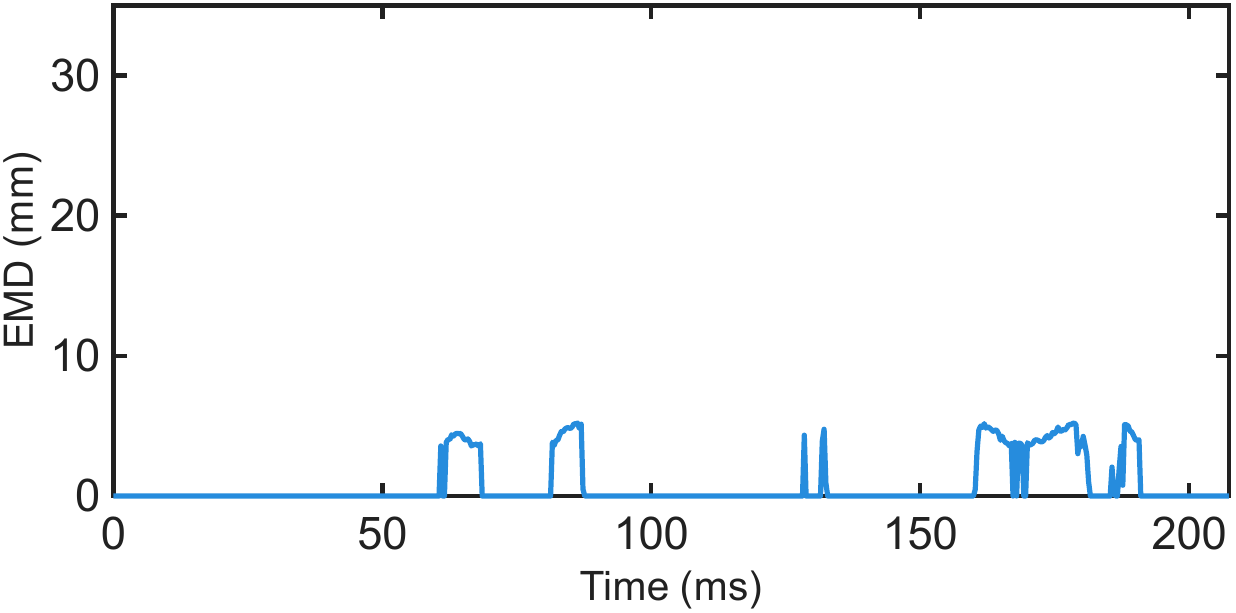}
        \end{minipage}

        \begin{minipage}{\sideWidth}
            \rotatebox{90}{right}
        \end{minipage}\begin{minipage}{\linewidth}
            \includegraphics[width=\linewidth]{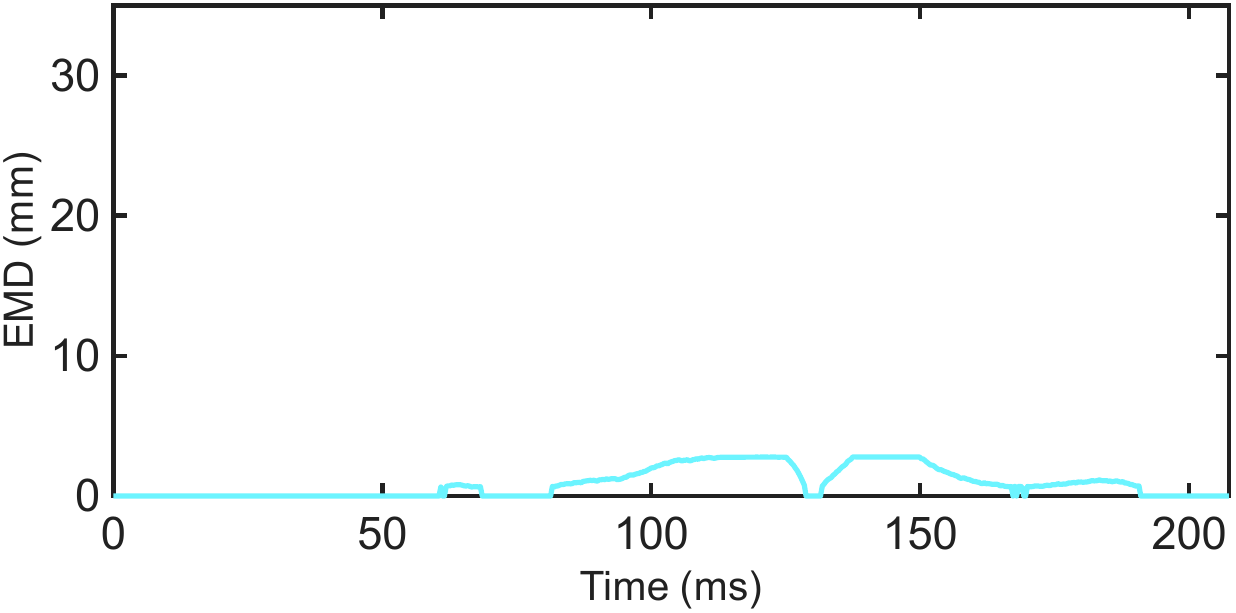}
        \end{minipage}
        
    \end{minipage}\hspace{0.5cm}\begin{minipage}{\imageWidth}
        \centering
        \small{Superior temporal}
        \begin{minipage}{\linewidth}
        \centering
        
            \includegraphics[width=\linewidth]{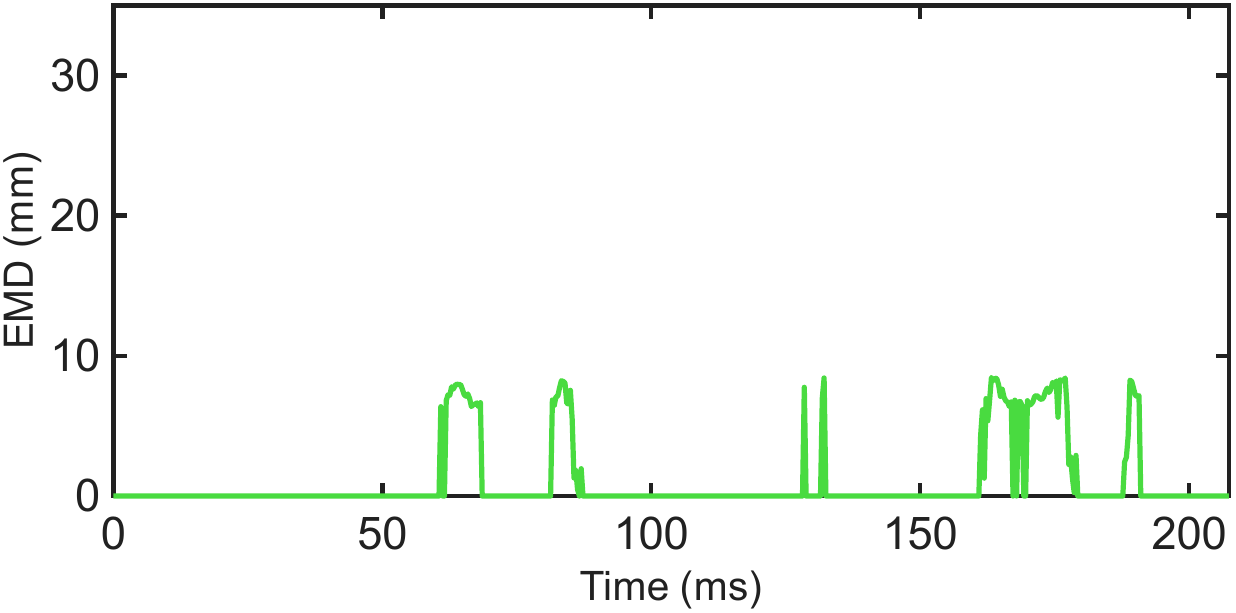}
        \end{minipage}

        \begin{minipage}{\linewidth}
            \includegraphics[width=\linewidth]{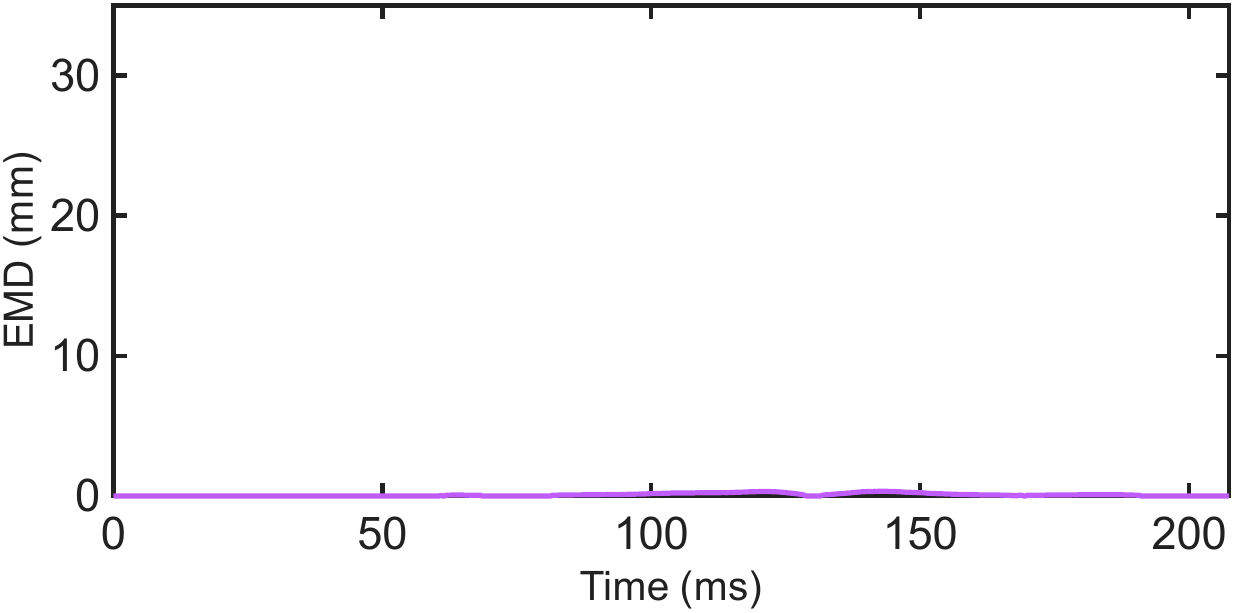}
        \end{minipage}
        
    \end{minipage}\begin{minipage}{\imageWidth}
        \centering
        \begin{minipage}{\linewidth}
        \centering
        \small{Transverse temporal}
            \includegraphics[width=\linewidth]{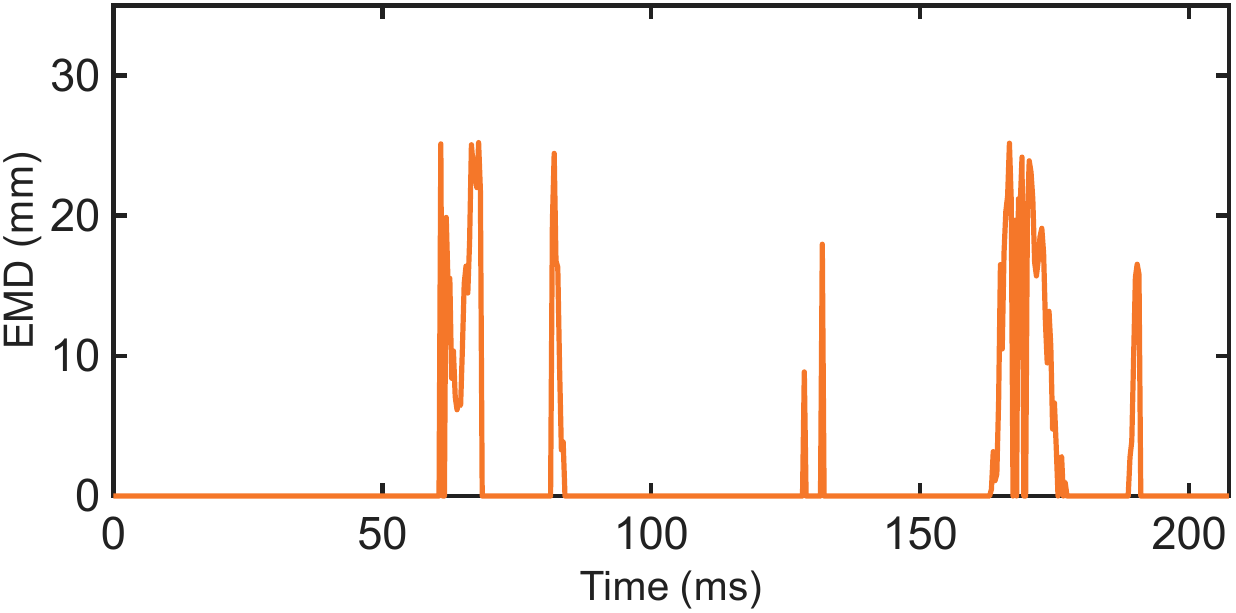}
        \end{minipage}

        \begin{minipage}{\linewidth}
            \includegraphics[width=\linewidth]{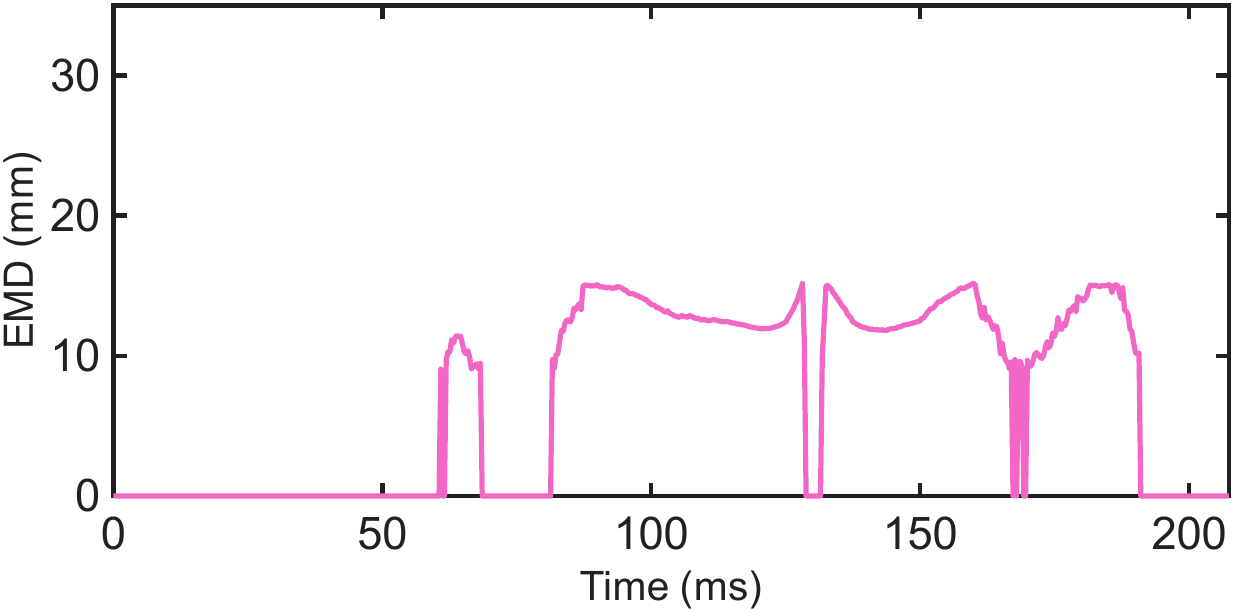}
        \end{minipage}
        
    \end{minipage}\begin{minipage}{\imageWidth}
        \centering
        \begin{minipage}{\linewidth}
            \includegraphics[width=\linewidth]{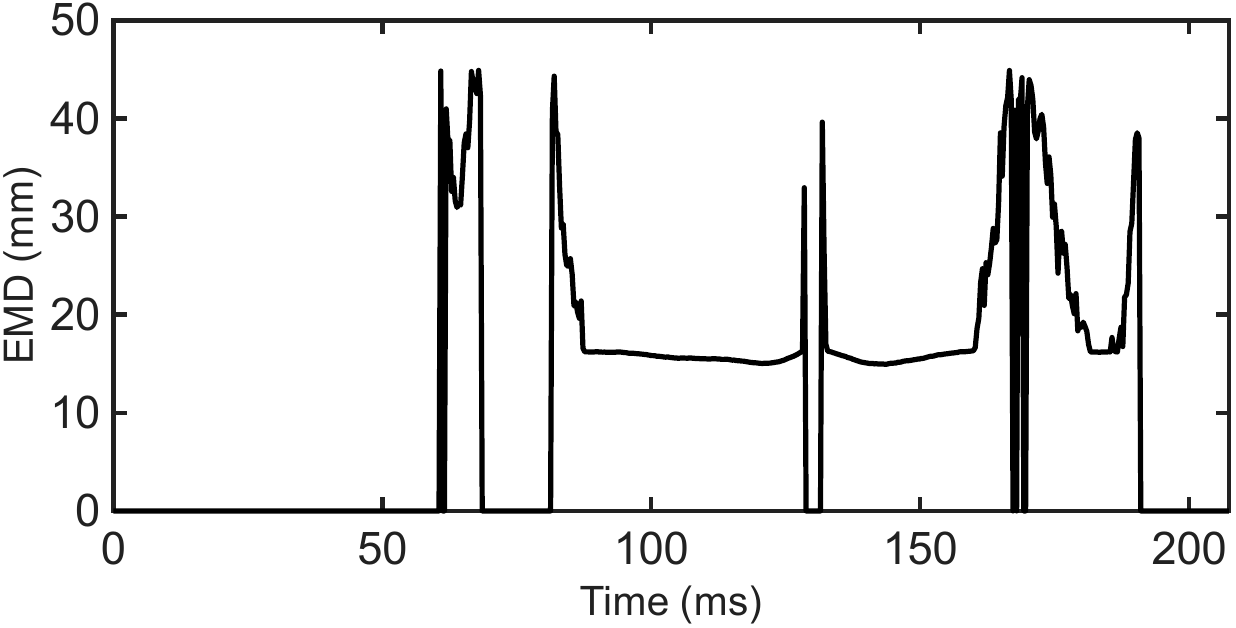}
        \end{minipage}
        
    \end{minipage}
    
    \begin{minipage}{\inverseMethodNameWidth}
        \rotatebox{90}{\bf SKF}
    \end{minipage}\begin{minipage}{\imageWidth}
        \centering
        \begin{minipage}{\sideWidth}
            \rotatebox{90}{left}
        \end{minipage}\begin{minipage}{\linewidth}
            \includegraphics[width=\linewidth]{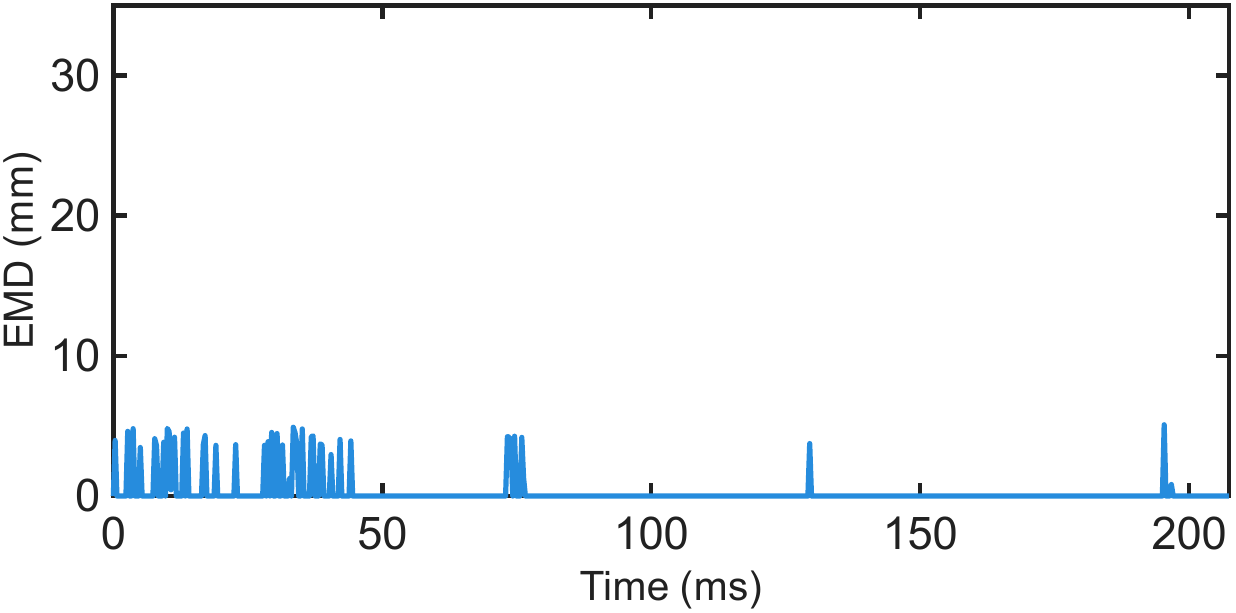}
        \end{minipage}

        \begin{minipage}{\sideWidth}
            \rotatebox{90}{right}
        \end{minipage}\begin{minipage}{\linewidth}
            \includegraphics[width=\linewidth]{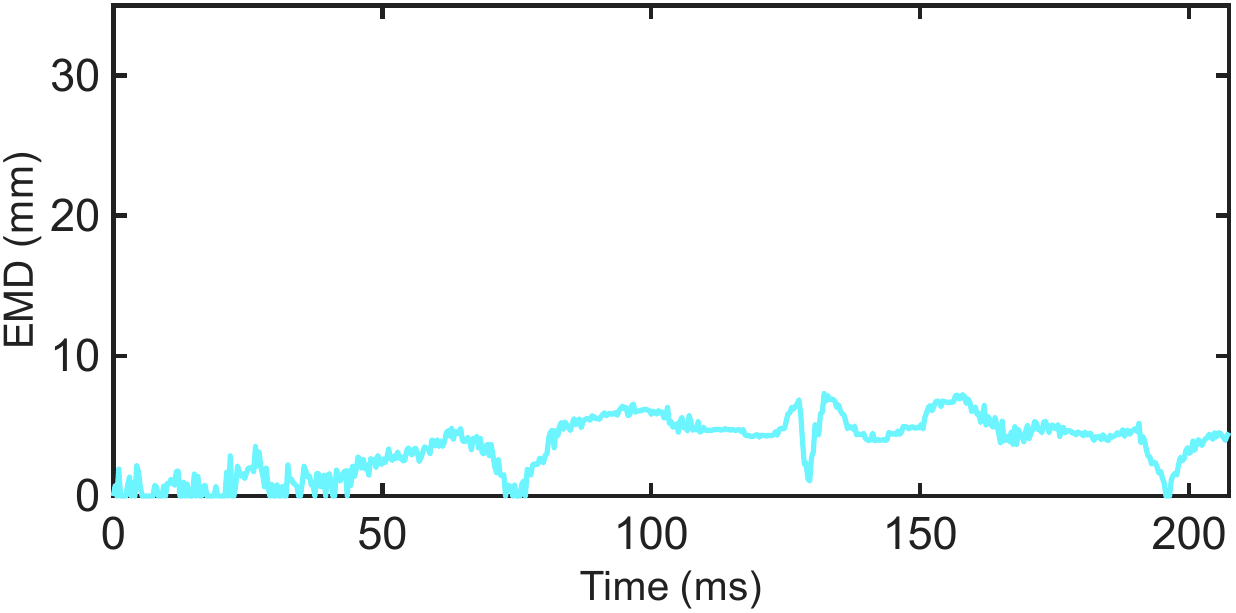}
        \end{minipage}
        
    \end{minipage}\hspace{0.5cm}\begin{minipage}{\imageWidth}
        \centering
        \begin{minipage}{\linewidth}
            \includegraphics[width=\linewidth]{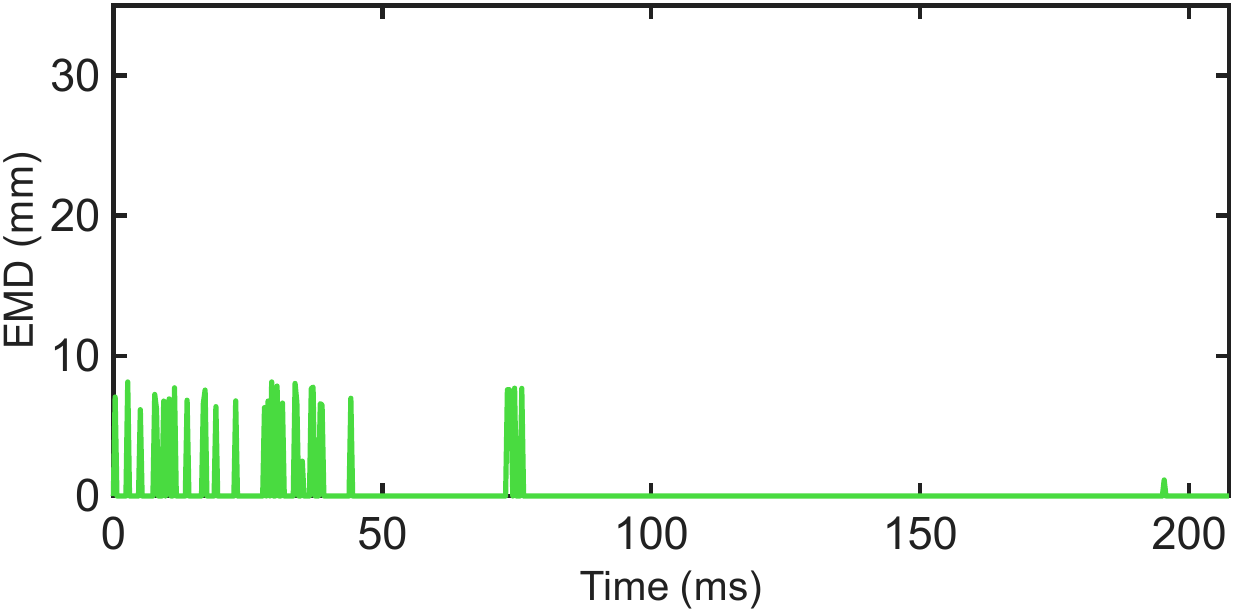}
        \end{minipage}

        \begin{minipage}{\linewidth}
            \includegraphics[width=\linewidth]{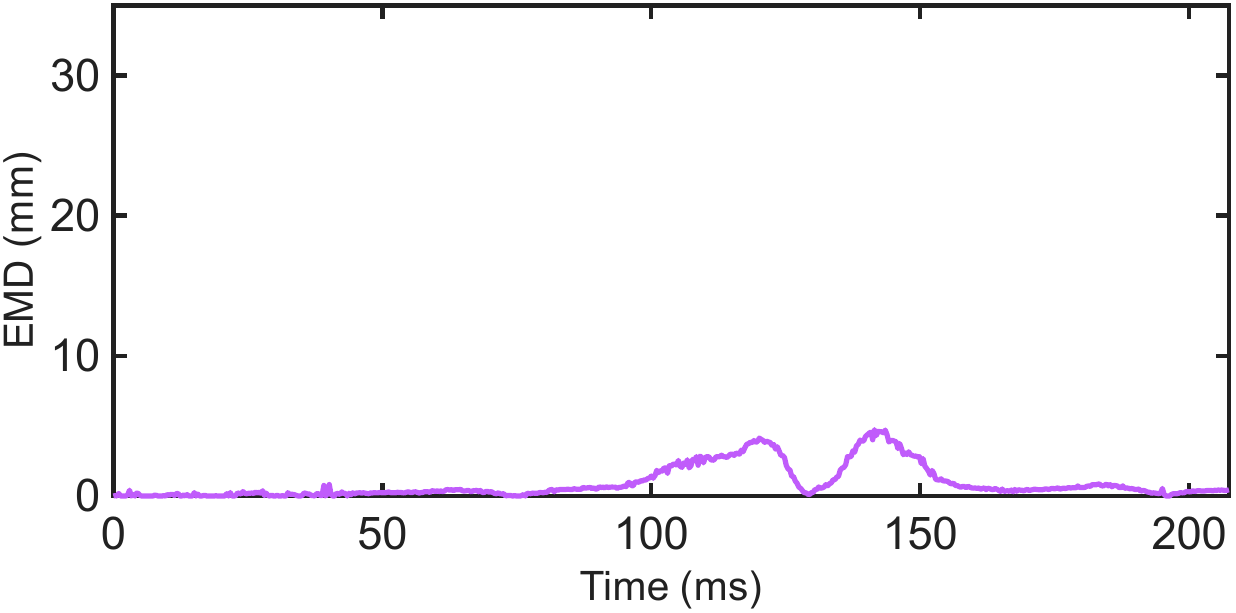}
        \end{minipage}
        
    \end{minipage}\begin{minipage}{\imageWidth}
        \centering
        \begin{minipage}{\linewidth}
            \includegraphics[width=\linewidth]{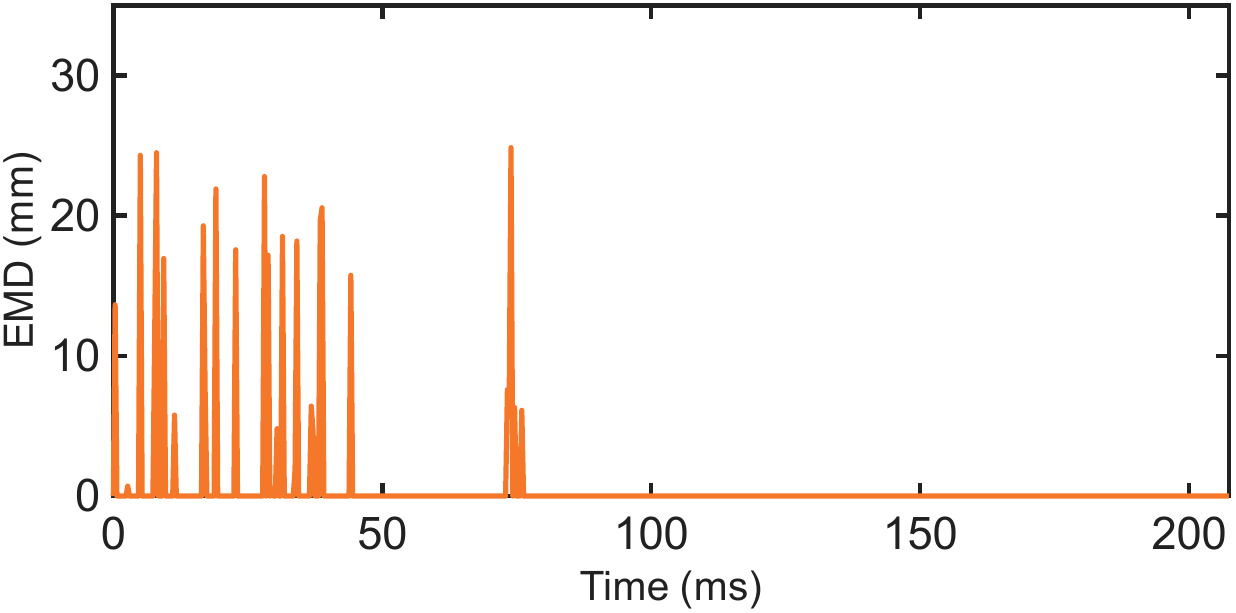}
        \end{minipage}

        \begin{minipage}{\linewidth}
            \includegraphics[width=\linewidth]{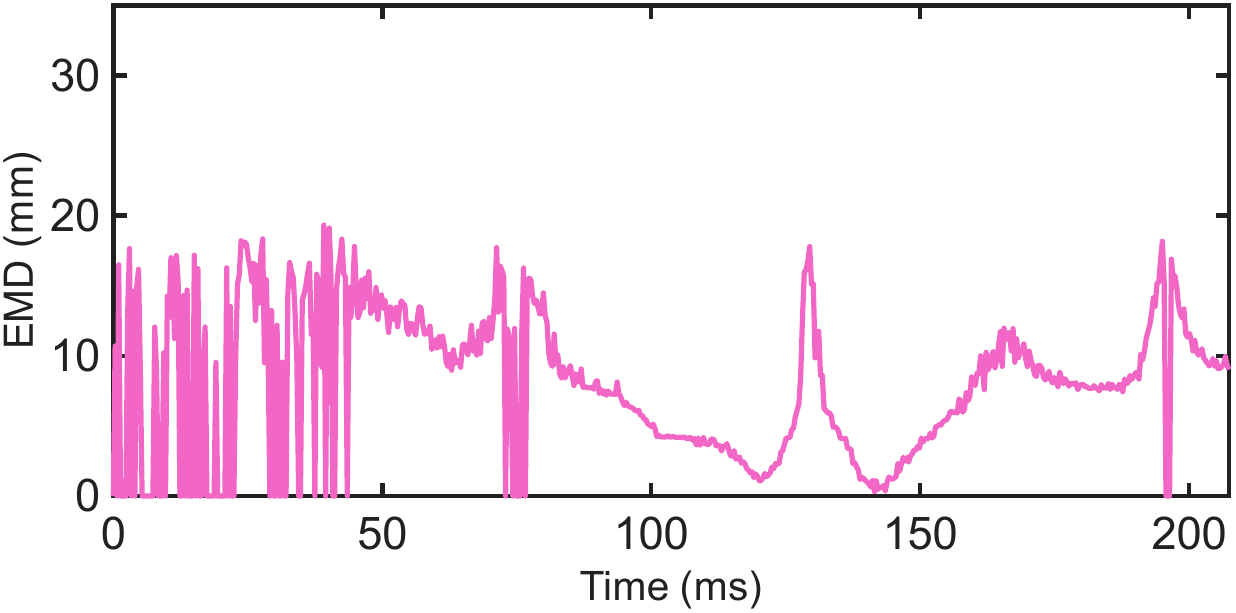}
        \end{minipage}
        
    \end{minipage}\begin{minipage}{\imageWidth}
        \centering
        \begin{minipage}{\linewidth}
            \includegraphics[width=\linewidth]{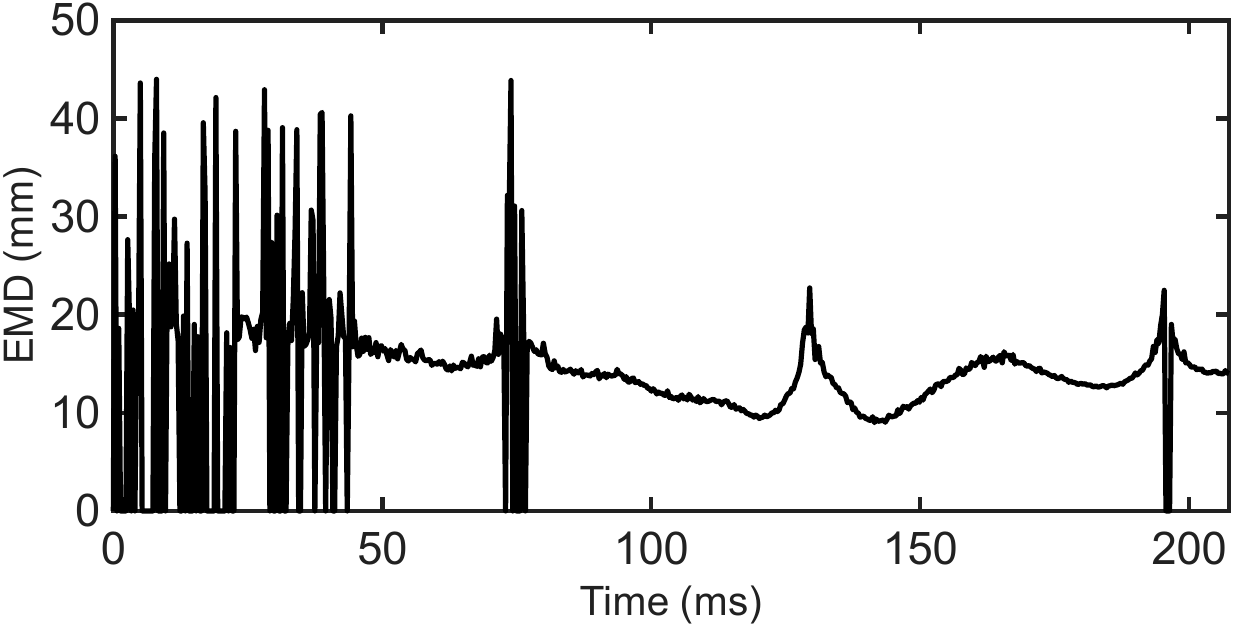}
        \end{minipage}
        
    \end{minipage}

    \begin{minipage}{\inverseMethodNameWidth}
        \rotatebox{90}{\bf DTI-KF}
    \end{minipage}\begin{minipage}{\imageWidth}
        \centering
        \begin{minipage}{\sideWidth}
            \rotatebox{90}{left}
        \end{minipage}\begin{minipage}{\linewidth}
            \includegraphics[width=\linewidth]{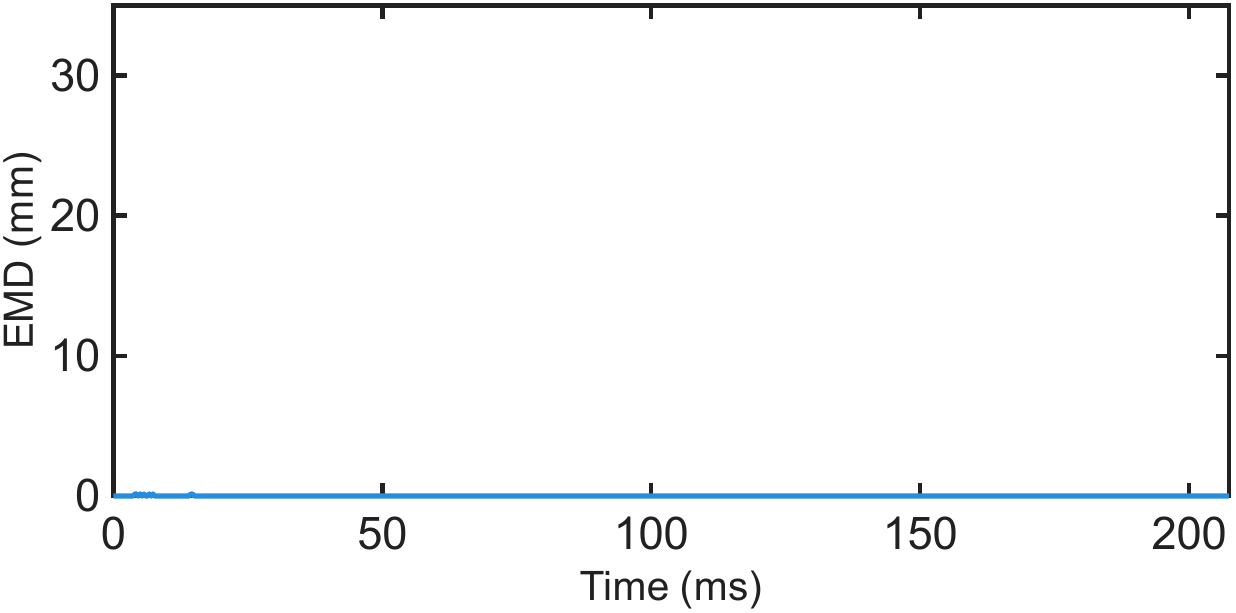}
        \end{minipage}

        \begin{minipage}{\sideWidth}
            \rotatebox{90}{right}
        \end{minipage}\begin{minipage}{\linewidth}
            \includegraphics[width=\linewidth]{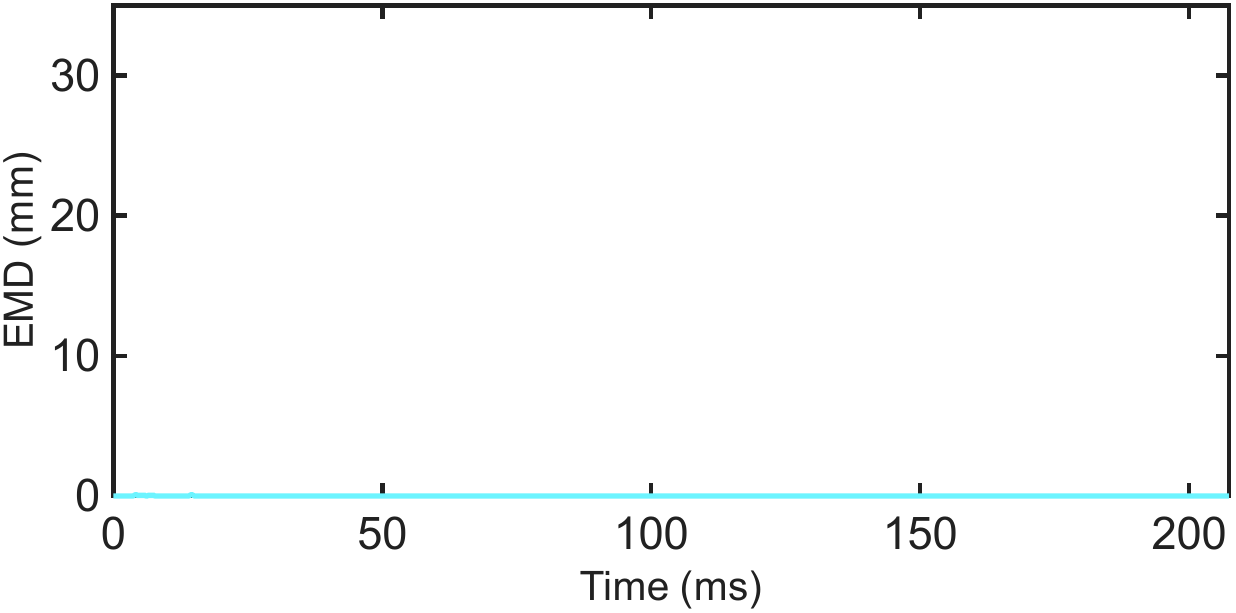}
        \end{minipage}
        
    \end{minipage}\hspace{0.5cm}\begin{minipage}{\imageWidth}
        \centering
        \begin{minipage}{\linewidth}
            \includegraphics[width=\linewidth]{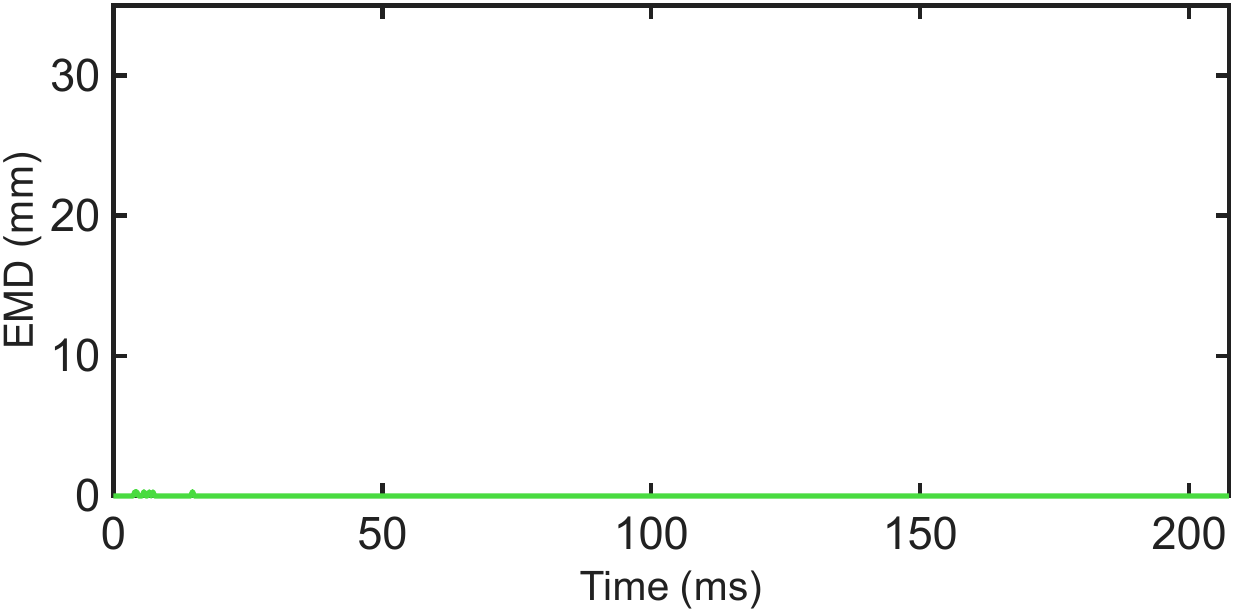}
        \end{minipage}

        \begin{minipage}{\linewidth}
            \includegraphics[width=\linewidth]{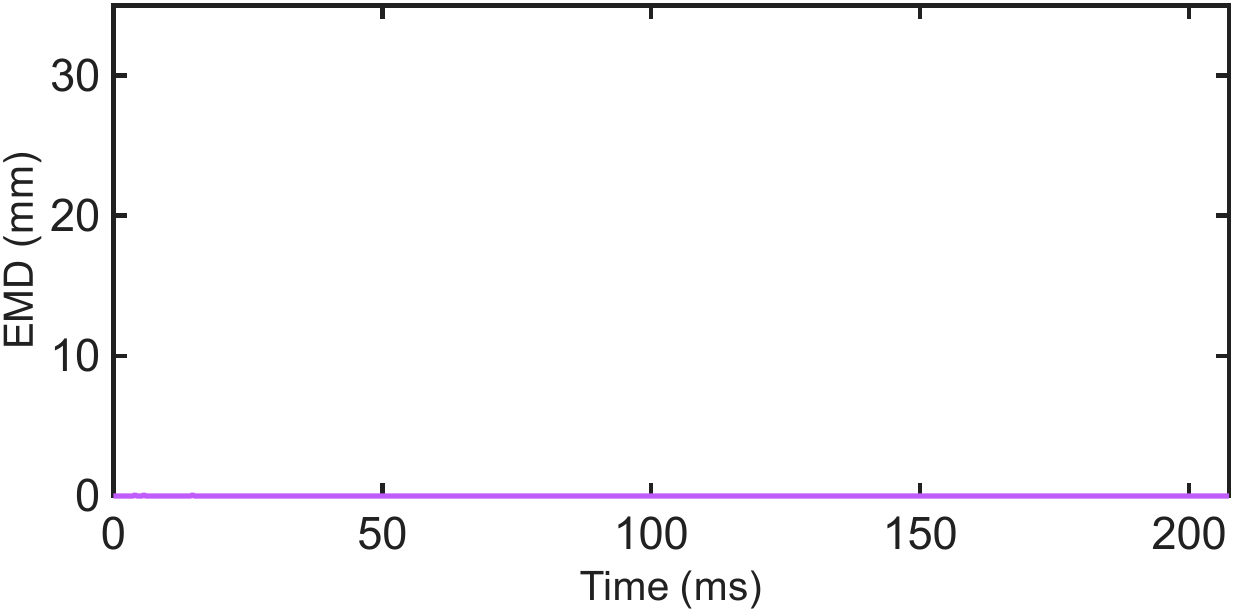}
        \end{minipage}
        
    \end{minipage}\begin{minipage}{\imageWidth}
        \centering
        \begin{minipage}{\linewidth}
            \includegraphics[width=\linewidth]{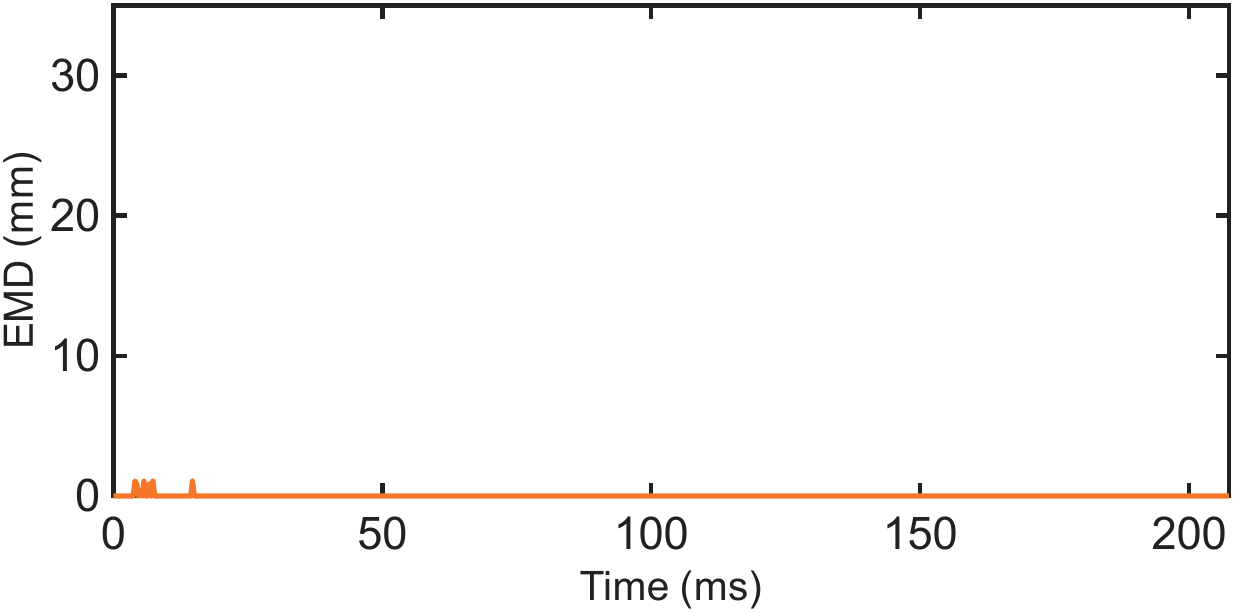}
        \end{minipage}

        \begin{minipage}{\linewidth}
            \includegraphics[width=\linewidth]{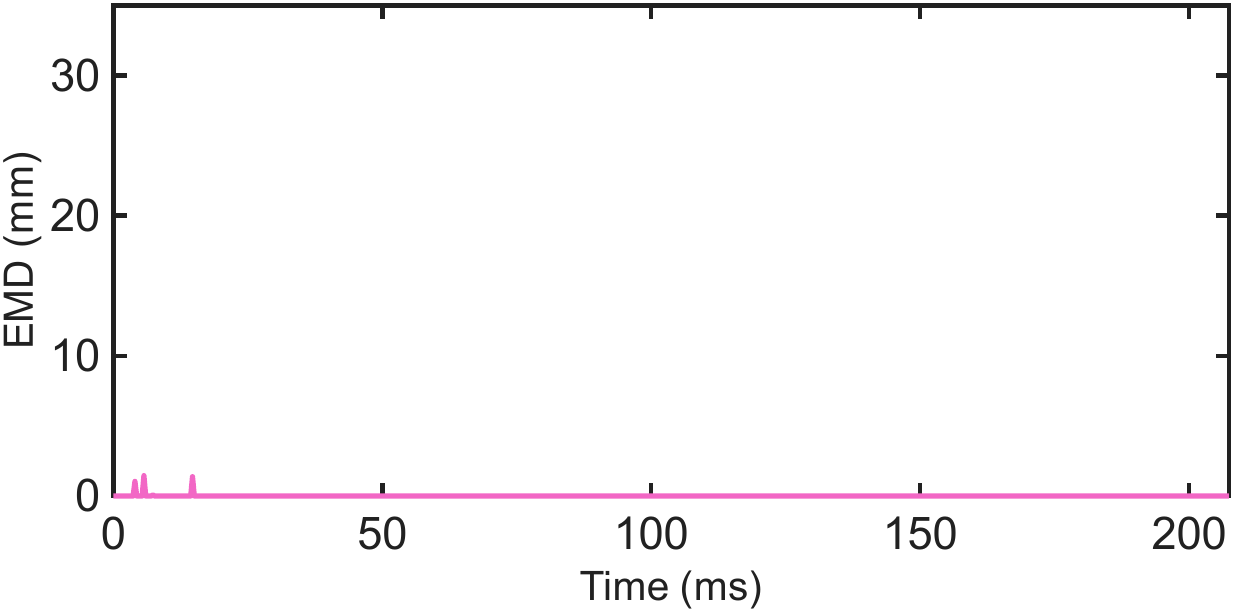}
        \end{minipage}
        
    \end{minipage}\begin{minipage}{\imageWidth}
        \centering
        \begin{minipage}{\linewidth}
            \includegraphics[width=\linewidth]{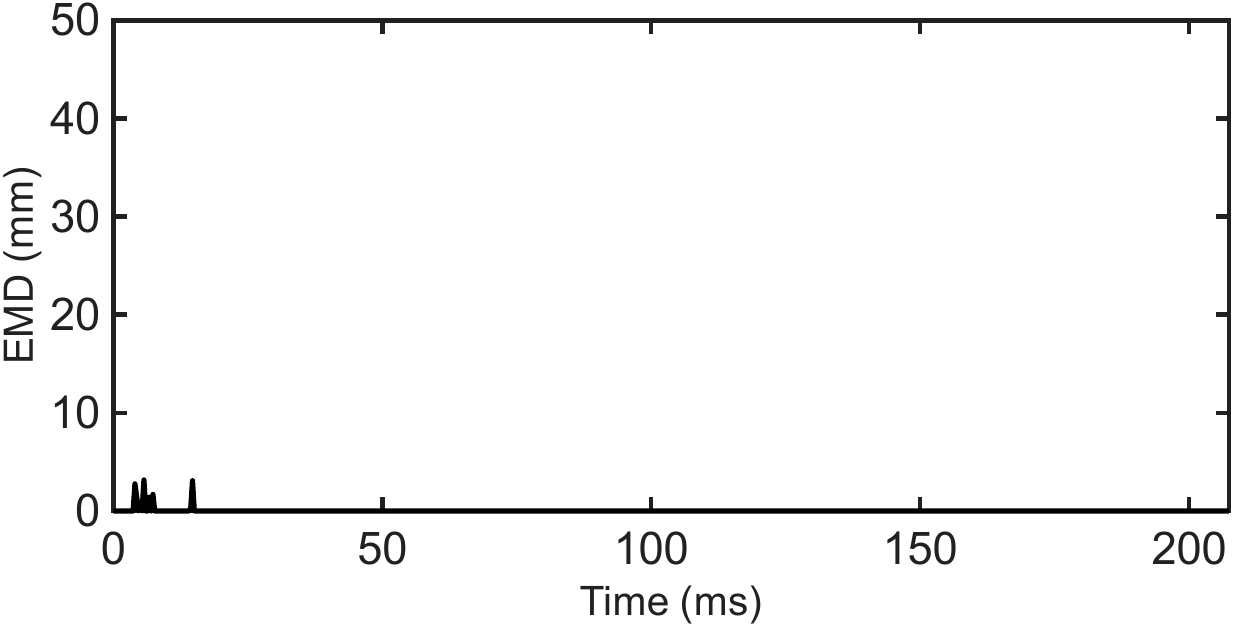}
        \end{minipage}
        
    \end{minipage}
    
    \begin{minipage}{\inverseMethodNameWidth}
        \rotatebox{90}{\bf DTI-SKF}
    \end{minipage}\begin{minipage}{\imageWidth}
        \centering
        \begin{minipage}{\sideWidth}
            \rotatebox{90}{left}
        \end{minipage}\begin{minipage}{\linewidth}
            \includegraphics[width=\linewidth]{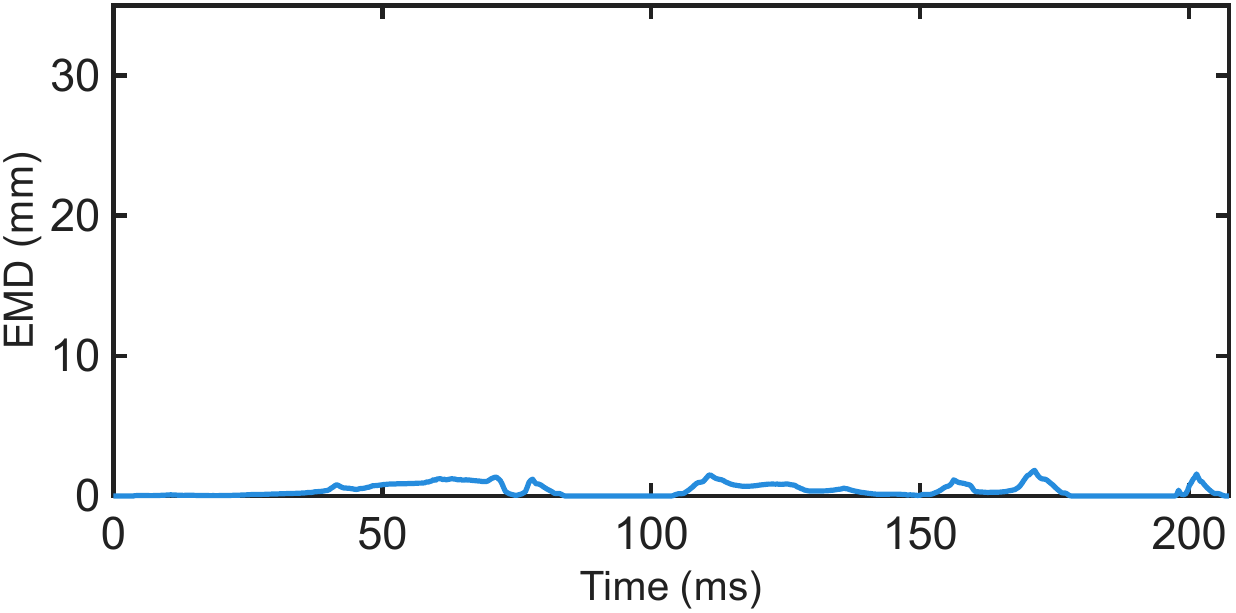}
        \end{minipage}

        \begin{minipage}{\sideWidth}
            \rotatebox{90}{right}
        \end{minipage}\begin{minipage}{\linewidth}
            \includegraphics[width=\linewidth]{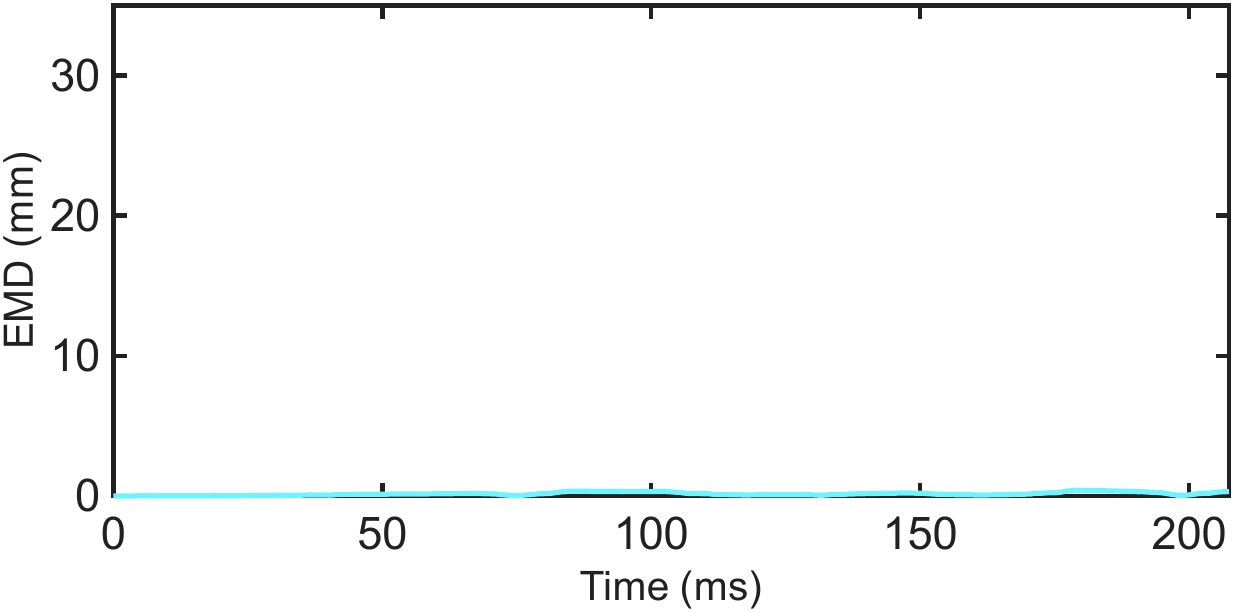}
        \end{minipage}
        
    \end{minipage}\hspace{0.5cm}\begin{minipage}{\imageWidth}
        \centering
        \begin{minipage}{\linewidth}
            \includegraphics[width=\linewidth]{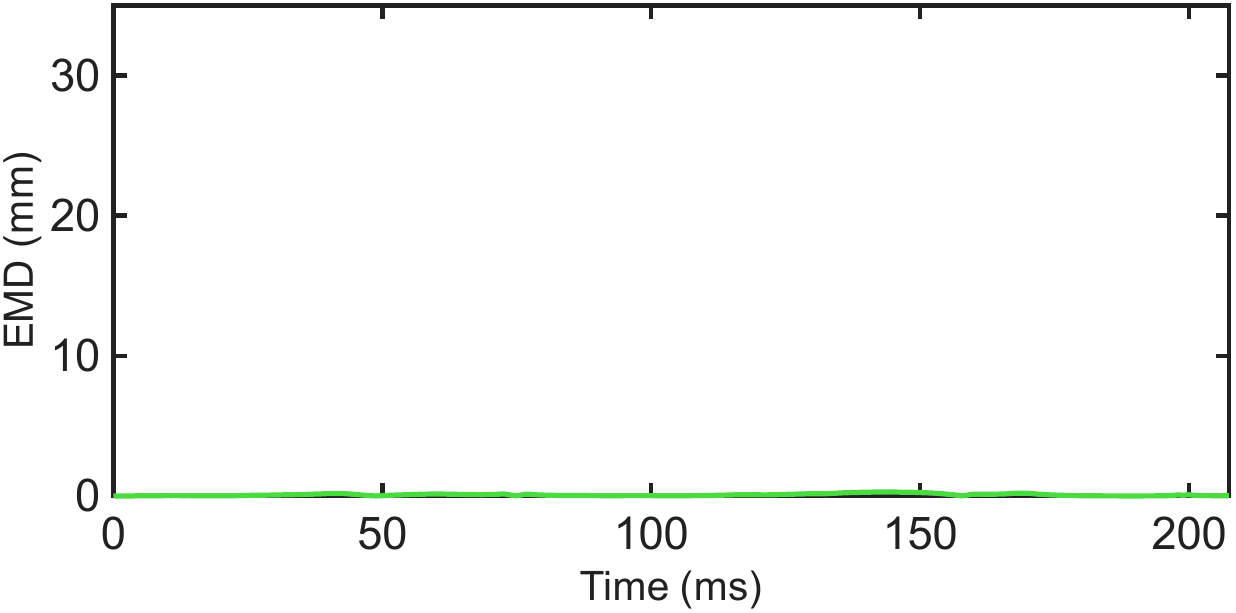}
        \end{minipage}

        \begin{minipage}{\linewidth}
            \includegraphics[width=\linewidth]{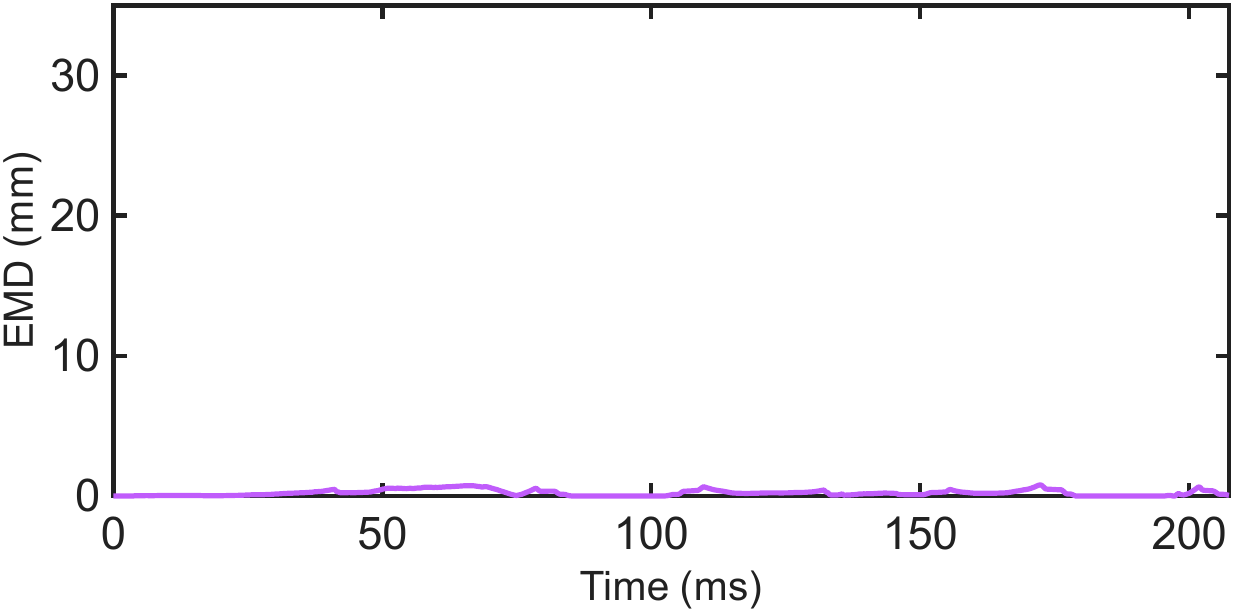}
        \end{minipage}
        
    \end{minipage}\begin{minipage}{\imageWidth}
        \centering
        \begin{minipage}{\linewidth}
            \includegraphics[width=\linewidth]{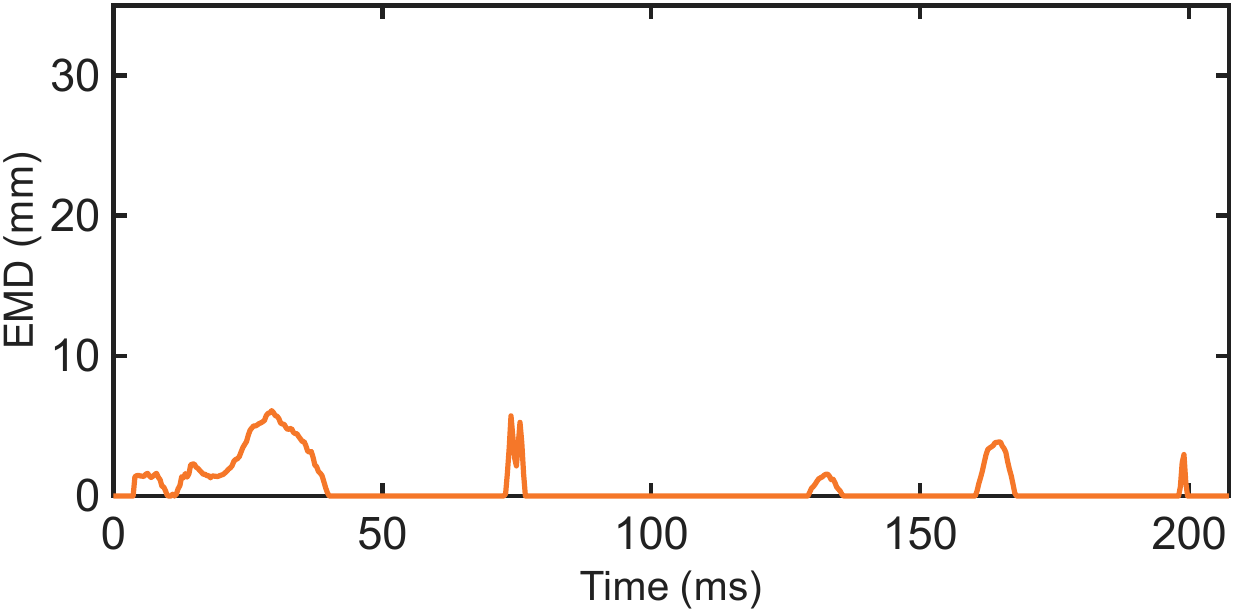}
        \end{minipage}

        \begin{minipage}{\linewidth}
            \includegraphics[width=\linewidth]{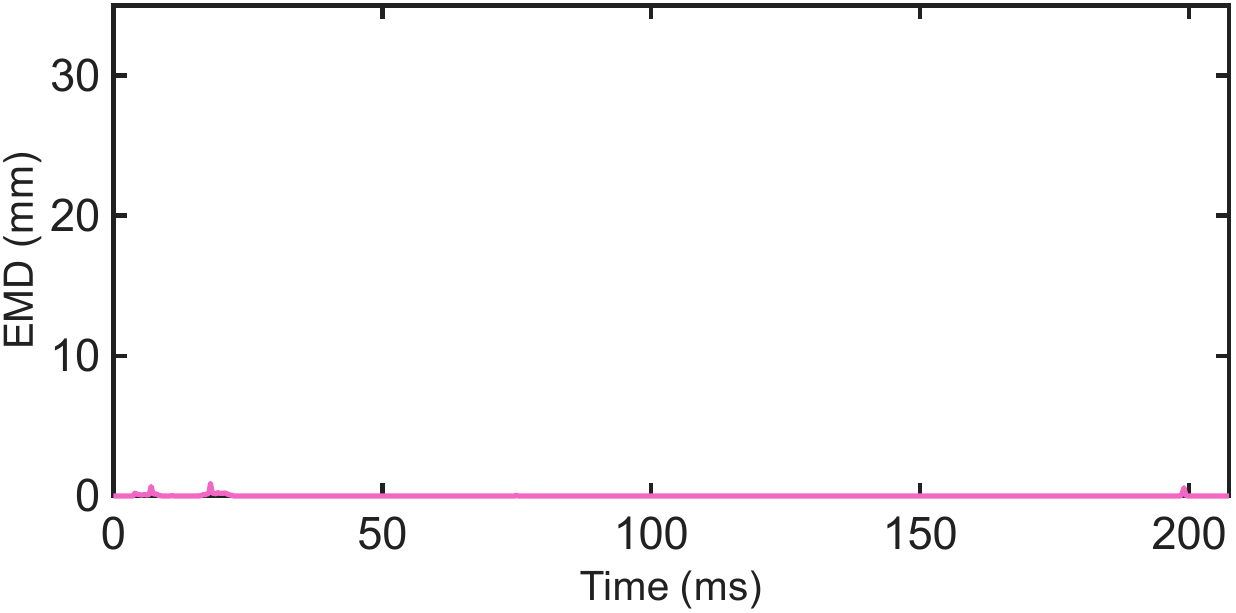}
        \end{minipage}
        
    \end{minipage}\begin{minipage}{\imageWidth}
        \centering
        \begin{minipage}{\linewidth}
            \includegraphics[width=\linewidth]{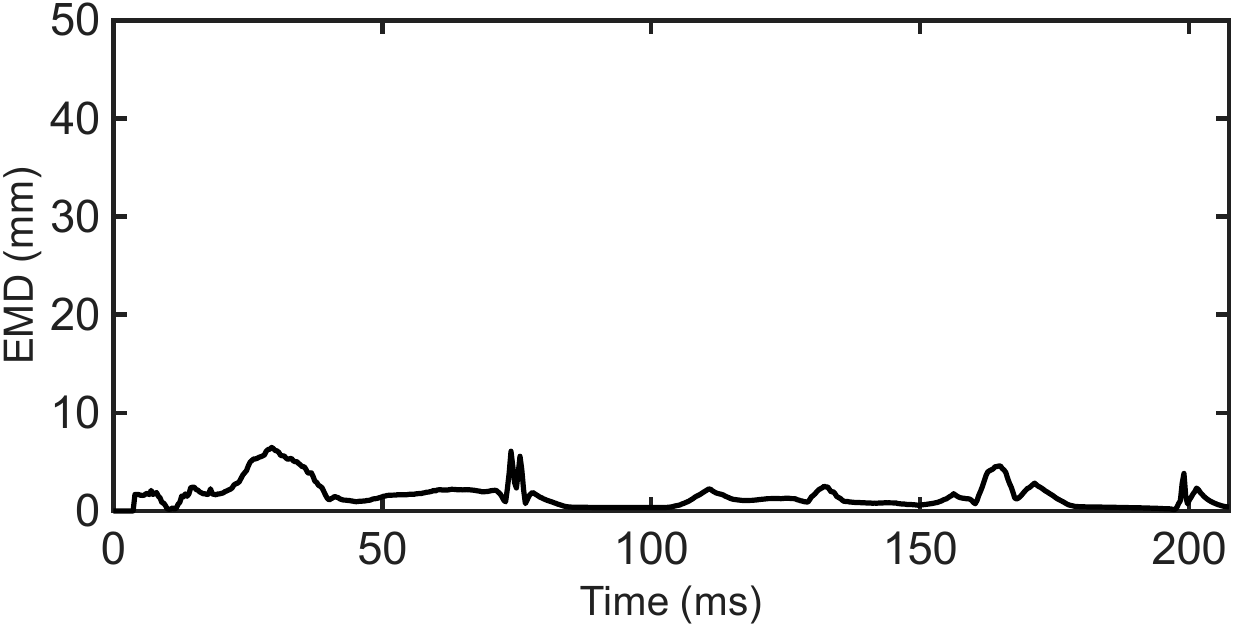}
        \end{minipage}
        
    \end{minipage}

    \caption{Earth Mover's Distances towards individual anatomical regions (insula, superior temporal, transverse temporal) for both hemispheres (left and right sub-rows) and the summed total EMD for random walk evolution-modeled Kalman filter (KF) and Standardized Kalman filter (SKF), and their counterparts using DTI-based transition model (DTI-KF and DTI-SKF).}
    \label{fig:AUDITORYEMD}
\end{figure*}

From the time courses of the estimation values in Figure \ref{fig:AUDITORYtimeseries}, we can see that the estimation of the hemisphere-wise estimation follows the summation trajectory of the true time course of the corresponding side; hence, the left and right activity can be separated from each other. However, the regional estimates from the random-walk KF and SKF correlate nearly perfectly, and the maximum estimate appears in the superior temporal gyrus (blue) or middle temporal gyrus (green). Therefore, we cannot localize the estimate beyond the correct hemisphere. With the DTI model, ideal localization is not guaranteed beyond the insula, the largest of the inspected regions by volume. The estimation in the right transverse temporal gyrus, while being overshadowed by values in the insula, follows the true trajectory relatively closely after 150 \unit{\milli\second}. By examining the time evolution in the insula more closely, we see that the model allows separation between regional estimates. EMDs (Figure \ref{fig:AUDITORYEMD}) show that the majority of mislocalized estimates are pushed to the transverse temporal gyri for KF and SKF. The DTI model causes a significant drop in the EMD values. While the activity does not peak at the right locations, it appears close to the active regions.

\clearpage

\section{Discussion}
In this paper, we demonstrated a strategy to apply raw, tractography-free DTI data to the known Kalman filter approaches used in EEG brain source imaging. DTI data were used to create a probability-based connectivity map that models brain fluid flow from one anatomical brain region to another. We show that this flow can be interpreted as a Riemann manifold, allowing geodesic algorithms to determine anatomical connectivity. In the algorithmic development, the emphasis was on efficiency, as the rigorous voxel-based computation scheme is computationally heavy. Namely, the obtained connectivity map is then applied to the transition matrix that describes the assumed dynamical evolution of the brain activity. The DTI-based model describes the transition as an instantaneous decay of activity from the starting brain region to equal activity gain in target regions. The Kalman filters use the observation time series to infer the rate and regions where activity will most likely be targeted. Based on the results, this model prevents the estimate from spreading far from the true source location, even in complex auditory evoked potentials. It allows a slight distinction of nearby active regions, which is not possible with the random walk evolution model, as shown by the results and earlier research \cite{Lahtinen2024SKF,Dilshanie2026SKF,piispa2025DSKF}. The main advantage is the ability to detect deep activity with a non-standardized Kalman filter. Based on the results, we conclude that the connectivity-based translation model is more beneficial for the Kalman filter than the standardized Kalman filter. This is due to the "low-resolution topography" aspect it inherits from the standardization \cite{PascualMarqui2002}. This makes the estimation more spatially widespread, making qualitative source analysis more open to interpretation and increasing the correlation of estimation values in neighboring anatomical regions. 

The introduced methodology has some limitations worth mentioning: while all parts of the algorithms are designed to be as light as possible, solving a DTI-conductivity model is still time-consuming. By pre-processing the DTI data, the diffusion field can be made much smoother and less pathological, reducing the need to consider special cases, such as violated non-degeneracy. A smoother diffusion field lets one generalize the geodesic path from one region to another with few voxel-to-voxel geodesic computations. This could also make the Monte Carlo approach plausible for geodesic computations. We see these as future research topics that make DTI-based Kalman filtering applicable in clinical settings.

\section*{Acknowledgement}
This study was supported by Tampere University and the Flagship of Advanced Mathematics for Sensing, Imaging and Modelling (FAME) (359185). The author would like to thank 
Dr.\ Narayan Puthanmadam Subramaniyam for granting access to his coupled Jansen--Rit neural mass model code and Prof.\ Sampsa Pursiainen for computational resources. 

\section*{Declaration of Competing Interests}
The author have no known competing financial interests or personal relationships that could have appeared to influence the work reported in this manuscript. 

\onecolumn
\appendix
\section{Fokker-Planck equation on Riemann manifolds and diffusion tensor--metric connection}
\label{app:FP-DTg}

Let $\lbrace {\bf x}_t\colon t\geq 0\rbrace$ be a time-homogeneous stochastic process. Consider the following derivative
\begin{equation}
    \int_\Omega h({\bf x})\frac{\partial p({\bf x},t\mid {\bf y})}{\partial t}\mathrm{d}{\bf x}=\frac{\partial }{\partial t}\int_\Omega h({\bf x})\partial p({\bf x},t\mid {\bf y})\mathrm{d}{\bf x},
\end{equation}
where $h$ is a smooth function with compact support. By writing the Fréchet derivative $D_t$ for the function, 
\begin{equation}
    D_t(s)=\int_\Omega h({\bf x})\left(p({\bf x},t+s\mid {\bf y})-p({\bf x},t\mid {\bf y})\right)\mathrm{d}{\bf x}
\end{equation}
and by using the Chapman-Kolmogorov identity, we obtain
\begin{equation}
    D_t(s)=\int_\Omega h({\bf x})\left(\int_\Omega p({\bf x},s\mid {\bf z})p({\bf z},t\mid {\bf y})\mathrm{d}{\bf z}-p({\bf x},t\mid {\bf y})\right)\mathrm{d}{\bf x}
\end{equation}
By changing the variable from the subtracted term from ${\bf y}$ to ${\bf z}$, we get
\begin{equation}
    D_t(s)=\int_\Omega p({\bf z},t\mid {\bf y})\int_\Omega p({\bf x},s\mid {\bf z})\left(h({\bf x})-h({\bf z})\right)\mathrm{d}{\bf x}\, \mathrm{d}{\bf z}.
\end{equation}
Using Taylor expansion up to order 2 for $h({\bf x})$ around ${\bf z}$, we obtain
\begin{equation}
    h({\bf x})-h({\bf z})=\frac{1}{2}\nabla h^T\nabla_g d({\bf x},{\bf z})^2+\frac{1}{4}\left(\nabla_g d({\bf x},{\bf z})^2\right)^TH_{g}^\sharp h\nabla_g d({\bf x},{\bf z})^2,
\end{equation}
where $\nabla_g$ is gradient on the manifold $(\mathcal{M},g)$ and $H_g^\sharp=g^{ik}(H_{ij}-\Gamma_{ij}^k\partial_k)$ is the Hessian operator. Integrating the first-order term by parts yields
\begin{equation}
    \frac{1}{2}\int_\Omega \psi\nabla h^T\nabla_g d({\bf x},{\bf z})^2\mathrm{d}{\bf z}=-\frac{1}{2}\int_\Omega h({\bf z})\nabla\cdot \left(\psi\nabla_g d({\bf x},{\bf z})^2\right)\mathrm{d}{\bf z})
\end{equation}
and the second order
\begin{equation}
\begin{split}
    \frac{1}{4}\int_\Omega \psi&\, \left(\nabla_g d({\bf x},{\bf z})^2\right)^TH_g^\sharp[h]\nabla_g d({\bf x},{\bf z})^2\mathrm{d}{\bf z}\\
    &=-\frac{1}{4}\int_\Omega (\nabla_g\psi)^T \nabla_g d({\bf x},{\bf z})^2\left(\nabla_g d({\bf x},{\bf z})^2\right)^T\nabla_g h \mathrm{d}{\bf z}\\
    &=\frac{1}{4}\int_\Omega h\left[\left(\left(\nabla_g d({\bf x},{\bf z})^2\right)^T\nabla\right)^2\psi \right]\mathrm{d}{\bf z}
\end{split}
\end{equation}
where $\psi=p({\bf z},t\mid x)p({\bf x},s\mid {\bf z})$. Now we make connections
\begin{align}
    \mu &= \frac{g_{ij}}{2}p({\bf x},\Delta t\mid {\bf z})\nabla_g d({\bf x},{\bf z})\\
    \sigma_{ki}\sigma_{kj}&=\frac{g_{ij}}{2}\nabla_g d({\bf x},{\bf z})^2\left(\nabla_g d({\bf x},{\bf z})^2\right)^T p({\bf x},\Delta t\mid {\bf z}).
\end{align}
This yields the following representation of the Fokker-Planck equation:
\begin{equation}
    \partial_t p=\nabla(g^{ij}\mu p)+\nabla^T(\sigma g^{ij}\sigma^T p)\nabla_L,
\end{equation}
where $\nabla_L$ is the Euclidean gradient acting on the left.

Now, assuming driftless movement, i.e., $\mu = 0$, by setting $g^{ji}=(\sigma_{ki}\sigma_{kj})^{\dagger}_{ij}$, where $(\cdot)^\dagger$ is the pseudo-inverse, the distribution on $(\mathcal{M},g)$ follows the distribution of a standard Wiener process on a flat space with holes at the points where $\sigma({\bf x})=O$. Therefore, $\sigma^T\sigma$ describes contravariant metric of $(\mathcal{M},g)$.

\section{Convergence of the shooting algorithm for geodesics}\label{app:shootingConverg}
As the shooting method provides a linear approximation for a path $\bm{\gamma}$, due to the Lagrange error bound for the Taylor series:
\begin{equation}
    \left|\int_{t}^{t+\Delta t}\left(\left\|\mathbf{x}_m-\mathbf{x}_p\right\|_g-\left\|\dot{\bm{\gamma}}\right\|_g\right)\mathrm{d}t\right|\leq \frac{\Delta t^2}{4}\max\left\|\ddot{\bm{\gamma}}\right\|_g
\end{equation}
Since for geodesics
\begin{equation}
   \Delta d_\gamma:= d_\gamma(\mathrm{x},\mathrm{x}+\Delta \mathrm{x})=\Delta t,
\end{equation}
we can write
\begin{equation}
    \left|\int_{t}^{t+\Delta t}\left(\left\|\mathbf{x}_m-\mathbf{x}_p\right\|_g-\left\|\dot{\bm{\gamma}}\right\|_g\right)\mathrm{d}t\right|\leq \frac{\Delta d_\gamma^2}{4}\max\left\|\ddot{\bm{\gamma}}\right\|_g.
\end{equation}
As we cut the distance of every sub-interval in half after $(n^2+n)/2$ steps:
\begin{equation}
    \left|\int_{t}^{t+\Delta t}\left(\left\|\mathbf{x}_m-\mathbf{x}_p\right\|_g-\left\|\dot{\bm{\gamma}}\right\|_g\right)\mathrm{d}t\right|\leq \frac{d^2}{4}\left\|\ddot{\bm{\gamma}}^*\right\|_g\left(\frac{1}{2}\right)^{2K},
\end{equation}
where $\left\|\ddot{\bm{\gamma}}^*\right\|_g=\max\left\|\ddot{\bm{\gamma}}\right\|_g$ and $K=(\sqrt{8n+1}-1)/2$, so in terms of $n$ we have
\begin{equation}
    \left|\int_{t}^{t+\Delta t}\left(\left\|\mathbf{x}_m-\mathbf{x}_p\right\|_g-\left\|\dot{\bm{\gamma}}\right\|_g\right)\mathrm{d}t\right|\leq \frac{d^2}{2}\left\|\ddot{\bm{\gamma}}^*\right\|_g\left(\frac{1}{2}\right)^{\sqrt{8n+1}}\leq \frac{d^2}{2\cdot 4^{\sqrt{2}}}\left\|\ddot{\bm{\gamma}}^*\right\|_g\left(\frac{1}{2}\right)^{\sqrt{n}}.
\end{equation}

\section{Derivation of the hitting probability formula}\label{app:probFormula}
First, we make the observation that for any unit vector $\mathbf{v}$:
\begin{equation}
    \mathbf{v}^T\mathbf{W}_t=\sum_{k=1}^d v_kW_t^{(k)}=\sum_{k=1}^d W_{t/v_k^2}^{(k)},
\end{equation}
and hence
\begin{equation}\label{eq:wiener}
    \mathbb{E}\left[\mathbf{v}^T\mathbf{W}_t\right]=0,\quad \mathrm{var}\left[\mathbf{v}^T\mathbf{W}_t\right]=\left\|\mathbf{v}\right\|_2^2t=t.
\end{equation}
So the resulting random process is a Wiener process.

Let us then have a path $\mathcal{C}(t):=(r(t),\theta,z(t))=\gamma(t)$ that induces a local, rotationally symmetric, cylindrical coordinate system. The integral can be written as:
\begin{equation}
    \int_{\gamma(0)}^{\gamma(t)}\int_{0}^{r(\gamma)} \int_{\theta_1}^{\theta_2}\mathrm{d}\theta\, p(r,\gamma,\cdot,t)\,\mathrm{d}r\mathrm{d}z\propto \int_{\gamma(0)}^{\gamma(t)}\int_{0}^{r(\gamma)} p(r,\gamma,t)\,\mathrm{d}r\mathrm{d}\gamma.
\end{equation}
Using integration by substitution
\begin{align}
    \int_{\gamma(0)}^{\gamma(T)}\int_{0}^{r(\gamma)} p(r,\gamma,t)\,\mathrm{d}r\mathrm{d}\gamma&=\int_{\gamma(0)}^{\gamma(T)}\int_{0}^{\gamma} p(r(\gamma),\gamma,t)\frac{\mathrm{d}r}{\mathrm{d}\gamma}\,\mathrm{d}^2\gamma\\
    &=\int_{\gamma(0)}^{\gamma(T)}\int_{0}^{t} p(r(\tau),\gamma,t)\dot{r}(\tau)\,\mathrm{d}\tau\mathrm{d}\gamma\\
    &=\int_{0}^{T}\int_{0}^{t} p(r(\tau),\gamma(t),t)\dot{r}(\tau)\dot{\gamma}(t)\,\mathrm{d}\tau\mathrm{d}t.
\end{align}
For the derivatives, we have
\begin{equation}
    \dot{r}=\frac{\mathrm{d}}{\mathrm{d} t}\left(\frac{d_r^2}{2t}\right)=d_r\frac{d_r'}{t}-\frac{d_r^2}{2t^2}=\frac{d_r^2}{2t}\left(2\frac{d_r'}{d_r}-\frac{1}{t}\right)=\frac{d_r^2}{2t}\frac{\mathrm{d}}{\mathrm{d} t}\left(\log(d_r^2)-\log(t)\right)=\frac{d_r^2}{2t}\frac{\mathrm{d}}{\mathrm{d} t}\log\left(\frac{d_r^2}{2t}\right)
\end{equation}
Therefore, the integral can be rewritten as
\begin{equation}
    \int_{0}^{T}\int_{0}^{\log q(t)} \hat{q}\,p(\hat{q},\gamma(t),t)\dot{\gamma}(t)\,\mathrm{d}\log(\hat{q})\,\mathrm{d}t=\int_{0}^{T}\int_{0}^{q(t)} p(\hat{q},\gamma(t),t)\dot{\gamma}(t)\,\mathrm{d}\hat{q}\,\mathrm{d}t,
\end{equation}
and assuming $\gamma$ is geodesic, then by definition, $\dot{\gamma}\equiv 1$ and thus
\begin{equation}\label{eq:integralId}
    \iiint p(r(\gamma),\gamma,\theta,t)\mathrm{d}V\propto \int_{0}^{T}\int_{0}^{q(t)} p(\hat{q},d_\gamma(t),\cdot,t)\,\mathrm{d}\hat{q}\,\mathrm{d}t=\int_{0}^{d_\gamma(T)}\int_{0}^{q(z)} p(\hat{q},z,\cdot ,t)\,\mathrm{d}\hat{q}\,\mathrm{d}z.
\end{equation}

Let $\gamma(T)$ be the geodesic of interest. We can define the probability distribution in a cylindrical, tube-like coordinate system as
\begin{equation}
    p(\mathbf{x},T)\propto \exp\left(\frac{d_\perp (\mathbf{x}_r)^2+d_\gamma^2}{2T}\right).
\end{equation}
We can then marginalize the variables orthogonal to the path at each point, so that we compute $P(\mathbf{x}_r(t))\in S_\mathrm{tube}(t))$. 

In local cylindrical coordinates:
\begin{equation}
    p(r,\theta,z,t)=\frac{r}{2\pi \sigma_r^2t\sqrt{2\pi \sigma_z^2 t}}\exp\left(-\frac{r^2}{2\sigma_r^2 t}\right)\exp\left(-\frac{z^2}{2\sigma_z^2 t}\right),
\end{equation}
Then
\begin{equation}
    \int_{\theta_1}^{\theta_2}\int_{0}^{r_0(T)} p(r,\theta,t)\mathrm{d}r\mathrm{d}\theta = \frac{\Delta\theta}{2\pi}\frac{1}{\sqrt{2\pi \sigma_z^2 t}}\left(1-\exp\left(-\frac{r_0^2}{2\sigma_r^2T}\right)\right)\exp\left(-\frac{z^2}{2\sigma_z^2 t}\right)
\end{equation}
and
\begin{equation}
    \int_{z_0}^\infty\int_{\theta_1}^{\theta_2}\int_{0}^{r_0} p(r,\theta,t)\mathrm{d}r\mathrm{d}\theta\mathrm{d}z=\frac{\Delta\theta}{4\pi}\left(1-\exp\left(-\frac{r_0^2}{2\sigma_r^2T}\right)\right)\left(1-\mathrm{erf}\left(\frac{z_0}{\sqrt{2\sigma_z^2 t}}\right)\right).
\end{equation}
Now, as the difference
\begin{equation}
    P_\infty-P_T=\frac{\Delta\theta}{4\pi}\left(1-\exp\left(-\frac{r_0^2}{2\sigma_r^2T}\right)\right)\mathrm{erf}\left(\frac{z_0}{\sqrt{2\sigma_z^2 t}}\right).
\end{equation}

According to the observed change of variables in Eq. (\ref{eq:integralId}), the result can be extended to general continuous paths. Moreover, by using the identity in Eq. (\ref{eq:Td-connection}) and $r_0=h/2$, we get

\begin{equation}
    P_\infty-P_T\propto \left(1-\exp\left(-\frac{h^2}{8T(\mathbf{x},\mathbf{x}_0)}\right)\right) \mathrm{erf}\left(\frac{d_\sigma(\mathbf{x},\mathbf{x}_0)}{\sqrt{2T(\mathbf{x},\mathbf{x}_0)}}\right).
\end{equation}

\twocolumn


 \bibliographystyle{elsarticle-num} 
 \bibliography{biblo}



\end{document}